\documentclass[12pt,a4paper,english]{article}
\usepackage{amsfonts,amssymb,latexsym,xspace,textcomp}
\usepackage{amsmath,enumerate}
\usepackage{amsthm}
\usepackage{graphicx}
\usepackage{color}
\usepackage[english]{babel}
\usepackage[latin1]{inputenc}
\usepackage[linktocpage]{hyperref}
\numberwithin{equation}{subsection}
\usepackage{titlesec}
\titleformat{\subsubsection}[runin]
{\normalfont\large\bfseries}{\thesubsubsection}{1em}{}
\theoremstyle{definition}
\newtheorem{theorem}[subsubsection]{Theorem}
\newtheorem{definition}[subsubsection]{Definition}
\newtheorem{example}[subsubsection]{Example}
\newtheorem{remark}[subsubsection]{Remark}
\newtheorem{proposition}[subsubsection]{Proposition}
\newtheorem{lemma}[subsubsection]{Lemma}
\newtheorem{corollary}[subsubsection]{Corollary}
\begin{document}

\begin{center}
\textbf{{\Large Singularities of Fredholm maps with one-dimensional kernels, I: \\ \vspace{10pt} 
A complete classification}}\footnote{This research was partially supported by MIUR project ``Elliptic and Hamiltonian Differential Problems and their applications''.} \\
\vspace{50pt}
\textbf{F. Balboni} \\ \vspace{8pt} 
\textit{\footnotesize Dipartimento di Matematica, Università degli Studi ``Tor Vergata'', Via della Ricerca Scientifica, 00133, Roma, Italy}\\ 
{\footnotesize E-mail: balboni@mat.uniroma2.it} \\ \vspace{18pt} 
and \\  \vspace{18pt} 
\textbf{F. Donati   } \textbf{\textdied}\\
\vspace{8pt} 
\textit{\footnotesize Dipartimento di Matematica, Università degli Studi ``Tor Vergata'', Via della Ricerca Scientifica, 00133, Roma, Italy}\\ \vspace{18pt} 
\end{center}
\textbf{\footnotesize Abstract}\bigskip  \\
\indent \footnotesize We study the simplest singular points of Fredholm maps of index zero between Banach spaces, i.e. when the kernel of the Fréchet derivative of the map has dimension one. Even in this relatively simple case we have a rich variety of singularities which are completely classified under the natural geometric assumption of \textit{transversality} of the map. In fact we have, locally, a suitable \textit{stratification} of the singular points that allows us to identify three kinds of singularities: a) the \textit{ordinary} or $k$-singularities (the infinite-dimensional analogues of the well-known Morin singularities) and  two new types, b) the \textit{maximal-transverse} singularities; c) the \textit{infinite-transverse} singularities (the latter ones can only occur in infinite dimensions). 

\normalsize
\section*{{\large INTRODUCTION}}
\quad A quite simple problem like the periodic problem for the Riccati equation
\begin{equation*}
(\text{P}_{1}) \begin{cases} u' + a(t)u^2 = h \;\; \text{in} \;\; (0,1) \\
u(0) = u(1) \end{cases}
\end{equation*}
has a very different behaviour whether the coefficient $a(t)$ changes sign or not. When $a(t) \in C^{0}([0,1]) \setminus \{0\}$ does not change sign it is proved in \cite{M-S} that the map naturally associated to problem  $(\text{P}_{1})$ is a global fold between the spaces $C_{\#}^{1}([0,1])$ and $C^{0}([0,1])$. Hence it is a first order analogue of the celebrated map studied by A. Ambrosetti and G. Prodi in \cite{A-P} (cf. also \cite{B-C}) which was associated to the Dirichlet problem for a second order semilinear PDE with a strictly convex nonlinearity. This implies that both problems can have at most two solutions for any right-hand side $h$. The singular points of the maps defined by these problems are all infinite-dimensional fold points and hence are the analogues of the first type of the singularities studied in the finite-dimensional case by B. Morin (cf. \cite{Mo}). We explicitely point out that Morin singularities can be considered as a reference point for the development of infinite-dimensional Singularity Theory. One of the basic reasons we refer to them is the so-called  Normal Form Theorem for Morin singularities which allows describing the behaviour of a given map near a Morin singularity. By means of the Normal Form Theorem it could be possible to study how the number of solutions to the equation $F(u) = h$ varies when $u$ is near a given Morin singularity $u_0$ for a map $F$ and $h$ is near $h_{0}$ = $F(u_0)$. More generally we can say that the Normal Form Theorem for Morin singularities guarantees the local finiteness of the solutions of the equation near a given Morin singularity.\\
On the other hand the multiplicity of solutions to $(\text{P}_{1})$ is deeply altered when $a(t)$ changes sign. It is proved in \cite{C-D 2} that when $a(t) \in C^{2}([0,1]) \setminus \{0\}$, $a(0) = a(1) = 0$ and $a(t)$ changes sign, there exists a suitable set of right-hand sides $h$ whose preimage consists of exactly one unbounded real-analytic curve of solutions. Of course the points on such curves are singular points for the map associated to the problem and, since there is no  local finiteness, these singular points cannot be infinite-dimensional analogues of the Morin singularities. This raises the question if it is also possible to describe the behaviour of the map near these new singular points. In \cite{C-D 2} a description of the global behaviour is obtained without a deep analysis of the nature of these singular points because the result follows from the rich structure of the Riccati equation (note that the proof of property (6) in \cite{C-D 2} is wrong but becomes correct by replacing $L^{2}(0,1)$ with $C^0$\(([0,1])\) and hence one has to replace everywhere Sobolev spaces with classical spaces).\\
Furthermore, we note that curves of solutions to a fixed right-hand side also appear in second order problems. A remarkable example is provided by the study of the periodic oscillations of a ``short'' pendulum without friction, i.e. the problem 
\begin{equation*}
(\text{P}_{2}) \begin{cases} u'' + A  \sin u = h \;\; \text{in} \;\; (0,T) \\
u(0) = u(T)\\
u'(0)= u'(T) \end{cases}
\end{equation*}
where $T>0$ and $A > (\frac{2\pi}{T})^{2}$. In fact under this assumption on $A$ it is well known that problem $(\text{P}_{2})$ with $h \equiv 0$ admits, because of conservation of energy, at least a non-constant solution $\widetilde{u}$ and hence a whole circle of solutions given by the time-translates of $\widetilde{u}$. 
Note that the existence of curves of solutions is not a phenomenon confined to periodic problems:  this also occurs in other boundary value problems for ODEs and PDEs. A quite natural example is given by the Neumann problem
\begin{equation*}\label{}
(\text{P}_{3})\begin{cases} u'' + u  u' = h \;\; \text{in} \;\; (0,1) \\
u'(0) = 0\\
u'(1)= 0 \end{cases}
\end{equation*}
which (of course) has the straight-line of the constant functions as a curve of solutions to  $h \equiv 0$.
Yet \textit{a priori} one of the curves of solutions to $(\text{P}_{1})$, a circle of solutions to $(\text{P}_{2})$ and the above straight-line of solutions for $(\text{P}_{2})$ could be made up of singular points with a very different nature. Hence it is important to characterize such singular points and understand the differences among them. In fact it is known (cf. \cite{C-D 2}) that fold points are the only singular points near the curves of solutions in $(\text{P}_{1})$ while other kinds of singular points could be identified near the curves described above for  $(\text{P}_{2})$ and $(\text{P}_{3})$. This characterization in turn should become useful in order to describe how the multiplicity of solutions to the equation $F(u) = h$ changes near these singular points.\\
\indent Though we believe that the existence of such curves of solutions is uncommon, we remark that this phenomenon occurs for ODEs and PDEs under quite different boundary conditions. As already said the singular points on these curves cannot be infinite-dimensional analogues of the Morin singularities; in order to analyze them, a strategy completely different from what has been done so far was needed. \\
In order to overcome the difficulties of relating to what was already known in the finite-dimensional setting (as sketched in the above paragraph) we were prompted to look for a new, infinite-dimensional, approach to the study of singularities (at least in the simple case where the involved maps are Fredholm of index zero and have one-dimensional kernels, see below). The main results presented in this series of papers allow us to develop a unified treatment both for the analogue of Morin singularities of order $k$, which we call here $k$-singularities, and for the new classes of singularities related to problems like $(\text{P}_{1})$, $(\text{P}_{2})$ and $(\text{P}_{3})$ which will be called maximal $k$-transverse singularities and $\infty$-transverse singularities. Specifically, the purpose of this first paper is to introduce the classification of the above-mentioned singularities. It is preliminarily important to note that the map naturally associated to all the problems considered above is a smooth Fredholm map of index zero between suitable infinite-dimensional Banach spaces: for this reason we confine our study to this class of maps. In fact, for the sake of simplicity, in these papers we only consider Fredholm maps with one-dimensional kernels and almost all statements are presented for $C^d$ maps with $d \geqslant 2$ as this could be useful for future developments, though our main interest is in $C^\infty$ maps. In contrast with the two known approaches in Singularity Theory, for which we refer to the historical part at the end of this Introduction, we adopted a local approach, partly suggested by the direct study of problems such as $(\text{P}_{1})$. This leads us to consider suitable nonlinear and linear functionals, called ``fibering functionals'', which are inductively defined in a neighbourhood of the studied singularity by means of a suitable analytic tool that we call ``fibering pair'', a sort of smooth representation of the kernel and cokernel of the Fréchet derivative of the map near the singular point. It is indeed interesting to note that the very definitions of the fibering functionals were partially suggested by the study of first order problems (e.g. see \cite{C-D 1}, \cite{C-D 2}). Then the fibering functionals allow us to introduce our classification of the singularities of a smooth Fredholm map of index zero.\\
\indent We would also like to stress that we deliberately adopted an elementary approach to the classification, based on well-known results in nonlinear analysis in infinite dimensions and a minimal amount of tools from differential topology, i.e. the basic notions of Banach manifolds, maps between manifolds and the Submersion Theorem as the more sophisticated result. Moreover, despite a few setbacks (e.g. the somehow lengthy proofs of the invariance of the stratification of singularities with respect to the chosen fibering pair), we hope that our mostly analytical and self-contained approach could be modified and extended to investigate further applications of Singularity Theory to differential problems, such as the study of boundary value problems whose related Fredholm map exhibits singularities with two-dimensional kernels. \\
\indent Finally, we point out that this is the first paper of a series of three. While here we are  mainly concerned with the classification of singularities, the other two articles of the series deal with two important features that have to be considered when one is interested in applying Singularity Theory to the study of nonlinear differential problems. In the second paper we analyze the local behaviour of a map $F$ near the different types of singular points by finding notable information about the number of solutions to the equation $F(u) = h$. Since the conditions given in the classification of singularities are often difficult to verify when studying nonlinear differential problems, it is also necessary to find operative conditions which are easy to use: this constitutes the main objective of the third paper of the series.\\
\indent
The present work is organized as follows (we refer to the end of this Introduction for the table of contents). In Section \ref{ss11} we first recall that, given a regular Fredholm map of index zero $F:X\rightarrow Y,X,Y$ spaces of Banach, a point $u_{o}\in X$ is called a (simple) singularity if dim$\,N(F'(u_{o}))=1$. We then define the main tool that is needed for the classification of singularities, that is the \textit{fibering pair} or \textit{f-pair}$(\varphi,\psi)$ near a singular point $u_{o}$, and we study the relationship between different fibering pairs. Moreover we define the associated families of functionals $J_{k}(\varphi,\psi),\, I_{k}(\varphi,\psi)$ which allow us to classify the singularities in Chapter \ref{s2}. In Section \ref{ss12}, we state and prove the Local Representation Theorem, the basic tool of nonlinear analysis we use here. This also provides the prototype of maps on which some of our proofs are based on (e.g. the  existence of fibering pairs near a singularity). These are maps of the form $G : \mathbb{R}\times Z \rightarrow \mathbb{R}\times Z, G(t,\xi)=(g(t,\xi),\xi)$, and because of their intrinsic importance here and in the next papers (e.g. for the Normal Form Theorem in \cite{B-D 2}) it is convenient to name them \textit{Lyapunov-Schmidt maps} or \textit{LS-maps}. Finally, in Sections \ref{ss13} and \ref{ss14}, we study the f-pairs in the simple case of LS-maps and, in the general case, we look at how f-pairs are affected by changes of coordinates and in which way different f-pairs are related to each other near a fixed singularity: incidentally, this will also prove that f-pairs exist.\\
\indent In Section \ref{ss21}, given a singularity $u_{o}$ and a fixed f-pair $(\varphi,\psi)$ near $u_{o}$, we use the related functionals $J_{k}(\varphi,\psi),\, I_{k}(\varphi,\psi)$ (shortly $J_{k}, \,I_{k}$) to define four possible types of ``behaviour'' as a singular point for $u_{o}: k-transverse\;singularity,\; k-singularity,\;$$ maximal\; k-transverse\; singularity,\; \infty-transverse\; singularity$. In Section \ref{ss22} we study the linear independence of the functionals $I_{1},\ldots,I_{k}$ and we determine when the zero-sets of the functionals $J_{1},\ldots,J_{k}$  are a submanifold of $X$ of codimension $k$, using the \textit{Submersion Theorem} under suitable conditions of \textit{transversality}. In this way, the definition of $k$-transverse singularity proves to be a useful tool for the description and classification of singularities that follows next. In fact, we show that, near a $k$-transverse singularity, there exists a suitable, nested \textit{stratification of manifolds} of singular points (cf. Section \ref{ss24}); in particular, near a 1-transverse singularity the singular set is a one-codimensional submanifold (this is preliminarily shown in Section \ref{ss23}). Thus, under the basic assumption of 1-transversality, i.e. smoothness of the singular set, we partition the set of (simple) singularities for a given map $F$ in three kinds: $k$-singularities, maximal $k$-transverse singularities and $\infty$-transverse singularities (cf. Section \ref{ss25}). This is what we call the \textit{classification} of singularities; in \cite{B-D 2} it is shown that the $k$-singularities are the analogues of Morin singularities of order $k$ in finite dimensions, while maximal $k$-transverse singularities and $\infty$-transverse singularities are introduced here for the first time. The $\infty$-transverse singularities can only occur in infinite dimensions and when the map $F$ is smooth: some general properties of these singularities are also investigated in \cite{B-D 2}. In Section \ref{ss25} we prove that the classification of singularities is well-posed (i.e. independent of the f-pair chosen in the definition) and \textit{invariant} under diffeomorphisms. Finally, in Section \ref{ss26} we consider theoretical examples of all singularities for suitable polynomial-type LS-maps. It is worthwhile to emphasize that we also provide differential examples of all kinds of singularities in Section \ref{ss21}, though a comprehensive study with full proofs will be found in other papers by the same authors. \\
\indent To conclude, we find it useful to recall below the main steps in the development of infinite-dimensional Singularity Theory and its application to differential problems.\\
\indent As is well known this method was originally proposed by A. Ambrosetti and G. Prodi in \cite{A-P}, where a global multiplicity result for a class of second order nonlinear Dirichlet problems was proved. A more geometrical approach to the same problem was then given by M. Berger and P. Church in \cite{B-C}, where it was essentially shown that the notion of ``ordinary singular point'' introduced in \cite{A-P} is the infinite-dimensional generalization of the ``fold point'' considered by H. Whitney in \cite{Wh} for mappings of the plane into the plane. In \cite{Wh}, under suitable assumptions on the second derivative of the considered maps, an explicit construction of local changes of coordinates showed the local equivalence with the fold map  $p(x,y) = (x^{2},y)$. An analogous construction was employed in \cite{B-C}, by using global changes of coordinates, in order to prove that the Fredholm map associated with the differential problem in \cite{A-P} is globally equivalent to an infinite-dimensional fold map. These results prompted H. McKean and J. Scovel to study first and second order problems with a quadratic nonlinearity, cf \cite{M-S}, \cite{Sc}. The authors showed that the first order periodic problem behaves like the one considered in \cite{A-P} and \cite{B-C}. Interestingly, the study of the second order Dirichlet problem appearing in \cite{M-S} and \cite{Sc} can be also regarded as the starting point for the analysis of singularities more complicated than folds. Almost simultaneously M. Berger, P. Church and J. Timourian proved the normal form theorems for fold and cusp singularities, see \cite{B-C-T}. For $C^\infty$ maps the proof follows the approach used by B. Morin in the finite-dimensional case, see \cite{Mo}, which does not require an explicit construction of local changes of coordinates. Approximately at the same time V. Cafagna and F. Donati, inspired by the problem studied in \cite{Sc}, presented in \cite{C-D 1} an example of map whose singularities are only folds and cusps. Such a map is associated with a first order periodic problem where the nonlinearity is cubic. In \cite{C-D 1} (cf also \cite{C-D 2}) a global multiplicity result for the considered problem is obtained by combining the normal form theorems for the singularities with the global behaviour of the map. A similar result, but of a semi-global nature, was then proved by F. Lazzeri and A.M. Micheletti in \cite{L-M} for an asymptotically linear Dirichlet problem. Among the works which show significant examples of cusp singularities we mention \cite{C-T}, \cite{Ruf}, \cite{C-D-T}, \cite{M-S-T} and \cite{C-L} where global results for proper maps are obtained, while \cite{Do} reports on a local result for a non-proper map.
The next step for a better understanding of the local structure of smooth Fredholm maps near a singularity was proposed by J. Damon in \cite{Da}. The algebraic-geometric approach in \cite{Da} is based on the classification of the singularities by means of their local ring, in a way that closely resembles the finite-dimensional case; this allows studying \textit{stable} singularities, such as the Morin singularities. Another important paper devoted to the extention of the algebraic singularity theory to infinite dimensions was then written by P. Church and J. Tmourian, cf \cite{C-Ti}. We remark that the maximal k-transverse singularities are not stable, i.e. a slight perturbation of the map can locally alter the nature of the singularity (cf. Remark \ref{Rem270}). The stability and other features of $\infty$-transverse singularities are  discussed in \cite{B-D 2}. We finally point out a conjecture by B. Ruf, cf. \cite{Ruf bis}, on the existence of elliptic boundary value problems with Neumann conditions whose associated maps are \textit{global} $k$-singularities for $k$ even. This conjecture has some affinities with problem $(\text{P}_{6})$ in Section \ref{ss21} and it was partially proved, locally near the origin, in \cite{C-Ti}.

\tableofcontents
\vspace{10 pt}
\section{Fibering pairs}
\subsection{Basic Definitions} \label{ss11}
\subsubsection{} Let $U$, $V$ be open subsets of the real, infinite-dimensional Banach spaces (or $B$-spaces) $X$, $Y$ and let $F: U \subseteq X \rightarrow V \subseteq Y$ be a $C^{d}$ map, $d = 1, 2, 3, ...$ or $\infty$ or $\omega$ where $C^{\omega}$ refers to real analytic maps. We recall that $F$ is a \textit{Fredholm map of index} $0$, or simply a $0$-\textit{Fredholm map}, if dim$\,N(F'(u)) = \text{dim}\,Y/R(F'(u)) < +\infty$, $\forall \, u \in U$, where $F'(u) \in L(X,Y)$ is the Fréchet derivative of $F$ at $u$, $L(X,Y)$ being the $B$-space of the bounded linear operators from $X$ into $Y$, and $N(F'(u))$, $R(F'(u))$, $Y/R(F'(u))$ are the \textit{kernel}, \textit{range} and \textit{cokernel} of $F'(u)$ respectively. Since $F'(u)$ is a continuous operator and codim$\,R(F'(u)) := \text{dim}\,Y/R(F'(u)) < +\infty$ then $R(F'(u))$ is closed in $Y$ (cf. \cite{Ze1}, proposition 8.14). Let $Y^\ast$ be the (topological) dual space of $Y$ and let $R(F'(u))^\perp := \{\gamma \in Y^\ast: \gamma(h) = 0 \; \forall \; h \in R(F'(u))\}$. Since $R(F'(u))^\perp \cong (Y/R(F'(u)))^\ast$ (cf. \cite{Rud}, theorem 4.9) and dim$\, Y/R(F'(u)) < +\infty$, then $R(F'(u))^\perp$ is isomorphic to the cokernel  $Y/R(F'(u))$. Thus dim$\, R(F'(u))^\perp = \text{codim}\,R(F'(u))$.\\
\indent Addtionally, it is useful to consider the \textit{adjoint  operator} $ F'(u)^\ast : Y^\ast \rightarrow X^\ast$ of $F'(u)$, $u \in U$, which is defined as $F'(u)^\ast (\gamma) := \gamma \circ F'(u)$ for all $\gamma \in Y^\ast$. \\
It is well known (cf. \cite{Br}, corollaire II.17 and théorème II.18) that $ F'(u)^\ast\in L(Y^\ast,X^\ast),\\
N(F'(u)^\ast) = R(F'(u))^\perp$ and, being $R(F'(u))$ closed, $R(F'(u)^\ast) = N(F'(u))^\perp = \{\delta \in X^\ast : \delta(u) = 0, \forall \, u \in N(F'(u))\}$. Since $X^\ast /N(F'(u))^\perp \cong N(F'(u))^\ast$ (cf. \cite{Rud}, theorem 4.9), then $X^\ast/N(F'(u))^\perp$ is isomorphic to the kernel  $N(F'(u))$ and so codim$\,N(F'(u))^\perp = \text{dim}\, N(F'(u))$. Summarizing:
\begin{equation*}
\begin{split}
\text{dim}\,N(F'(u)^\ast) = \text{dim}\,R(F'(u))^\perp = \text{codim}\,R(F'(u)) = \\
=\text{dim}\,N(F'(u)) =\text{codim}\,N(F'(u))^\perp =\text{codim}\,R(F'(u)^\ast).
\end{split}
\end{equation*}
\begin{definition}
\label{d112} Let $F : U \subseteq X \rightarrow V \subseteq Y$ be a $C^d$ $0$-Fredholm map between open subsets $U$, $V$ of the $B$-spaces $X$, $Y$. The point $u \in U$ is said to be a \textit{singularity} for the map $F$ if dim$\,N(F'(u)) \geq 1$ and a \textit{simple singularity} for $F$ if dim$\, N(F'(u)) = 1$. The set of all singularities will be denoted by $S(F) := \{u \in U : \text{dim}\, N(F'(u)) \geq 1\}$, while $S_1(F ):= \{u \in U : \text{dim}\,N(F'(u)) = 1\}$ is the subset of simple singularities.
\end{definition}
\indent Since we will only be concerned with simple singularities we will often refer to them as \textit{singularities}. In the same way we shall write $S_1$ instead of $S_1(F)$ if there is no risk of confusion. Note that if $u \in S_1(F)$  then
\begin{equation*}
\text{dim}\,N(F'(u)) = \text{codim}\,R(F'(u)) = \text{dim}\,R(F'(u))^\perp = \text{dim}\,N(F'(u)^\ast) = 1.
\end{equation*}\\
\indent
We now introduce a basic tool for the study of singularities for Fredholm maps, which was inspired by the arguments developed in the unpublished paper \cite{C-D 3}. In fact, this notion is useful to classify the singular points and, as we show in Chapter 1 of \cite{B-D 3}, it also provides an operative characterization of singularities. 

\begin{definition}
\label{d113} Let $F : U \subseteq X \rightarrow V \subseteq Y$ be a $C^d$ 0-Fredholm map between open subsets $U$, $V$ of the $B$-spaces $X$, $Y$ such that dim$\,N(F'(u))\leq 1$, $\forall \, u \in U$. We say that $(\varphi,\psi)$ is a \textit{fibering pair}, or \textit{f-pair}, for $F$ on $U$ if:
\begin{align*}
&\text{i)} \;\, \varphi \in C^{d-1}(U,X \setminus \{0\}),\; \psi \in C^{d-1}(U,Y^\ast \setminus \{0\});\\
&\text{ii)} \; \forall \, u \in S_1(F) \Rightarrow \varphi(u)\in N(F'(u)),\;\psi(u)\in R(F'(u))^\perp = N(F'(u)^\ast).
\end{align*}
\end{definition}

The mappings $\varphi$ and $\psi$ will be also called \textit{kernel} and \textit{cokernel fibering maps} for $F$, respectively.\\
\par
The existence of a fibering pair for $F$, even on a suitable neighbourhood $U$ of a simple singularity, is not obvious. We postpone the proof to Section 1.4.  \\

\begin{remark}
We refer to $\varphi$ and $\psi$ as \textit{fibering maps} for the following reason. If $S_1(F)$ is a submanifold of $X$ then it can be shown that the \textit{set sum} or  \textit{disjoint union} $\mathbf{\bigvee}_{u\in S_1(F)} N(F'(u))$  is a  $C^{d-1}$ \textit{vector bundle} with \textit{base space} $S_1(F)$ and one-dimensional \textit{fibers}  $N(F'(\cdot))$. Hence the map $\varphi$ is, at least locally, a \textit{non-vanishing section} of the bundle, i.e. $\varphi(u)$ is a non-zero vector generating the fiber  $N(F'(u))$. An analogous statement is true for the set sum $\mathbf{\bigvee}_{u\in S_1(F)} R(F'(u))^\perp$  and the map $\psi$. In fact it can be proved that such bundles are $C^{d-1}$ \textit{pull-back} vector bundles, \textit{via} $F'$, of suitable one-dimensional analytic vector bundles in $L(X,Y)$. For a proof of these statements we refer to \cite{Ba3}, where another proof of the local existence of fibering maps is given.
\end{remark}
\vspace{4pt}
The example below shows that there can be several f-pairs for a given map $F$. Hence it is important to know how different f-pairs are related to each other. We refer to Sections \ref{ss13} and \ref{ss14} for results on this subject.\\

\begin{example}\label{Ex115}
Let us consider the following problem:
\begin{equation*}
(\text{P}_{4}) \begin{cases} u' + a(t)g(u) = h \;\; \text{in} \;\; (0,1) \\
u(0) = u(1), \end{cases}
\end{equation*}
which generalizes the model problem $(\text{P}_{1})$ mentioned in the Introduction. Here $a \in C^0([0,1]) \setminus \{0\}, h \in C^0([0,1]), g \in C^\infty(\mathbb{R})$  and $u \in C_{\#}^{1}([0,1]):= \{u \in C^1([0,1]) : u(0) = u(1)\}$. Hence we can associate a map $F$ with problem $(\text{P}_{4})$, namely the $0$-Fredholm map $F$ of class $C^\infty,F : C_{\#}^{1} ([0,1]) \rightarrow C^0([0,1])$, defined as $F(u) = u' + a(t)g(u)$. Then it is not difficult to check that for $v \in C_{\#}^{1}([0,1])$ one has $F'(u)v = v' + a(t)g'(u)v$  and  $S_1(F) = \{u \in C_{\#}^{1}([0,1]): \int_0^1 a(t)g'(u(t))dt =0\}$. Moreover, for all $u \in S_1(F)$ we have:
\begin{equation*}
N(F'(u)) = \{v \in C_{\#}^{1}([0,1]):v = c\cdot \text{exp}[- \int _0^t a(\tau)g'(u(\tau))d\tau],\;\; c \in \mathbb{R}\},
\end{equation*}
and
\begin{equation*} 
R(F'(u)) = \{h \in C^0([0,1]): \int_0^1 h(t)\cdot \text{exp}[\int _0^t a(\tau)g'(u(\tau))d\tau]dt = 0\}.
\end{equation*}
\end{example}
It is also easy to verify that the pairs
\begin{align*}
&(\varphi(u), \psi(u))=\\
&=\big(\text{exp}[t\cdot \int_0^1 ag'(u)]\cdot \text{exp}[-\int_0^t\! ag'(u)], \text{exp}[-t\cdot \int_0^1 ag'(u)]\cdot \text{exp}[\int_0^t ag'(u)]\big),\\
&(\widetilde{\varphi}(u),\widetilde{\psi}(u))=\big(\text{exp}[t\cdot \int_0^1 ag'(u)]\cdot \text{exp}[- \int_0^t ag' (u)], \text{exp}[\int_0^t ag'(u)]\big)\\
&(\overline{\varphi}(u), \overline{\psi}(u))=\big(\{1-t+t\cdot \text{exp}[\int_0^1 ag'(u)]\}\cdot \text{exp}[-\int_0^t ag'(u)], \text{exp}[\int_0^t ag'(u)]\big)
\end{align*}
are all global f-pairs for $F$, i.e. f-pairs defined for any $u\in C_{\#}^{1}([0,1])$. Note that while the kernel fibering maps $\varphi(u), \widetilde{\varphi}(u),\overline{\varphi}(u)$ are all elements of $C_{\#}^{1}([0,1])$, the cokernel fibering maps $\psi(u), \widetilde{\psi}(u),\overline{\psi}(u)$ are identified with their respective representatives. For example, the map $\psi(u) \in C^0([0,1])^\ast$ is explicitly given by the functional
\begin{equation*}
\psi(u)h = \int_0^1 \text{exp}[-t\cdot \int_0^1 ag'(u)]\cdot \text{exp}[\int_0^t ag'(u)]\cdot h(t) dt, \; \,  h \in C^0([0,1]).
\end{equation*}
Similar conclusions can be obtained for the general problem
\begin{equation*} 
(\text {P}_{4}^{\,\prime})\begin{cases} u' + g(t,u) = h \; \; \text{in} \;\;(0,1) \\
u(0) = u(1) \end{cases} 
\end{equation*}
where $g \in C^\infty([0,1]\times \mathbb{R})$, provided one defines $F(u):= u' + g(\cdot ,u)$ and replaces $ag'(u)$ with $\frac{\partial g}{\partial u}(t,u)$ everywhere. For instance, the first pair defined above naturally extends to the f-pair
\begin{align*}
&(\varphi(u),\psi(u))=\\
&\big(\text{exp}[t\cdot\int_0^1\frac{\partial g}{\partial u}(t,u)]\cdot \text{exp}[-\int_0^t \frac{\partial g}{\partial u}(t,u)], \text{exp}[-t\cdot \int_0^1\frac{\partial g}{\partial u}(t,u)]\cdot \text{exp}[\int_0^t\frac{\partial g}{\partial u}(t,u)]\big).
\end{align*}
\par
We are now able to introduce the main tools that will be used in the following to give the classification of singularities of a 0-Fredholm map. We refer to Chapter 2 for a comprehensive study of this classification.\\
\begin{definition}\label{D116}
Let  $F : U\subseteq X \rightarrow V \subseteq Y$ be a $C^d$ 0-Fredholm map between open subsets $U, V$ of the $B$-spaces $X, Y$ and let $(\varphi,\psi)$ be a fibering pair for $F$ on $U$. We define
\begin{equation*}
J_0(\varphi,\psi) \equiv J_0:U\rightarrow\mathbb{R}\; \text{  by  }\;J_0(\varphi,\psi)(u):=\psi(u)F'(u)\varphi(u), u \in U, 
\end{equation*}
and inductively, $\forall \ k \in \{1,\ldots,d-1\} \;\text{and} \; \forall \  u \in U$,
\begin{equation*}
I_k(\varphi,\psi)(u) \equiv I_k(u):X \rightarrow\mathbb{R}\; \text{   by   }\;I_k(\varphi,\psi)(u):= J'_{k-1}(\varphi,\psi)(u)\in X^\ast
\end{equation*}
and  
\begin{equation*}
J_k(\varphi,\psi) \equiv J_k : U \rightarrow\mathbb{R}\;\text{   by   }\; J_k(\varphi,\psi)(u):=I_k(u)\varphi(u) = J'_{k-1}(u)\varphi(u).
\end{equation*} 
All these maps are called \textit{fibering functionals} for $F$ (which, of course, depend on the pair $(\varphi,\psi)$).\\

For the regularity of the maps  $J_k, I_k$  we have that
\begin{equation*}
J_k(\varphi,\psi) \in C^{d-k-1}(U,\mathbb{R}) \text{  for  } k = 0, 1,\ldots,d-1
\end{equation*}
and
\begin{equation*}
I_k(\varphi,\psi) \in C^{d-k-1}(U,X^\ast) \text{  for  } k = 1,\ldots,d-1.
\end{equation*}
\end{definition}
\vspace{4pt}
Since the approach we adopted to define the functionals $J_k$ is just a special case of a more general procedure, which we will use again later, it is worthwhile to give the following\\
\begin{definition}\label{D117}
Let $X, Z$  be  $B$-spaces, $U$ an open subset of  $X, G : U \subseteq X \rightarrow Z$  a  $C^d$  map, $d \geq 1$, and  $\xi:U \subseteq X \rightarrow X$  a  $C^{d-1}$  map. For each integer $k,0 \leq k \leq d$, we define inductively the $C^{d-k}$  maps  $G_k(\xi) \equiv G_k : U \subseteq X \rightarrow Z$  as   
\[G_0 :=G, G_k := G'_{k-1}\xi, 1\leq k \leq d,\text{ i.e. } G_k(u) := G'_{k-1}(u)\xi(u), u \in U.\]
The maps $G_k$ are called \textit{iterated derivatives} of  $G$  along  $\xi$ or \textit{iterated} $\xi-$\textit{derivatives} of $G$, sometimes abbreviated as $\xi-$\textit{derivatives} of $G$.
\end{definition}
\vspace{4pt}
We can also look at the maps $G_k$ in terms of Vector Field Theory, thus justifying their name too. Indeed, $\xi$ is a vector field on $U$ and the map $G_k$  is simply obtained by deriving $G_{k-1}$ along $\xi$. Then  $G_k$ is the so-called  \textit{Lie Derivative} of  $G_{k-1}$  with respect to the vector field  $\xi$,  i.e. $G_k =G'_{k-1}\xi:=\mathcal{L}_\xi G_{k-1}$  (cf. \cite{La}, chapter V, §2). It is now clear that the functionals $J_k :J'_{k-1}(\varphi,\psi)\varphi, 1\leq k \leq d-1$, introduced in Definition \ref{D116}, are  $\varphi \mbox{-derivatives}$ of $J_0(\varphi,\psi) = \psi F'\varphi$, i.e. $J_k(\varphi,\psi) = \mathcal{L}_\varphi J_{k-1}(\varphi,\psi) = (J_0(\varphi,\psi))_k$.\\
\par
It may now be interesting to see, at least for the easy case considered in Example \ref{Ex115}, how the fibering functionals depend on the f-pairs and how they are related to the set of simple singularities.\\
\begin{example}
\label{Ex118}Let us consider the smooth 0-Fredholm map $F: C_{\#}^{1}([0,1])
\rightarrow C^0([0,1])$, defined as $F(u) = u' + a(t)g(u)$ and associated to problem $(\text{P}_{4})$ (as shown in Example \ref{Ex115}). Then we can explicitly write the first fibering functional, corresponding to each of the f-pairs given in Example \ref{Ex115}, in the following way:\\
\begin{align*}
& J_0(\varphi,\psi)(u)\equiv J_0(u)=\psi(u)F'(u)\varphi(u)=  \\
&=\int_0^1\{\text{exp}[-t\cdot\int_0^1 ag'(u)]\cdot \text{exp}[\int_0^t ag'(u)]\}\cdot\{\text{exp}[t\cdot\int_0^1 ag'(u)]\cdot \\
&\cdot \text{exp}[-\int_0^t ag'(u)]\cdot\int_0^1 ag'(u)\}dt= \int_0^1 ag'(u),\\
\\
&J_0(\widetilde{\varphi},\widetilde{\psi})(u) \equiv \widetilde{J}_{0}(u)=\widetilde{\psi}(u)F'(u)\widetilde{\varphi}(u) = \\
&= \int_0^1\{\text{exp}[\int_0^t ag'(u)]\}\cdot \{\text{exp}[t\cdot \int_0^1 ag'(u)]\cdot \text{exp}[-\int_0^t  ag'(u)]\cdot \int_0^1 ag'(u)\}dt = \\
&= -1 + \text{exp}[\int_0^1 ag'(u)],\\
\\
&J_0(\overline{\varphi},\overline{\psi}) (u) \equiv \overline{J}_{0}(u) = \overline{\psi}(u)F'(u)\overline{\varphi}(u) =\\
&= \int_0^1\{\text{exp}[\int_0^t ag'(u)]\}\cdot \{(-1 + \text{exp}[\int_0^1 ag'(u)])\cdot \text{exp}[-\int_0^t ag'(u)]\}dt =\\
&= - 1 + \text{exp}[\int_0^1 ag'(u)].
\end{align*}
\end{example}
\vspace{4pt}
Although the functional $J_0$  is quite different from  $\widetilde{J}_0 \equiv \overline{J}_0$, it is evident that all three functionals coincide on $S_1(F) = \{u \in C_{\#}^{1}([0,1]): \int_0^1  a(t)g'(u(t))dt = 0\}$ since they are identically zero on $S_1(F)$. Note that in this case $S_1(F)$ coincides with the zero-sets of the functionals $J_0, \widetilde{J}_0, \overline{J}_0$; in general it is only possible to say that $S_1(F)$ is a subset of the zero-set of the first fibering functional (see Example \ref{Ex136} below).\\
\indent
We can finally remark that, if one replaces $ag'(u)$ with $\frac{\partial g}{\partial u}(t,u)$, similar formulas hold for the map $F(u) = u' + g(\cdot,u)$ associated to problem $(\text{P}_4^{\prime})$: for instance, if one naturally modifies $(\varphi,\psi)$ as suggested at the end of Example \ref{Ex115} we obtain that $J_0(\varphi,\psi)(u) = \int_0^1 \frac{\partial g}{\partial u}(t,u)$. In Section \ref{ss21} we will study a concrete problem of the form $(\text{P}_4^{\prime})$ and explain how to compute the next iterated functionals $I_{k}$ and $J_{k}, k=1,2,3$, starting from the above $J_0(\varphi,\psi)$.\\
\par
Moreover Example \ref{Ex118} also suggests that it is important to understand the relations between fibering functionals corresponding to different f-pairs. For this we refer to Sections \ref{ss22} and \ref{ss23}. \vspace{2pt}
\subsection{The Local Representation Theorem}\label{ss12}
\subsubsection{} Let $U, U', V, V'$ be open subsets in the  $B$-spaces $X$, $X'$, $Y$, $Y'$, respectively, and let $F : U \subseteq X \rightarrow V \subseteq Y$, $\varPhi : U' \subseteq X' \rightarrow V' \subseteq Y'$ be mappings of class  $C^d\, (d = 1, 2, 3, \ldots \mbox{ or } \infty \mbox{ or }\omega)$. We recall that $F$ and $\Phi$ are $C^d$ \textit{locally equivalent} near $u_o \in U $ and $u'_o \in U'$ if there exist diffeomorphisms $\alpha$, $\beta$ of class $C^d$ defined on some neighbourhoods of $u_o$ and $F(u_o)$, respectively, such that $\alpha(u_o) = u_o^\prime$ and near $u_o$ one has $\beta F = \Phi \alpha$, i.e. the following diagram locally commutes:
$\\$
\begin{center}
\begin{tabular}{ccc}
&$F \quad $&  \medskip \\
$u_o \in U \subseteq X \quad $&$ \rightarrow \quad $&$ V \subseteq Y $ \medskip  \\
$\alpha \downarrow $&&$ \downarrow \beta$ \medskip \\ 
$u'_o \in U' \subseteq X' \quad $&$ \rightarrow  \quad $&$ V' \subseteq Y'$ . \medskip \\
&$\varPhi \quad $&\\
\end{tabular}
\end{center}

The above diagram is referred to as a \textit{local commutative diagram of class} $C^d$ ($C^d$ \textit{l.c.d.} in short).\\
\indent We also recall that a $C^d$ map $F : U \subseteq X \rightarrow V \subseteq Y$, where $U, V$  are open subsets of the $B$-spaces $X,Y$y, is said to be \textit{double-splitting} at $u_o \in U$ if $N(F'(u_o)), R(F'(u_o))$ are closed subspaces with closed complementary subspaces in $X$ and $Y$  respectively. Of course smooth 0-Fredholm maps are double-splitting at each point of the domain; a quite important property of double-splitting maps is given by the following theorem (cf. \cite{A-M-R}, theorem 2.5.14). We give a direct proof of this result because it is a fundamental tool to prove the existence of f-pairs and we will use it extensively in this series of papers.\\
\begin{theorem}
\label{Teo122}(Local Representation Theorem). \textit{Let} $U,V$ \textit{be open subsets of the}  $B$-spaces $X,Y$\textit{ respectively, let}  $F : U \subseteq X \rightarrow V \subseteq Y$ \textit{be a }$C^d$ \textit{map which is double-splitting at}  $u_o \in U$ \textit{and let} $X_o, Y_o$ \textit{be closed complements of} $N(F'(u_o)), R(F'(u_o))$. \textit{Then there exists a} local commutative diagram \textit{of class} $C^d$
\begin{center}
\begin{tabular}{cccc}
&$F$&  \medskip \\
$\hspace{3 cm} u_o \in U  \subseteq X $&$ \rightarrow $&$\hspace{-1cm} V \subseteq Y $ \medskip  \\
$\text{(D)} \hspace{5 cm} \alpha \downarrow \hspace{2 cm}  $&&$\hspace{-0.5cm} \downarrow \beta$ \medskip \\ 
$\hspace{2 cm}\quad (0,0)\in N(F'(u_o)) \times R(F'(u_o))$&$ \rightarrow $&$ \quad Y_o \times  R(F'(u_o)),$ \medskip \\
&$\varPhi$&\\
\end{tabular}
\end{center}
\textit{with }$\alpha$, $\beta$, $\Phi$ \textit{depending on }$X_o, Y_o$ \textit{and such that }$\varPhi$ \textit{has the form }$\Phi(n,r)=(f(n,r),r),$   \textit{for} $(n,r)$
\textit{near} $(0,0)=\alpha(u_o),$ \textit{where} $f:N(F'(u_o)) \times R(F'(u_o)) \rightarrow Y_o$ \textit{is a  suitable} $C^d$ \textit{map defined near} $(0,0)$ \textit{such that} $f(0,0)=0$ \textit{and} $f'(0,0)=0$.
\end{theorem}
\vspace{4pt}
\textbf {Proof.} Let $p \in L(X,X), \pi \in L(Y,Y)$ be the projections of $X = N(F'(u_o)) \oplus X_o$ and $Y=R(F'(u_o)) \oplus Y_o$ on $N(F'(u_o))$ and $R(F'(u_o))$ respectively. We define the map
\begin{equation*}
\alpha : U\subseteq X \rightarrow N(F'(u_o)) \times R(F'(u_o))
\end{equation*}
as
\begin{equation*}
\alpha(u):=(p u-p u_o, \pi F(u)-\pi F(u_o)),\medskip
\end{equation*}
which is a $C^d$ map with $\alpha(u_o)=(0,0)$. Let us show that
\begin{equation*}
\alpha'(u_o):X \rightarrow N(F'(u_o)) \times R(F'(u_o))
\end{equation*}\medskip
is an isomorphism. Clearly $\alpha'(u_o)v=(p v,\pi F'(u_o)v)=(p v, F'(u_o)v),\, \forall \, v \in X$; moreover, $F'(u_o)\!\!\mid_{X_o}: X_o\rightarrow R(F'(u_o))$  is bijective by construction and it is easy to see that  $\alpha'(u_o)^{-1}$ is given by the map \medskip
\begin{equation*}
(n,r) \in N(F'(u_o))\times R(F'(u_o))\mapsto n + (F'(u_o)\mid_{X_o})^{-1}r.\medskip
\end{equation*}
By virtue of Banach's Open Mapping Theorem $\alpha'(u_o) \in GL(X,N(F'(u_o))\times R(F'(uo)))$ (where $GL(X,Y)$ is the group of invertible continuous linear maps from $X$ to $Y$), hence $\alpha$ is a  $C^d$  diffeomorphism near $u_o$ thanks to the Inverse Function Theorem.
Finally, let us define the map
\begin{equation*}
\beta:Y\rightarrow Y_o\times R(F'(u_o))
\end{equation*}
as
\begin{equation*}
\beta(y):=((\mathbf{1_Y}-\pi)[y-F(u_o)], \pi[y-F(u_o)]),
\end{equation*}\medskip
which is affine, hence of class $C^\omega$.\\
Let us introduce $\varPhi:=\beta F\alpha^{-1}$, which is well-defined near  $\alpha(u_o)=(0,0)$, and let us show that such a  $C^d$ map has the required properties.\\
By definition $\alpha(\alpha^{-1}(n,r))=(n,r)$, i.e. $(p \alpha^{-1}(n,r)-p u_o,\pi F(\alpha^{-1}(n,r)) -\pi F(u_o))=(n,r)$, therefore $\pi F(\alpha^{-1}(n,r))- \pi F(u_o)=r$. From this we have  
\begin{align*}
&\varPhi(n,r)=\beta(F(\alpha^{-1}(n,r))) =\\
& \qquad \quad =((\mathbf{1_Y}-\pi)[F(\alpha^{-1}(n,r))-F(u_o)],\pi[F(\alpha^{-1}(n,r))-F(u_o)]) =\\
& \qquad \quad =((\mathbf{1_Y}-\pi)[F(\alpha^{-1}(n,r)) -F(u_o)],r).
\end{align*}
Hence it suffices to define  $f(n,r):=(\mathbf{1_Y}-\pi)[F(\alpha^{-1}(n,r))-F(u_o)]$ for $(n,r)$ near $\alpha(u_o)=(0,0)$. Then $f$ is a $C^d$ map such that  $f(0,0)=0$  and
\begin{equation*}
f'(0,0) = (\mathbf{1_Y}-\pi)F'(\alpha^{-1}(0,0))(\alpha^{-1})'(0,0)=(\mathbf{1_Y}-\pi)F'(u_o)(\alpha'(u_o))^{-1}=0 ,
\end{equation*}
since  $\mathbf{1_Y}-\pi$ is the projection on  $Y_o\ldotp \bracevert$\\
\par
\begin{remark}
\label{Rem123} When $F$ is a 0-Fredholm map and $u_o \in S_1(F)$ then 

\hspace{4 cm} dim$\, N(F'(u_o))=\text{codim}\,R(F'(u_o))=1$.\\
Since $R(F'(u_o))$ is a closed subspace of  $Y$, then $Z:=R(F'(u_o))$ is a $B$-space. If, as above, $X_o, Y_o$ are closed complements of $N(F'(u_o)), R(F'(u_o))$ then dim$\,Y_o=1$ and there exist isomorphisms $\sigma :N(F'(u_o))\rightarrow \mathbb{R},\tau :Y_o\rightarrow \mathbb{R}$. From Theorem \ref{Teo122} and by using the isomorphisms $\sigma,\tau$ it is easy to deduce the $C^d$  l.c.d. \medskip
\begin{center}
\begin{tabular}{ccc}
&$F \quad $&  \medskip \\
$u_o \in U \subseteq X \quad $&$ \rightarrow \quad $&$ V \subseteq Y $ \medskip  \\
$\gamma \downarrow $&&$ \downarrow \delta$ \medskip \\ 
$(0,0)\in \mathbb{R} \times Z \quad $&$ \rightarrow  \quad $&$ \mathbb{R} \times Z$  \medskip \\
&$\varPhi \quad $&\\
\end{tabular}
\end{center} \medskip
where $\gamma, \delta, \varPhi$ depend on $X_o, Y_o, \sigma,\tau$ and $\varPhi$ has the form  $\varPhi(t,z)=(f(t,z),z)$, for $(t,z)$ near $(0,0) =\gamma(u_o)$. As seen above, $f:\mathbb{R}\times Z \rightarrow \mathbb{R}$ is a $C^d$ function defined near $(0,0)$ such that $f(0,0)=0$ and $f'(0,0)=0$.
\end{remark}
\vspace{4pt}
In the rest of the paper we will often use maps such as $\varPhi$; it is thus convenient to introduce the following
\begin{definition}
Given a $B$-space $Z$, we say that a $C^d$ map $G:U\subseteq \mathbb{R}\times Z\rightarrow V\subseteq \mathbb{R}\times Z$ is a \textit{Lyapunov-Schmidt map} or \textit{LS-map} if it is of the form $G(t,\xi)=(g(t,\xi),\xi)\;\, \forall \,(t,\xi) \in U$, where $g:U\subseteq \mathbb{R}\times Z \rightarrow \mathbb{R}$ is a $C^d$ function.
\end{definition} \vspace{4pt}

\subsection{Fibering Pairs for LS-maps}\label{ss13}
\subsubsection{}\label{sss131} Remark (\ref{Rem123}) suggests that a possible way to construct an f-pair for the considered map $F$ could be to obtain an f-pair for the simpler map $\varPhi$ and then, through suitable diffeomorphisms, to lift it to an f-pair for $F$. This approach can be actually carried out as we will show in Section \ref{ss14}; of course this strategy requires a preliminary study of the f-pairs for LS-maps, which is why the current section lists  some of the properties of these f-pairs.\\
\par
Let $U, V$ be open subsets in $\mathbb{R}\times \varXi$, where $\varXi$ is a $B$-space, and let $F:U\subseteq \mathbb{R} \times \varXi \rightarrow V\subseteq \mathbb{R}\times \varXi$ be a $C^d$  LS-map, i.e. such that $F(t,\xi)=(f(t,\xi),\xi),\; \forall \,(t,\xi)\in U$, with $f:U\subseteq \mathbb{R} \times \varXi \rightarrow \mathbb{R}$ a $C^d$ function. Let us show that: 
\begin{equation}\label{131}
F \text{ is a  0-Fredholm  map and dim}\,N(F'(t,\xi))\leqslant 1,\;\, \forall \,(t,\xi) \in U.
\end{equation}
Note that the Fréchet derivative of $F$ at $(t,\xi) \in U$, i.e. $F'(t,\xi):\mathbb{R}\times \varXi \rightarrow \mathbb{R}\times \varXi$, can be expressed in a matricial form as
\begin{equation*}
F'(t,\xi)= \left[ \begin{array}{cc}
\frac{\partial f}{\partial t}(t,\xi)&\frac{\partial f}{\partial\xi}(t,\xi)\\\\

0 &\mathbf{1_\varXi}
\end{array} \right].
\end{equation*}
In other words,$\;\; \forall \;\;(r,v) \in \mathbb{R}\times \varXi$ one has
\begin{equation*}
F'(t,\xi)(r,v) = \left[ \begin{array}{cc}
\frac{\partial f}{\partial t}(t,\xi)&\frac{\partial f}{\partial\xi}(t,\xi)\\\\
0 &\mathbf{1_\varXi}
\end{array} \right]
\left[ \begin{array}{c}
r \\
\\
v \end{array} \right]
= (r \frac{\partial f}{\partial t}(t,\xi) + \frac{\partial f}{\partial \xi}(t,\xi)v,v).
\end{equation*}
\\
Here and in the following we consider $\frac{\partial f}{\partial t}(t,\xi)$ as an element of $\mathbb{R}$, thanks to the natural isomorphism  $L(\mathbb{R},\mathbb{R}) = \mathbb{R}^\ast \cong \mathbb{R}$ given by $\varphi \in L(\mathbb{R},\mathbb{R})\mapsto \varphi(1)\in \mathbb{R}$, while 
$\frac{\partial f}{\partial\xi}(t,\xi) \in L(\varXi,\mathbb{R})=\varXi^{\ast}$. It is thus clear that $F'(t,\xi)$ is an isomorphism \textit{iff} $\;\frac{\partial f}{\partial t}(t,\xi) \neq 0$. Note that when $\frac{\partial f}{\partial t}(t,\xi)=0$ we obtain
\begin{equation}\label{132}
\begin{split}
N(F'(t,\xi))&=\mathbb{R}\times\{0\},\;\;\;  \text{where}\;\;\; 0 \in \varXi,\\
R(F'(t,\xi))&=\{(\frac{\partial f}{\partial\xi}(t,\xi)v,v):v\in \varXi\} .
\end{split}
\end{equation}
Moreover, $R(F'(t,\xi))$ is complemented in $\mathbb{R}\times \varXi$  by $\mathbb{R}\times \{0\}$ because every $(r,v) \in \mathbb{R}\times \varXi$ can be written in a unique way as $(r,v)=(\frac{\partial f}{\partial\xi}(t,\xi)v,v)+(r-\frac{\partial f}{\partial\xi}(t,\xi)v,0)$ and hence codim$\,R(F'(t,\xi))=1$. Therefore, $\;\, \forall \,(t,\xi) \in U$ we have that either $F'$ is an isomorphism or dim$\,N(F'(t,\xi))=\text{codim}\,R(F'(t,\xi))=1$, hence (\ref{131}) holds.\\
\par
In the above proof we also deduced that $F$ can only admit simple singularities and $(t,\xi)$  is a singular point if and only if $\frac{\partial f}{\partial t}(t,\xi)=0$. Consequently
\begin{equation}\label{133}
S_1(F)=\{(t,\xi) \in U :\frac{\partial f}{\partial t}(t,\xi)=0\}.
\end{equation}
\par
By virtue of (\ref{132}) we shall construct, in a natural way, a (global) $C^{d-1}$  f-pair for $F$ on $U$: we call such a pair the \textit{canonical fibering pair} $(\varphi_C,\psi_C)$ for $F$ (of course when $F$ is an LS-map).\\
For each $(t,\xi)\in U$ we define $\varphi_C(t,\xi)\in \mathbb{R}\times\varXi$ and the functional $\psi_C(t,\xi):\mathbb{R}\times \varXi \rightarrow \mathbb{R}$ as
\begin{equation}\label{134}
\begin{split}
\varphi_C(t,\xi)&:=(1,0), 1 \in \mathbb{R}, 0 \in \varXi,\\
\psi_C(t,\xi)(r,v)&:=r-\frac{\partial f}{\partial\xi}(t,\xi)v,(r,v) \in \mathbb{R}\times \varXi.
\end{split}
\end{equation}
It is clear that $\psi_C(t,\xi) \in (\mathbb{R}\times \varXi)^\ast$. On the other hand, by means of the natural isomorphisms $(\mathbb{R}\times \varXi)^\ast \cong \mathbb{R}^\ast \times \varXi^\ast \cong \mathbb{R}\times \varXi^\ast$, we can think of $\psi_C(t,\xi)\in \mathbb{R}\times \varXi^\ast$ as given by $(1,-\frac{\partial f}{\partial \xi}(t,\xi))\in \mathbb{R}\times \varXi^\ast$. Then by the very definition $\varphi_C$ and $\psi_C$ are non-zero maps from $U$ into $\mathbb{R}\times \varXi$ and $(\mathbb{R}\times \varXi)^\ast\cong\mathbb{R}\times \varXi^\ast$, respectively. Furthermore, $\varphi_C \in C^{d-1}(U,\mathbb{R}\times \varXi)$ and $\psi_C \in C^{d-1}(U,(\mathbb{R}\times \varXi)^\ast)$ by construction.\\
Finally, if $(t,\xi)\in S_1(F)$, that is $\frac{\partial f}{\partial t}(t,\xi)=0$, we have from (\ref{132}) that \medskip \\
- $\varphi_C(t,\xi)=(1,0)\in \mathbb{R}\times \{0\}=N(F'(t,\xi))$;\medskip\\
- $\psi_C(t,\xi)R(F'(t,\xi))=0$\; i.e.\; $\psi_C(t,\xi) \in R(F'(t,\xi))^\perp$, \medskip \\
since for $(\frac{\partial f}{\partial \xi}(t,\xi)v,v) \in R(F'(t,\xi)), v \in \varXi$, it follows that 
$\psi_C(t,\xi)(\frac{\partial f}{\partial \xi}(t,\xi)v,v)) = \frac{\partial f}{\partial \xi}(t,\xi)v- \frac{\partial f}{\partial \xi}(t,\xi)v=0$. \\
Therefore we have established that (\ref{134}) defines a $C^{d-1}$ f-pair for $F$ on $U$.\\
\indent
It is now easy to determine the fibering functionals related to the canonical pair $(\varphi_C,\psi_C)$: they will be called the \textit{canonical fibering functionals} for the LS-map $F$.\\
From the definition of fibering functionals (see Definition \ref{D116}) we have that\\
$J_0(\varphi_C,\psi_C) \equiv J_{0,C}=\psi_{C}F'\varphi_C$, i.e. $\; \forall \,(t,\xi)\in U$\\
\begin{equation*}
J_{0,C}(t,\xi)=\psi_{C}(t,\xi)F'(t,\xi)\varphi_{C}(t,\xi)=\psi_{C}(t,\xi) \left[ \begin{array}{cc}
\frac{\partial f}{\partial t}(t,\xi)&\frac{\partial f}{\partial\xi}(t,\xi)\\\\
0 &\mathbf{1_\varXi}
\end{array} \right]
\left[ \begin{array}{c}
1 \\
\\
0 \end{array} \right] = 
\end{equation*}
\begin{equation*}
\; = \psi_{C}(t,\xi)(\frac{\partial f}{\partial t}(t,\xi),0)=\frac{\partial f}{\partial t}(t,\xi)-\frac{\partial f}{\partial \xi}(t,\xi)(0)=\frac{\partial f}{\partial t}(t,\xi).
\end{equation*}
\\\\
Moreover $I_{1}(\varphi_C,\psi_C)(t,\xi)\equiv I_{1,C}(t,\xi)=J'_{0,C}(t,\xi)$, i.e. $\forall \;(t,\xi)\in U,\;\forall \;\;(r,v)\in \mathbb{R}\times \varXi$
\begin{equation*}
I_{1,C}(t,\xi)(r,v)=J'_{0,C}(t,\xi)(r,v)=(\frac{\partial f}{\partial t})'(t,\xi)(r,v)=r\frac{\partial^2 f}{\partial  t^2}(t,\xi) + \frac{\partial^2 f}{\partial t\partial \xi}(t,\xi)v.
\end{equation*}
Then $J_{1}(\varphi_C,\psi_C)(t,\xi)\equiv J_{1,C}=I_{1,C}\varphi_C$, i.e. $\;\forall \,(t,\xi)\in U$ 
\begin{equation*}
J_{1,C}(t,\xi)=I_{1,C}(t,\xi)\varphi_{C}(t,\xi)=I_{1,C}(t,\xi)(1,0)=\frac{\partial^2 f}{\partial t^2}(t,\xi).
\end{equation*}
In the same fashion, arguing by induction it is easy to show that the canonical fibering functionals  $\{J_{h}(\varphi_C,\psi_C)\equiv J_{h,C}:h=0,\ldots,d-1\},\{I_{h}(\varphi_C,\psi_C)(t,\xi)\equiv I_{h,C}(t,\xi):h= 1,\ldots,d-1,(t,\xi)\in U\}$ have the form
\begin{align}\label{135}
J_{h,C}(t,\xi)&=\frac{\partial^{h+1} f}{\partial t^{h+1}}(t,\xi),(t,\xi)\in U, \nonumber\\
\\
I_{h,C}(t,\xi)(r,v)&=r \frac{\partial^{h+1} f}{\partial t^{h+1}}(t,\xi)+\frac{\partial^{h+1} f}{\partial t^{h}\partial \xi}(t,\xi) v, (t,\xi) \in U, (r,v) \in \mathbb{R}\times \varXi, \nonumber\\
\text{i.e.} \quad  I_{h,C}(t,\xi)&=( \frac{\partial^{h+1} f}{\partial t^{h+1}}(t,\xi), \frac{\partial^{h+1} f}{\partial t^{h}\partial \xi}(t,\xi))\in \mathbb{R}\times \varXi^\ast \cong (\mathbb{R}\times\varXi)^\ast. \nonumber
\end{align}
\par
\subsubsection{}\label{sss132} As in \ref{sss131}, let $F:U\subseteq \mathbb{R}\times \varXi \rightarrow V\subseteq \mathbb{R}\times\varXi, F(t,\xi)=(f(t,\xi),\xi)$, be a $C^d$  LS-map and let $(\varphi,\psi)$ be an assigned f-pair for $F$ on $U$. Here we are interested in describing the relation between $(\varphi,\psi)$ and any other f-pair $(\widetilde{\varphi},\widetilde{\psi})$ for $F$ on $U$, at least on a neighbourhood of $S_1(F)=\{(t,\xi)\in U:\frac{\partial f}{\partial t}(t,\xi)=0\}$.\\
If $(\varphi,\psi)$ is a $C^{d-1}$ f-pair on $U$ then 
\begin{equation*}
\varphi:U\subseteq \mathbb{R}\times \varXi \rightarrow (\mathbb{R}\times \varXi)\setminus \{(0,0)\}
\end{equation*}
and
\begin{equation*}
\psi:U\subseteq \mathbb{R}\times\varXi \rightarrow (\mathbb{R}\times\varXi)^\ast \setminus\{0\}\equiv (\mathbb{R}\times\varXi^*)\setminus \{(0,0)\}
\end{equation*}
have the form  $\varphi(t,\xi)=(\alpha(t,\xi),a(t,\xi))$ and $\psi(t,\xi)=(\beta(t,\xi),b(t,\xi))$, for suitable maps $\alpha:U\rightarrow \mathbb{R}, a:U\rightarrow \varXi, \beta:U\rightarrow\mathbb{R}, b:U\rightarrow\varXi^\ast$ which are of class $C^{d-1}$.\\
If $(t,\xi)\in S_1(F)$,we have that $\varphi(t,\xi)=(\alpha(t,\xi),a(t,\xi))\in N(F'(t,\xi))=\mathbb{R}\times \{0\}$ and thus $a(t,\xi)=0$ on $S_{1}(F)$. Then $\varphi(t,\xi)\neq 0$ implies that $\alpha(t,\xi)\neq 0$ on $S_{1}(F)$.\\
\indent In the same manner, for $(t,\xi)\in S_{1}(F)$, one has that $\psi(t,\xi)=(\beta(t,\xi),b(t,\xi)) \in R(F'(t,\xi))^\perp$, hence $\psi(t,\xi)\{\frac{\partial f}{\partial \xi}(t,\xi)v,v):v \in \varXi\}=0$. Therefore, $\;\forall \, v \in \varXi$,
\begin{equation*}
0=(\beta(t,\xi),b(t,\xi))({\frac{\partial f}{\partial \xi}}(t,\xi)v,v)=\beta(t,\xi){\frac{\partial f}{\partial \xi}}(t,\xi)v+b(t,\xi)v,
\end{equation*}
i.e. $\beta(t,\xi){\frac{\partial f}{\partial \xi}}(t,\xi)+ b(t,\xi)=0\in \varXi^\ast$. It is important to note that $\beta(t,\xi)$ is not zero, otherwise this would imply $b(t,\xi)=0$ and in turn $\psi(t,\xi)=0$, which would contradict the definition of $\psi$. Hence $\beta(t,\xi)\neq 0$ and $b(t,\xi)=-\beta(t,\xi)\frac{\partial f}{\partial \xi}(t,\xi)$  on $S_{1}(F)$.\\
If $(\widetilde{\varphi},\widetilde{\psi})$ is another $C^{d-1}$ f-pair  for $F$ on $U$ the same argument shows that $\widetilde{\varphi}$ and $\widetilde{\psi}$ have the form $\widetilde{\varphi}(t,\xi)=(\widetilde{\alpha}(t,\xi),\tilde{a}(t,\xi))$ and $\widetilde{\psi}(t,\xi)=(\widetilde{\beta}(t,\xi),\tilde{b}(t,\xi))$, with suitable  $C^{d-1}$ maps $\widetilde{\alpha},\tilde{a},\widetilde{\beta},\tilde{b}$. Moreover $\widetilde{\alpha}(t,\xi)\neq 0$ and $\tilde{a}(t,\xi)=0$ hold on $S_{1}(F)$ and, analogously, $\widetilde{\beta}(t,\xi)\neq 0$ and $\tilde{b}(t,\xi)= -\widetilde{\beta}(t,\xi)\frac{\partial f}{\partial\xi}(t,\xi)$ are true on  $S_1(F)$.\\
From the continuity on $U$ of the maps $\alpha,\widetilde{\alpha},\beta,\widetilde{\beta}$ there exists an open neighbourhood $U'$ of $S_{1}(F)$ such that they are simultaneously different from zero on $U'$. Hence on   $U'$ we have that
\begin{equation*}
\widetilde{\varphi}=(\widetilde{\alpha},\tilde{a}) =\alpha^{-1}\widetilde{\alpha}\cdotp (\alpha, a)+(0,\tilde{a}-\alpha^{-1}\widetilde{\alpha}a)=\gamma\varphi+c,
\end{equation*}
where we write $\alpha^{-1}=\alpha^{-1}(t,\xi)$ for the reciprocal of $\alpha(t,\xi),(t,\xi)\in U'$.\\
Note that $\gamma:=\alpha^{-1}\widetilde{\alpha}:U'\rightarrow \mathbb{R}$ is a non-zero $C^{d-1}$ function on $U'$, while the $C^{d-1}$ map $c:=(0,\tilde{a}-\alpha^{-1}\widetilde{\alpha}a):U'\rightarrow \mathbb{R}\times \varXi$ is equal to zero on $S_{1}(F)$ because the same is true for $a,\tilde{a}$.\\
Analogously, on $U'$ we can write 
\begin{equation*}
\widetilde{\psi}=(\widetilde{\beta},\tilde{b})=\beta^{-1}\widetilde{\beta}\cdotp (\beta,b)+(0,\tilde{b}-\beta^{-1}\widetilde{\beta}b)=\delta\psi+d
\end{equation*}
with  $\delta:=\beta^{-1}\widetilde{\beta}:U'\rightarrow\mathbb{R}$ a non-zero $C^{d-1}$ function on $U'$  and $d:=(0,\tilde{b}-\beta^{-1}\widetilde{\beta}b):U'\rightarrow \mathbb{R}\times \varXi^\ast\equiv( \mathbb{R}\times \varXi)^\ast$ a $C^{d-1}$ map. Once again, the map $d$ is zero on $S_{1}(F)$ because  $b=-\beta\frac{\partial f}{\partial \xi}$ and $\tilde{b}=-\widetilde{\beta}\frac{\partial f}{\partial \xi} $ on $S_{1}(F)$.\medskip \\
The above discussion can be summarized in the following

\begin{proposition}
\label{Pro133} \textit{If} $(\varphi,\psi),(\widetilde{\varphi},\widetilde{\psi})$ \textit{are} $C^{d-1}$ \textit{fibering pairs for} $F$ \textit{on} $U$\textit{ then there exist an open neighbourhood }$U'$ \textit{of} $S_{1}(F)$ \textit{in } $U$ \textit{and } $C^{d-1}$ \textit{maps}
\begin{equation*}
\gamma:U'\rightarrow\mathbb{R},\delta:U'\rightarrow\mathbb{R},c:U'\rightarrow \mathbb{R}\times \varXi,d:U'\rightarrow\mathbb{R}\times \varXi^\ast \cong (\mathbb{R}\times\varXi)^\ast
\end{equation*}
\textit{such that}
\begin{equation*}
\widetilde{\varphi}=\gamma\varphi+c, \widetilde{\psi}=\delta\psi+d \;\; on \;\; U'.
\end{equation*}
\textit{Moreover,} $\gamma\neq 0,\delta \neq 0 $ \textit{on} $U'$ \textit{and} $c=0,d=0$ \textit{on} $S_{1}(F)$.
\end{proposition}
\vspace{4pt}
As the following example shows, the neighbourhood  $U'$ of $S_{1}(F)$ depends on the chosen f-pairs.

\begin{example}
\label{Ex134} Let  $F:\mathbb{R}\times \mathbb{R}\rightarrow \mathbb{R}\times \mathbb{R}$ be given by $F(t,\xi):=(t^2,\xi)$, that is  $f(t,\xi)=t^2$. Then  $S_{1}(F)=\{(t,\xi)\in \mathbb{R}\times \mathbb{R}:\frac{\partial f}{\partial t}(t,\xi)=2t=0\}$, i.e. $S_{1}(F)$ coincides with the  $\xi$-axis. Consider $\varphi(t,\xi):=(\varepsilon-t^2,t)$ for a fixed $\varepsilon>0$. Then  $\varphi\neq 0$ everywhere, and for $(0,\xi)\in S_{1}(F)$ we have that $\varphi(0,\xi)=(\varepsilon,0)$ generates $N(F'(t,\xi))=\mathbb{R}\times \{0\}$. Hence $\varphi$ is a (global) kernel fibering map for $F$. Notice that $\varphi(t,\xi)=(\varepsilon-t^2,t)=(\varepsilon-t^2)(1,0)+(0,t)=(\varepsilon-t^2)\varphi_C+(0,t)$. Consequently, for the kernel fibering maps $\varphi$ and $\varphi_C$ one has $\gamma=\varepsilon-t^2,c=(0, t)$. Therefore, the largest neighbourhood $U'$ of $S_{1}(F)$ where $\gamma$ is different from 0 is the set $\{(t,\xi)\in \mathbb{R}\times\mathbb{R}:t^2<\epsilon\}$.
\end{example}

\subsubsection{} At this stage we have shown the existence of a fibering pair (the canonical one) for the LS-maps and in Proposition \ref{Pro133} we have studied the relation between two different f-pairs for such a class of maps. The f-pairs are also used to define the fibering functionals on which we will then build the classification of the singularities. It is thus worthwhile to give here a first account of the basic relationship between fibering functionals and singular set.\\
\indent
Let $F$ be the LS-map considered in Subsection \ref{sss131}, and let $(\varphi_C,\psi_C)$ be the canonical f-pair associated with $F$. In \ref{sss131} we saw that $S_1(F)=\{(t,\xi)\in U:\frac{\partial f}{\partial t}(t,\xi)=0\}$ and $J_{0,C}(t,\xi)=\frac{\partial f}{\partial t}(t,\xi)$. Hence the following relation holds:
\begin{equation}\label{136}
J_{0,C}(t,\xi)=0 \;\;\Leftrightarrow\;\;(t,\xi)\in S_{1}(F).
\end{equation}
\par
However, in general this is not true for all f-pairs $(\varphi,\psi)$ associated with $F$. In fact, if  $(\varphi,\psi)$ is an f-pair for $F$ and $J_0\equiv J_{0}(\varphi,\psi)$ is the related first fibering functional, i.e. $J_0(t,\xi)=\psi(t,\xi)F'(t,\xi)\varphi(t,\xi)$, then
\begin{equation}\label{137}
S_{1}(F)\subseteq\{(t,\xi)\in \mathbb{R}\times \varXi:J_{0}(\varphi,\psi)(t,\xi)=0\}
\end{equation}
since $F'(t,\xi)\varphi(t,\xi)=0$ on $S_{1}(F)$. On the other hand, if $J_{0}(\varphi,\psi)(t,\xi)=0$  then it is not necessarily true that $(t,\xi)\in S_{1}(F)$ as the following example illustrates.\\

\begin{example}
\label{Ex136} If $\widetilde{F}:\mathbb{R}\times \mathbb{R} \rightarrow \mathbb{R}\times \mathbb{R}$ is given by $\widetilde{F}(t,\xi):=(t^3,\xi)$, that is $\tilde{f}(t,\xi)=t^3$, one has
\begin{equation*}
\widetilde{F}'(t,\xi) = \left[ \begin{array}{cc}
\frac{\partial \tilde{f}}{\partial t}(t,\xi)&\frac{\partial \tilde{f}}{\partial\xi}(t,\xi)\\\\
0 &\mathbf{1_\mathbb{R}}
\end{array} \right] =
\left[ \begin{array}{cc}
3t^2 & 0 \\
\\
0 & 1 \end{array} \right]
\end{equation*}
and so $S_{1}(\widetilde{F})=\{(t,\xi)\in \mathbb{R}\times \mathbb{R}:\frac{\partial \tilde{f}}{\partial t}(t,\xi)=3t^2=0\}$, i.e. $S_{1}(\widetilde{F})$ is the $\xi$-axis. \\
From the characterization (\ref{132}) we know that $N(\widetilde{F}'(t,\xi))=\mathbb{R}\times \{0\}$ and $\mathbb{R}(\widetilde{F}'(t,\xi))=\{0\}\times \mathbb{R}$, for $(t,\xi)\in S_{1}(\widetilde{F})$. It is easy to see that, for all $(t,\xi)\in \mathbb{R}\times \mathbb{R}$, the maps $\widetilde{\varphi}(t,\xi):=(1,t) \in \mathbb{R}\times \varXi$ and $\widetilde{\psi}(t,\xi):=(1,-3t) \in \mathbb{R}\times \varXi^\ast \cong (\mathbb{R}\times \varXi)^\ast$ can  be taken as a fibering pair for $\widetilde{F}$. If $\widetilde{J}_0:= J_{0}(\widetilde{\varphi},\widetilde\psi)$ then for all $(t,\xi)\in \mathbb{R}\times \mathbb{R}$ we obtain 
\begin{equation*}
\widetilde{J}_0(t,\xi)=\widetilde\psi(t,\xi)\widetilde{F}'(t,\xi)\widetilde{\varphi}(t,\xi)=
\left [  \begin{array}{cc}
1,-3t\end{array} \right ]
\left [ \begin{array}{cc}
3t^2&0\\
\\
0 & 1\end{array} \right]
\left[ \begin{array}{c}
1\\
\\
t \end{array} \right]
=3t^2-3t^2=0.
\end{equation*}
\end{example}
\vspace{4pt}
However, there exist maps for which the inclusion in (\ref{137}) is an equality for all f-pairs (at least locally near a given singular point), as shown in the following example.\\
\begin{example}
 Let  $F(t,\xi):=(t^2,\xi)$  be the map of Example \ref{Ex134}. In this case
\begin{equation*}
F'(t,\xi)= \left[ \begin{array}{cc} 
2t & 0\\
\\
0 & 1 \end{array} \right]
\end{equation*}
and we know that $S_{1}(F)$ coincides with the $\xi$-axis.\\
By \ref{sss132}, each fibering pair $(\varphi,\psi)$ has the form  
\begin{equation*}
\varphi(t,\xi) = (\alpha(t,\xi),a(t,\xi))\quad,\quad\psi(t,\xi)=(\beta(t,\xi),b(t,\xi)),
\end{equation*} 
where $\alpha:\mathbb{R}\rightarrow\mathbb{R},a:\mathbb{R}\rightarrow\mathbb{R},\beta:\mathbb{R} \rightarrow\mathbb{R},b:\mathbb{R}\rightarrow\mathbb{R}^\ast\cong\mathbb{R}$ are (smooth) maps such that one has $\alpha \neq 0,a=0,\beta\neq 0$ and  $b=-\beta \frac{\partial f}{\partial \xi}$ on $S_{1}(F)$. Since  $\frac{\partial f}{\partial \xi}\equiv 0$ then $b=0$ on $S_{1}(F)$. Hence
\begin{equation*}
\begin{split}
J_{0}(t,\xi)=\psi(t,\xi)F'(t,\xi)\varphi(t,\xi)&= [\alpha(t,\xi),a(t,\xi)]
\left[ \begin{array}{cc}
2t & 0\\
\\
0 & 1 \end{array} \right]
\left[ \begin{array}{c}
\beta(t,\xi)\\
\\
b(t,\xi)  \end{array} \right] =\\
&=2t\alpha(t,\xi)\beta(t,\xi)+a(t,\xi)b(t,\xi).
\end{split}
\end{equation*}
\indent Now we show that, for any f-pair $(\varphi,\psi)$ for $F$, in a neighbourhood of $S_{1}(F)=\{(t,\xi) \in \mathbb{R}\times\mathbb{R}:t=0\}$ one has that the equality holds in (\ref{137}). In other words we claim that, near the $\xi$-axis,
\begin{equation*}
2t\alpha(t,\xi)\beta(t,\xi)+a(t,\xi)b(t,\xi)=0 \quad \Leftrightarrow \quad t=0.
\end{equation*}
In fact, since  $a(0,\xi)=0$ for all $\xi\in \mathbb{R}$, by the fundamental theorem of calculus we have that
\begin{equation*}
a(t,\xi)=a(t,\xi)-a(0,\xi)=\int_0^1[\frac{d}{ds}a(st,\xi)]ds=\int_0^1 \frac{\partial a}{\partial t}(st,\xi)t\;ds=tc(t,\xi),
\end{equation*}
where $c(t,\xi):=\int_0^1 \frac{\partial a}{\partial t}(st,\xi)ds$ is a smooth function. In the same way $b(t,\xi)=td(t,\xi)$ with $d(t,\xi)$ a suitable smooth  function. It follows that 
\begin{equation*}
J_{0}(t,\xi)=2t\alpha(t,\xi)\beta(t,\xi)+a(t,\xi)b(t,\xi)=t[2\alpha(t,\xi)\beta(t,\xi)+tc(t,\xi)d(t,\xi)].
\end{equation*}
Since $2\alpha(t,\xi)\beta(t,\xi)+tc(t,\xi)d(t,\xi)\neq 0$ on the $\xi$-axis, by continuity we can choose a neighbourhood $U$ of $S_{1}(F)$ such that  $2\alpha(t,\xi)\beta(t,\xi)+tc(t,\xi)d(t,\xi)\neq 0$, for $(t,\xi)\in U$. On such a neighbourhood
\begin{equation*}
J_{0}(t,\xi)=t[2\alpha(t,\xi)\beta(t,\xi)+tc(t,\xi)d(t,\xi)]=0 \quad \Leftrightarrow \quad t=0, 
\end{equation*}
which is indeed the result we wanted to prove.
\end{example}
\vspace{4pt}
As we will see in Chapter \ref{s2}, the different behaviours exhibited by the maps $\tilde{F},F$ cited in the above examples are related to the notion of 1-transversality (see Definition \ref{D211} and Section \ref{ss23}). Precisely, all the singular points of $F$ are 1-transverse while $\tilde{F}$ has no 1-transverse singularities. \vspace{2pt}
\subsection{Fibering Pairs in the General Case}\label{ss14}
\subsubsection{}\label{sss141} In order to study the local behaviour of a map near a singularity it is useful to introduce the \textit{germs} of fibering maps. Let  $F:U\subseteq X \rightarrow V \subseteq Y$  be a  $C^d$, 0-Fredholm map between open subsets $U, V$ of the $B$-spaces  $X, Y$  and let $u \in U$ be a simple singularity for $F$, i.e. dim$\,N(F'(u))=1$. Let $(\varphi,\psi),(\widetilde{\varphi},\widetilde{\psi})$ be two $C^{d-1}$ f-pairs for $F$, defined on open neighbourhoods $U'$ and $U''$ of $u$ respectively. We say that $(\varphi,\psi),    (\widetilde{\varphi},\widetilde{\psi})$ are \textit{equivalent} if there exists a neighbourhood $V$ of $u$ such that $V \subseteq U' \cap U''$ and $(\varphi,\psi)=(\widetilde{\varphi},\widetilde{\psi})$ on $V$.\\
\indent
This is indeed an equivalence relation on the set of all $C^{d-1}$ f-pairs for $F$ which are defined on a neighbourhood of $u$. An equivalence class $[(\varphi,\psi)]$ is called a \textit{germ of the fibering pair $(\varphi,\psi)$ at u}. Hence this class consists of all f-pairs for $F$ which are defined and coincide near $u$. We denote by $\mathcal{P}(F,u)$ the set of all germs of $C^{d-1}$ f-pairs for $F$ at $u$; however, we will usually omit the brackets when discussing germs of f-pairs near $u$ and use the same symbols as for f-pairs.\\
\indent Of course, f-pairs in the same germ give raise to fibering functionals which are equivalent, i.e. they coincide near $u$, and thus \textit{germs of the fibering functionals at $u$} are also well-defined.\\
\indent In this section we prove that $\mathcal{P}(F,u)$ is not empty, i.e. a fibering pair exists near a simple singularity (see Theorem \ref{Teo146}).
\subsubsection{}\label{sss142} Let us study how f-pairs and fibering functionals are modified by changes of coordinates: this allows us to see how local diffeomorphisms act on the set $\mathcal{P}(F,u)$.\\
\indent Let us first consider the $C^d$ commutative diagram
\medskip
\begin{equation}\label{141}
\begin{tabular}{ccc}
&$F \quad $&  \medskip \\
$U \subseteq X \quad $&$ \rightarrow \quad $&$ V \subseteq Y $ \medskip  \\
$\gamma \downarrow $&&$ \downarrow \delta$ \medskip \\ 
$\widetilde{U}\subseteq \widetilde{X} \quad $&$ \rightarrow  \quad $&$ \widetilde{V}\subseteq\widetilde{Y}$,  \medskip \\
&$\widetilde{F}\quad $&\\
\end{tabular}
\end{equation} \medskip
where $F:U\subseteq X\rightarrow V\subseteq Y$ and $\widetilde{F}:\widetilde{U}\subseteq \widetilde{X}  \rightarrow \widetilde{V}\subseteq \widetilde{Y}$ are $C^d$ maps, with $U,V,\widetilde{U}$ and $\widetilde{V}$ open subsets in the $B$-spaces $X,Y,\widetilde{X}$ and $\widetilde{Y}$ respectively, and where $\gamma, \delta$ are $C^d$ diffeomorphisms such that $\gamma(U)=\widetilde{U}$ and $\delta(V)=\widetilde{V}$. As usual, we assume that F is a 0-Fredholm map with dim$\,N(F'(u))\leq 1,\,\forall\,u\in U$. From (\ref{141}) we have that
\begin{equation*}
\widetilde{F}(\tilde{u})=\delta(F(\gamma^{-1}(\tilde{u}))), \,\forall\,\tilde{u} \in \widetilde{U}.
\end{equation*}
Since $\gamma\hspace{1pt}'(u)\in GL(X,\widetilde{X})$ and $\delta'(h)\in GL(Y,\widetilde{Y}), \,\forall \,u\in U$ and  $\forall \,h\in V$, differentiation of the above relation gives
\begin{align}\label{142}
\begin{split}
\widetilde{F}'(\tilde{u})&=\delta'(F(\gamma^{-1}(\tilde{u})))F'(\gamma^{-1}(\tilde{u}))(\gamma^{-1})'(\tilde{u})=\\
&=\delta'(F(\gamma^{-1}(\tilde{u})))F'(\gamma^{-1}(\tilde{u}))(\gamma\,'(\gamma^{-1}(\tilde{u})))^{-1}.
\end{split}
\end{align}
It is then easy to deduce,$\, \forall \,\tilde{u}\in \widetilde{U}$, the following relations between kernels and ranges of the Fréchet derivatives of $F$ and $\widetilde{F}$
\begin{eqnarray}\label{143}
\begin{split}
N(\widetilde{F}'(\tilde{u})&=\gamma\hspace{1pt}'(\gamma^{-1}(\tilde{u}))N(F'(\gamma^{-1}(\tilde{u})))\\
R(\widetilde{F}'(\tilde{u}))&=\delta'(F(\gamma^{-1}(\tilde{u})))R(F'(\gamma^{-1}(\tilde{u}))), 
\end{split}
\end{eqnarray}
hence we obtain
\begin{eqnarray*}
\begin{split}
\text{dim}\,N(\widetilde{F}'(\tilde{u}))&=\text{dim}\,N(F'(\gamma^{-1}((\tilde{u}))),\\
\text{codim}\,R(\widetilde{F}'(\tilde{u}))&=\text{codim}\,R(F'(\gamma^{-1}((\tilde{u}))).
\end{split}
\end{eqnarray*}
Since $F$ is a 0-Fredholm map on $U$, from the previous relations it follows that $\widetilde{F}$ is also a 0-Fredholm map on $\widetilde{U}$ with dim$\,N(\widetilde{F}'(\tilde{u}))\leq 1,\,\forall \, \tilde{u}\in \widetilde{U}$, and dim$\,N(\widetilde{F}'(\tilde{u}))=1$ \textit{iff} dim$ \,N(F'(\gamma^{-1}(\tilde{u})))=1$, i.e.
\begin{equation}\label{144}
S_{1}(\widetilde{F})=\gamma(S_{1}(F)).
\end{equation}
\indent We are now able to show that, given a $C^{d-1}$ f-pair $(\varphi,\psi)$ for $F$ on $U$, there is a natural way to define a $C^{d-1}$ f-pair $(\widetilde{\varphi},\widetilde{\psi})$ for $\widetilde{F}$ on $\widetilde{U}$ related to $(\varphi,\psi)$. It is sufficient to define, $\forall \,\tilde{u}\in \widetilde{U}$, the maps
\begin{equation}\label{145}
\begin{split}
\widetilde{\varphi}(\tilde{u}):&=\gamma\hspace{1pt}'(\gamma^{-1}(\tilde{u}))\varphi(\gamma^{-1}(\tilde{u}))=((\gamma^{-1})'(\tilde{u}))^{-1}\varphi(\gamma^{-1}(\tilde{u})),\\
\widetilde{\psi}(\tilde{u}):&=\psi(\gamma^{-1}(\tilde{u}))(\delta'(F(\gamma^{-1}(\tilde{u}))))^{-1}=\\
&=\psi(\gamma^{-1}(\tilde{u}))(\delta^{-1})'(\delta(F(\gamma^{-1}(\tilde{u}))))=\psi(\gamma^{-1}(\tilde{u}))(\delta^{-1})'(\widetilde{F}(\tilde{u})).
\end{split}
\end{equation}
Since $\varphi\in C^{d-1}(U,X \setminus \{0\}),\psi\in C^{d-1}(U,Y^\ast \setminus \{0\})$ and $\varphi(u)\in N(F'(u)),\psi(u)\in R(F'(u))^{\perp},\,\forall \,u \in S_1(F)$ it is trivial to see that $\widetilde{\varphi}(\tilde{u})\in \widetilde{X}$ and $\widetilde{\psi}(\tilde{u})\in \widetilde{Y}^\ast$. Moreover, since $\gamma\hspace{1pt}'$ and $\delta'$ are maps of class $C^{d-1}$ then $\widetilde{\varphi}\in C^{d-1}(\widetilde{U},\widetilde{X}), \widetilde{\psi} \in C^{d-1}(\widetilde{U},\widetilde{Y}^\ast)$ (because the map $A\in GL(M, N)\mapsto A^{-1}\in GL(N, M)$ is real analytic if $M, N$ are $B$-spaces, see \cite{Ze1}, exercise 8.21, chapter 8). Finally, as $\gamma\hspace{1pt}'(\cdotp)$ and $\delta'(\cdotp)$ are isomorphisms, then $\widetilde{\varphi}(\tilde{u})$ and $\widetilde{\psi}(\tilde{u})$ are always different from zero; moreover, (\ref{143}) and (\ref{145}) imply that $\widetilde{\varphi}(\tilde{u})\in N(\widetilde{F}'(\tilde{u}))$ and $\widetilde{\psi}(\tilde{u})\in R(\widetilde{F}'(\tilde{u}))^{\perp},\,\forall \,\tilde{u} \in S_1(\widetilde{F})$. Hence $(\widetilde{\varphi},\widetilde{\psi})$ is an f-pair for $\widetilde{F}$ on $\widetilde{U}$. Note that, in Vector Field Theory, the map $\widetilde{\varphi}$ is nothing but the \textit{``push-forward} of $\varphi$ by $\gamma$'' (cf. \cite{A-M-R}, Definitions 4.2.1 (ii)).\\
\indent
Let us show how fibering functionals derived from $(\widetilde{\varphi},\widetilde{\psi})$ are related to the ones obtained by $(\varphi,\psi)$. By Definition \ref{D116} and (\ref{142}) we get\\
\begin{align*}
J_0(\widetilde{\varphi},\widetilde{\psi})(\tilde{u})&=\widetilde{\psi}(\tilde{u})\widetilde{F}'(\tilde{u})\widetilde{\varphi}(\tilde{u})=  \\ 
&=\psi(\gamma^{-1}(\tilde{u}))(\delta'(F(\gamma^{-1}(\tilde{u}))))^{-1} \delta'(F(\gamma^{-1}(\tilde{u})))\,\circ\\
& \quad \circ F'(\gamma^{-1}(\tilde{u}))(\gamma\hspace{1pt}'(\gamma^{-1}(\tilde{u})))^{-1}\gamma\hspace{1pt}'(\gamma^{-1}(\tilde{u}))\varphi(\gamma^{-1}(\tilde{u}))= \\
&=\psi(\gamma^{-1}(\tilde{u}))F'(\gamma^{-1}(\tilde{u}))\varphi(\gamma^{-1}(\tilde{u}))=J_{0}(\varphi,\psi) (\gamma^{-1}(\tilde{u}));  \\  \\
I_{1}(\widetilde{\varphi},\widetilde{\psi})(\tilde{u})&=J_{0}(\widetilde{\varphi},\widetilde{\psi})'(\tilde{u})=\\
&=J_{0}(\varphi,\psi)'(\gamma^{-1}(\tilde{u}))(\gamma^{-1})'(\tilde{u})=I_{1}(\varphi,\psi)(\gamma^{-1}(\tilde{u}))(\gamma^{-1})'(\tilde{u});  \\ \\
J_{1}(\widetilde{\varphi},\widetilde{\psi})(\tilde{u})&=I_{1}(\widetilde{\varphi},\widetilde{\psi})(\tilde{u})\widetilde{\varphi}(\tilde{u})=\\
&=I_{1}(\varphi,\psi)(\gamma^{-1}(\tilde{u}))(\gamma^{-1})'(\tilde{u})((\gamma^{-1})'(\tilde{u}))^{-1}\varphi(\gamma^{-1}(\tilde{u}))=  \\
&=I_{1}(\varphi,\psi)(\gamma^{-1}(\tilde{u}))\varphi(\gamma^{-1}(\tilde{u}))=J_{1}(\varphi,\psi)(\gamma^{-1}(\tilde{u})).
\end{align*}
By induction, $\;\forall \, \tilde{u} \in \widetilde{U}$ we obtain 
\begin{equation}\label{146}
\begin{split}
J_{h}(\widetilde{\varphi},\widetilde{\psi})(\tilde{u})&=J_h(\varphi,\psi)(\gamma^{-1}(\tilde{u})), \qquad 0\leq h \leq d-1;\\
I_h(\widetilde{\varphi},\widetilde{\psi})(\tilde{u})&= I_h(\varphi,\psi)(\gamma^{-1}(\tilde{u}))(\gamma^{-1})'(\tilde{u}), \qquad 1\leq h \leq d-1.
\end{split}
\end{equation}
\indent
It is now convenient to give a name to the fibering pair $(\widetilde{\varphi},\widetilde{\psi})$ obtained from $(\varphi,\psi)$. Since $\widetilde{\varphi}$ and $\widetilde{\psi}$ are obtained from $\varphi$ and $\psi$ by using the diffeomorphisms $\gamma$ and $\delta$ it is useful to rewrite (\ref{145}) as $\widetilde{\varphi}=\mathcal{T}_{N} [\gamma,\delta] \varphi$, $\widetilde{\psi}=\mathcal{T}_{R}[\gamma,\delta] \psi$ where,$ \; \forall \, \tilde{u} \in \widetilde{U}$,
\begin{equation}\label{147}
\begin{split}
(\mathcal{T}_{N} [\gamma, \delta] \varphi)(\tilde{u}) &= \gamma\hspace{1pt}'(\gamma^{-1}(\tilde{u}))\varphi(\gamma^{-1}(\tilde{u})),\\
(\mathcal{T}_{R} [\gamma, \delta]\psi)(\tilde{u})&=\psi (\gamma^{-1}(\tilde{u})(\delta^{-1})'(\widetilde{F}(\tilde{u})).
\end{split}
\end{equation}
Given the commutative diagram (\ref{141}), we denote by $\mathcal{P}(F,U)$ the set \{$C^{d-1}$ f-pairs for $F$ on $U$\} and, analogously, we have that $\mathcal{P}(\widetilde{F},\widetilde{U}):=$\{$C^{d-1}$ f-pairs for $\widetilde{F}$ on $\widetilde{U}$\}. Then we can introduce the following\\

\begin{definition}
\label{D143} The \textit{pair-transform} of the f-pair $(\varphi,\psi)\in \mathcal{P}(F,U)$ by the diffeomorphisms $\gamma,\delta$ is the f-pair $\mathcal{T}[\gamma,\delta](\varphi,\psi) \in \mathcal{P}(\widetilde{F},\widetilde{U})$ defined by $\mathcal{T}[\gamma,\delta](\varphi,\psi):=(\mathcal{T}_{N} [\gamma, \delta] \varphi,\mathcal{T}_{R} [\gamma, \delta]\psi)$. The same term, \textit{pair-transform}, will be used for the mapping $\mathcal{T}[\gamma,\delta]:\mathcal{P}(F,U)\rightarrow \mathcal{P}(\widetilde{F},\widetilde{U})$ which takes an f-pair $(\varphi,\psi)$ to $\mathcal{T}[\gamma,\delta](\varphi,\psi).$
\end{definition} 
\vspace{4pt}
With the notation given in the above definition we can rewrite the formulas (\ref{146}) as
\begin{equation}\label{148}
\begin{split}
J_{h}\big(\mathcal{T}[\gamma,\delta](\varphi,\psi)\big)(\tilde{u})&=J_h(\varphi,\psi)(\gamma^{-1}(\tilde{u})), \quad 0\leq h \leq d-1; \bigskip  \\
I_{h}\big(\mathcal{T}[\gamma,\delta](\varphi,\psi)\big)(\tilde{u})&= I_{h}(\varphi,\psi)(\gamma^{-1}(\tilde{u}))(\gamma^{-1})'(\tilde{u}), \quad 1\leq h \leq d-1.
\end{split}
\end{equation}
\indent
A basic property of the pair-transform  $\mathcal{T}[\gamma,\delta]$ is established in the following \\
\begin{proposition}
\label{Pro144} \textit{Under the assumptions of the commutative diagram} (\ref{141}), \textit{the pair-transform} $\mathcal{T}[\gamma,\delta]:\mathcal{P}(F,U)\rightarrow \mathcal{P}(\widetilde{F},\widetilde{U})$ \textit{is a bijection.}
\end{proposition}
\par
\textbf{Proof.} We shall prove that the pair-transform $\mathcal{T}[\gamma^{-1},\delta^{-1}]:\mathcal{P}(\widetilde{F},\widetilde{U})\rightarrow \mathcal{P}(F,U) $, where $\gamma^{-1},\delta^{-1}$ refer to the ``inverted'' diagram (\ref{141}) (i.e. the diagram obtained from (\ref{141}) by inverting the vertical arrows), is the inverse of the pair-transform $\mathcal{T}[\gamma,\delta]$. This means that we have to show that 
\begin{align*}
\mathcal{T}[\gamma^{-1},\delta^{-1}](\mathcal{T}[\gamma,\delta](\varphi,\psi)) &= (\varphi,\psi), \;\, \forall \, (\varphi,\psi)\in \mathcal{P}(F,U),\\
\mathcal{T}[\gamma,\delta](\mathcal{T}[\gamma^{-1},\delta^{-1}](\widetilde{\varphi},\widetilde{\psi})) &= (\widetilde{\varphi},\widetilde{\psi}), \;\, \forall \,(\widetilde{\varphi},\widetilde{\psi})\in \mathcal{P}(\widetilde{F},\widetilde{U}). 
\end{align*}
By interchanging the roles of $F,\gamma,\delta$ and $\widetilde{F},\gamma^{-1},\delta^{-1}$ we have that it suffices to prove the first identity. Hence we have to prove that, for any $(\varphi,\psi)\in \mathcal{P}(F,U)$, one has
\begin{equation*}
\mathcal{T}_{N}[\gamma^{-1},\delta^{-1}](\mathcal{T}_{N}[\gamma,\delta]\varphi) = \varphi \quad,\quad \mathcal{T}_{R}[\gamma^{-1},\delta^{-1}](\mathcal{T}_{R}[\gamma,\delta]\psi) = \psi.
\end{equation*}
\indent For $u \in U$, by applying (\ref{147}) twice we obtain
\begin{align*}
\big(\mathcal{T}_{N}[\gamma^{-1},\delta^{-1}](\mathcal{T}_{N}[\gamma,\delta]\varphi)\big)(u)=& (\gamma^{-1})'(\gamma(u))(\mathcal{T}_{N}[\gamma,\delta]\varphi)(\gamma(u)) =\\
=(\gamma^{-1})'(\gamma(u))\gamma\hspace{1pt}'(\gamma^{-1}(\gamma(u))\varphi(\gamma^{-1}&(\gamma(u))=(\gamma\hspace{1pt}'(u))^{-1}\gamma\hspace{1pt}'(u)\varphi(u)=\varphi(u);
\end{align*}
in the same manner
\begin{align*}
\big(\mathcal{T}_{R}[\gamma^{-1},\delta^{-1}](\mathcal{T}_{R}[\gamma,\delta]\psi)\big)(u)&=(\mathcal{T}_{R}[\gamma,\delta]\psi)(\gamma(u))\delta'(F(u))=\\
=\psi(\gamma^{-1}(\gamma(u))(\delta^{-1})'(\widetilde{F}(\gamma(u)))\delta'(F(u)) &= \psi(u)(\delta'(\delta^{-1}(\widetilde{F}(\gamma(u)))))^{-1}\delta'(F(u))\\
=\psi(u)(\delta'(F(u)))^{-1}&\delta'(F(u)) =\psi(u) \,.
\end{align*}
This concludes the proof. $\bracevert$ \\
\par
It is now easy to reformulate the above result in terms of germs. For this purpose, let us suppose that $\gamma$ and $\delta$ in (\ref{141}) are only local diffeomorphisms near $u_o \in U$ and $F(u_o) \in V$ respectively. From now on we shall assume that $u_o$ is a simple singularity for $F$. Thanks to (\ref{144}) we also have that $\tilde{u}_o:=\gamma(u_o)\in \widetilde{U}$ is a simple singularity for $\widetilde{F}$. Let $\mathcal{P}(F,u_o),\mathcal{P}(\widetilde{F},\tilde{u}_o)$ be the collections of germs of $C^{d-1}$ f-pairs near $u_o$ and $\tilde{u}_o$ respectively. It is possible to consider the pair-transform $\mathcal{T}[\gamma,\delta]$ of Definition 1.4.3 as a mapping between germs of fibering maps, that is formulas (\ref{147}) induce the transform
\begin{equation}\label{149}
\mathcal{T}[\gamma,\delta]:\mathcal{P}(F,u_o) \rightarrow \mathcal{P}(\widetilde{F},\tilde{u}_o) 
\end{equation}
where, for the sake of simplicity, we use the same symbol as in  Definition \ref{D143}.
From Proposition \ref{Pro144} we thus get the following\\

\begin{corollary}
\label{Co145} \textit{Let us suppose that there exists a $C^{d}$ l.c.d.}
\\
\begin{equation*}
\begin{tabular}{ccc}
&$F \quad $&  \medskip \\
$u_o \in U \subseteq X \quad $&$ \rightarrow \quad $&$ V \subseteq Y $ \bigskip  \\
$\gamma \downarrow $&&$ \downarrow \delta$ \medskip \\ 
$\tilde{u}_o \in \widetilde{U} \subseteq \widetilde{X} \quad $&$ \rightarrow  \quad $&$ \widetilde{V} \subseteq \widetilde{Y}$ ,\medskip\\
&$ \widetilde{F} \quad $&  \bigskip \\
\end{tabular}
\end{equation*}
\textit{where $F$ is a {\rm 0-}Fredholm  map. Then $\mathcal{T}[\gamma,\delta]:\mathcal{P}(F,u_o) \rightarrow \mathcal{P}(\widetilde{F},\tilde{u}_o)$  is a bijection.}
\end{corollary}
\par
The above result provides the last tool we need to prove that $\mathcal{P}(F,u_o)$ is not empty when $u_o$  is a simple singularity for $F$.\\
\par
Let  $M,N$ be $B$-spaces and $u \in M$: from now on we shall denote by  $C_{u}^{d}(M,N)$ the germs of $C^d$ maps defined from a neighbourhood of $u$ into $N$, cf. Subsection \ref{sss141}. Once again, we adopt the same notation for germs and $C^d$  maps defined near $u$. \\

\begin{theorem}
\label{Teo146} \textit{Let $U,V$ be open subsets in the $B$-spaces X,Y and let $F:U \subseteq X \rightarrow V  \subseteq Y$  be a $C^d$ \text{0-}Fredholm map. If $u_o \in S_1(F)$ then}:
\begin{enumerate}
\item[1)]
$\mathcal{P}(F,u_o)$ \textit{is not empty, i.e. near} $u_o$ \textit{there exists  a} $C^{d-1}$ \textit{fibering pair } $(\varphi,\psi)$ \textit{for }$F$;
\item[2)]
\textit{for any }$(\varphi,\psi)\in \mathcal{P}(F,u_o)$\textit{ and any pair of maps} $\widetilde{\varphi} \in C_{u_o}^{d-1}(X,X), \widetilde{\psi} \in C_{u_o}^{d-1}(X,Y^\ast)$ \textit{one has that }$(\widetilde{\varphi},\widetilde{\psi}) \in \mathcal{P}(F,u_o)$\textit{ if and only if} $\widetilde{\varphi},\widetilde{\psi}$ \textit{have the form}
\begin{center}
$\widetilde{\varphi}=\alpha\varphi + a \;\;,\;\;\widetilde{\psi}=\beta\psi + b \;\;$,
\end{center}
\textit{with} $\alpha,\beta \in C_{u_o}^{d-1}(X,\mathbb{R}), a \in C_{u_o}^{d-1}(X,X), b \in C_{u_o}^{d-1}(X,Y^\ast)$ \textit{such that} $\alpha,\beta$ \textit{are non-vanishing functions and }$a,b$ \textit{vanish on} $S_1(F)$;
\item[3)]
\textit{there exists }$(\varphi,\psi) \in  \mathcal{P}(F,u_o)$ \textit{such that, near }$u_0$,
\begin{center}
$u \in S_1(F) \Leftrightarrow J_{0}(\varphi,\psi)(u) = 0$.
\end{center}
\textit{In general, for }$(\widetilde{\varphi},\widetilde{\psi}) \in  \mathcal{P}(F,u_o)$,\textit{ one has } 
\begin{center}
$u \in S_{1}(F) \Rightarrow J_{0}(\widetilde{\varphi},\widetilde{\psi})(u) = 0$.
\end{center}
\end{enumerate}
\end{theorem}
\begin{remark}
The maps $\alpha,a,\beta,b$ mentioned in the above result are not unique, in general. To show this, let us consider once more the map $F$ of Example \ref{Ex134}, i.e. $F:\mathbb{R}\times \mathbb{R}\rightarrow \mathbb{R}\times \mathbb{R}, F(t,\xi)=(f(t,\xi),\xi):=(t^{2},\xi)$. We saw that $S_{1}(F)=\{(t,\xi)\in \mathbb{R}\times \mathbb{R}:t=0\}$, i.e. $S_{1}(F)$ is the $\xi$-axis. If we consider the canonical fibering map $\varphi_C = (1,0)$ on $\mathbb{R}\times \mathbb{R}$ then $\varphi_C = (1,0) = 1\cdotp \varphi_C + (0,0)$, i.e. $\alpha=1$ and $a=0$. Yet we also have that $\varphi_C=(1,0)=(1+t^{2},0)+(-t^{2},0)=(1+t^{2})\varphi_C+(-t^{2},0)$, i.e. $\alpha=1+t^2$ and $a=(-t^{2},0)$. Note that $\alpha \neq 0$ on $\mathbb{R}\times \mathbb{R}$ and $a = 0$ on $S_{1}(F)$. Of course, an analogous argument holds for $\psi_C$  and, actually, for each f-pair for $F$.\\
\indent On the other hand, we have that $1 = 1 + t^2$ on $S_{1}(F)$. This follows from a general fact: the maps $\alpha,a,\beta,b$ in Theorem 1.4.6 are uniquely determined on $S_{1}(F)$. Indeed, it must be $a(u)=0$ for $u\in S_{1}(F)$, hence $\widetilde{\varphi}(u)=\alpha(u)\varphi(u)$. Since $\widetilde{\varphi}(u)$ and $\varphi(u)$ span the one-dimensional space $N(F'(u))$ it follows that $\alpha(u)$ is uniquely determined for $u \in S_{1}(F)$. The same holds for $\beta$ and $b$.
\end{remark}
\begin{remark}
It may be interesting to see how Theorem \ref{Teo146} applies to the f-pairs (and their related functionals) introduced in Examples \ref{Ex115} and \ref{Ex118}. For the map $F$ associated to problem (P$_{4}$) in \ref{Ex115} we gave three different (global) f-pairs for $F$ and thus point 1) of Theorem \ref{Teo146} is (globally) verified. As for point 2) it is easy to see that, with respect to the f-pair $(\varphi,\psi)$, the other f-pairs $(\widetilde{\varphi},\widetilde{\psi}),(\overline{\varphi},\overline{\psi})$ are obtained as
\begin{center}
$\widetilde{\varphi}(u)=\widetilde{\alpha}(u)\varphi(u)\;,\;\widetilde{\psi}(u)=\widetilde{\beta}(u)\psi(u)$ 
\end{center}
with $\widetilde{\alpha}(u)=1,\;\widetilde{\beta}(u)=\text{exp}[t\cdotp \int_{ 0}^{1}ag'(u)]$ and 
\begin{center}
$\overline{\varphi}(u)=\overline{\alpha}(u)\varphi(u)\;,\;\overline{\psi}(u)=\overline{\beta}(u)\psi(u)$ 
\end{center}
with $\overline{\alpha}(u)=\{1-t+t\cdotp\text{exp}[\int_{0}^{1}ag'(u)]\}\text{exp}[-t\cdotp \int_{ 0}^{t}ag'(u)],\,\overline{\beta}(u)=1$. We notice that $\widetilde{\beta}$ and $\widetilde{\alpha}$ are globally different from zero because they are strictly positive. Finally, for each of the three f-pairs, the associated first functional $J_{0}$ is such that $u\in S_{1}(F)  \Leftrightarrow J_{0}(u)=0$, as seen in Example \ref{Ex118}. Such a good choice is not always possible for an arbitrary f-pair since point 3) states that, in general, the zero-set of the functional $J_{0}$ does not coincide with the singular set $S_{1}(F)$. In the next chapter, in order to get (at least locally) such a useful identification of subsets for any f-pair we will be concerned with the geometrical condition of 1-\textit{transversality} (see Theorem \ref{Teo233}). 
\end{remark}
\par
\textbf{Proof.} 1) Since $u_o \in S_{1}(F)$, by the Local Representation Theorem (see Theorem \ref{Teo122} and Remark \ref{Rem123}) we get the following  $C^d$ l.c.d.
\medskip
\begin{center}
\begin{tabular}{ccc}
&$F \quad $&  \medskip \\
$u_o \in U \subseteq X \quad $&$ \rightarrow \quad $&$ V \subseteq Y $ \medskip  \\
$\gamma \downarrow $&&$ \downarrow \delta$ \medskip \\ 
$(0,0)\in \mathbb{R} \times Z \quad $&$ \rightarrow  \quad $&$ \mathbb{R} \times Z,$  \medskip \\
&$\varPhi \quad $&\\
\end{tabular}
\end{center} \medskip
where $Z=R(F'(u_o))$ and $\varPhi$ is a $C^d$ map defined near $(0,0)$ such that  $\varPhi(t,z)=(f(t,z),z)$, for $(t,z)$ near $(0,0)=\gamma(u_o)$ and with  $f$ a $C^d$ function defined near $(0,0)$.\\
The $C^d$  map $\varPhi:\mathbb{R}\times Z\rightarrow \mathbb{R}\times Z$ is an LS-map. Therefore, near $(0,0)\in \mathbb{R}\times Z$ there exists the $C^{d-1}$ \textit{canonical f-pair} $(\varphi_{C},\psi_{C})$ for $\varPhi$  defined by formulas (\ref{134}). This implies that $(\varphi_{C},\psi_{C}) \in \mathcal{P}(\varPhi,(0,0))$. Furthermore, by Corollary \ref{Co145}, the pair-transform $\mathcal{T}[\gamma^{-1},\delta^{-1}]:\mathcal{P}(\varPhi,(0,0))\rightarrow \mathcal{P}(F,u_o)$ is a bijection. This shows that $\mathcal{P}(F,u_o)$ is not empty since it contains the f-pair 
\begin{equation*}
(\varphi,\psi):=\mathcal{T}[\gamma^{-1},\delta^{-1}](\varphi_{C},\psi_{C}) .
\end{equation*}
\quad 2) We start by proving the \textquotedblleft only if \textquotedblright statement. Since the pair-transform $\mathcal{T}[\gamma^{-1},\delta^{-1}]$ is a bijection for each $(\varphi,\psi), (\widetilde{\varphi},\widetilde{\psi})\in \mathcal{P}(F,u_o)$ there exist $(\widehat{\varphi},\widehat{\psi}),(\overline{\varphi},\overline{\psi})\in \mathcal{P}(\varPhi,(0,0))$ such that 
\begin{center}
$(\varphi,\psi)=\mathcal{T}[\gamma^{-1},\delta^{-1}](\widehat{\varphi},\widehat{\psi})\;,\;(\widetilde{\varphi},\widetilde{\psi})=\mathcal{T}[\gamma^{-1},\delta^{-1}](\overline{\varphi},\overline{\psi}).$\end{center} 
We studied the relation between two different f-pairs for a given LS-map in Proposition 1.3.3; for $(\widehat{\varphi},\widehat{\psi}),(\overline{\varphi},\overline{\psi})\in \mathcal{P}(\varPhi,(0,0))$ there exist maps $\overline{\alpha},\overline{\beta}\in C_{(0,0)}^{d-1} (\mathbb{R} \times Z, \mathbb{R}),\overline{a} \in C_{(0,0)}^{d-1} (\mathbb{R} \times Z, \mathbb{R}\times Z)$ and $\overline{b} \in C_{(0,0)}^{d-1} (\mathbb{R} \times Z, \mathbb{R}\times Z)^\ast)$ such that
\begin{equation}\label{1410}
\overline{\varphi}=\overline{\alpha}\widehat{\varphi}+\overline{a} \,, \quad \overline{\psi}=\overline{\beta}\widehat{\psi}+\overline{b}\, , \quad \text{near \; }(0,0).
\end{equation}
Moreover, $\overline{\alpha}\neq 0, \overline{\beta} \neq 0$ near $S_{1}(\varPhi)$ and $\overline{a}=0,\overline{b}=0$ on 
$S_{1}(\varPhi)$.\\
Then by  formulas (\ref{147}) and (\ref{1410}), for $u$ near $u_o$ we get
\begin{align*}
\widetilde{\varphi}(u)&=\big(\mathcal{T}_{N}[\gamma^{-1},\delta^{-1}]\overline{\varphi}\big)(u)=(\gamma^{-1})'(\gamma(u))\overline{\varphi}(\gamma(u))=\\
&=(\gamma^{-1})'(\gamma(u))[\overline{\alpha}(\gamma(u))\widehat{\varphi}(\gamma(u))+\overline{a}(\gamma(u))]=\overline{\alpha}(\gamma(u))(\gamma^{-1})'(\gamma(u))\widehat{\varphi}(\gamma(u))+\\
&+(\gamma^{-1})'(\gamma(u))\overline{a}(\gamma(u))=\overline{\alpha}(\gamma(u))\big(\mathcal{T}_{N}[\gamma^{-1},\delta^{-1}]\widehat{\varphi}\big)(u)+ (\gamma^{-1})'(\gamma(u))\overline{a}(\gamma(u))=\\
&=\overline{\alpha}(\gamma(u))\varphi(u) + (\gamma^{-1})'(\gamma(u))\overline{a}(\gamma(u));\\
\\
\widetilde{\psi}(u)&=\big(\mathcal{T}_{R}[\gamma^{-1},\delta^{-1}]\overline{\psi}\big)(u)=\overline{\psi}(\gamma(u))\delta'(F(u))=\\
&=[\overline{\beta}(\gamma(u))\widehat{\psi}(\gamma(u))+\overline{b}(\gamma(u))]\delta'(F(u))=\overline{\beta}(\gamma(u))\widehat{\psi}(\gamma(u))\delta'(F(u))+\\
& \quad +\overline{b}(\gamma(u))\delta'(F(u))=\overline{\beta}(\gamma(u))\big(\mathcal{T}_{R}[\gamma^{-1},\delta^{-1}]\widehat{\psi}\big)(u)+\overline{b}(\gamma(u))\delta'(F(u))=\\
&=\overline{\beta}(\gamma(u))\psi(u) + \overline{b}(\gamma(u))\delta'(F(u)).
\end{align*}
Let us set, for each  $u$  near  $u_o$,
\begin{align*}
\alpha(u):=&\overline{\alpha}(\gamma(u)),\quad a(u):=(\gamma^{-1})'(\gamma(u))\overline{a}(\gamma(u)),\\
\beta(u):=&\overline{\beta}(\gamma(u)),\quad b(u):=\overline{b}(\gamma(u))\delta'(F(u)) ,\
\end{align*}
and let us recall that, by (\ref{144}), $u \in S_{1}(F) \Leftrightarrow \gamma(u) \in S_{1}(\varPhi), u \,\:\text{near} \,\: u_o.$\\
Because $\overline{\alpha}\neq 0, \overline{\beta}\neq 0$ at $(0,0)\in S_{1}(\varPhi)$  then we can suppose, by continuity, that $\alpha,\beta$ are non-zero near $u_o$. Finally, since $\overline{a}=0, \overline{b}=0$ on $S_{1}(\varPhi)$ near $(0,0)$ it follows that $a=0, b=0$ on $S_{1}(F)$ near $u_o$. Hence $\widetilde{\varphi}=\alpha\varphi+a, \widetilde{\psi}=\beta\psi+b$ and  $\alpha,\beta,a,b$ have the required properties.\\
\indent
For the \textquotedblleft if \textquotedblright part of 2), let us define  $\widetilde{\varphi}:=\alpha\varphi+a,\widetilde{\psi}:=\beta\psi+b$, with $\alpha,\beta,a,b$ satisfying the properties in the statement. Since $\widetilde{\varphi}(u_o)=\alpha(u_o)\varphi(u_o)+a(u_o)=\alpha(u_o) \varphi(u_o)\neq 0$ then, by continuity, $\widetilde{\varphi}\neq 0$ in a suitable neighbourhood of $u_o$ and the same is true for $\widetilde{\psi}$. Moreover, if $u \in S_{1}(F)$ then $\varphi(u)\in N(F'(u))$ because $\varphi$ is a kernel fibering map for $F$, hence $\widetilde{\varphi}(u)=\alpha(u)\varphi(u)+a(u)=\alpha(u) \varphi(u)\in N(F'(u))$ for $u \in S_1(F)$. Analogously, $\widetilde{\psi}(u)\in R(F'(u))^\perp$ for $u \in S_1(F)$. Thus $(\widetilde{\varphi},\widetilde{\psi}) \in \mathcal{P}(F,u_o)$.\\
\indent
3) We have to show that there exists $(\varphi,\psi) \in \mathcal{P}(F,u_o)$ such that $S_1(F)$ coincides, near   $u_o$, with the zero-set of the functional $J_0$ associated with $(\varphi,\psi)$. To this end it suffices to choose $(\varphi,\psi)=\mathcal{T}[\gamma^{-1},\delta^{-1}](\varphi_{C}, \psi_{C})$, as in the proof of 1). In fact, from formulas (\ref{148}) we get, near $u_o$,
\begin{equation*}
J_{0}(\varphi,\psi)(u)=J_{0}(\mathcal{T}[\gamma^{-1},\delta^{-1}](\varphi_{C}, \psi_{C}))(u)=J_{0}(\varphi_{C}, \psi_{C})(\gamma(u))=J_{0,C}(\gamma(u)).
\end{equation*}
The characterization (\ref{136}) implies that
\begin{equation*}
J_{0}(\varphi,\psi)(u)=0 \: \Leftrightarrow \: J_{0,C}(\gamma(u))=0 \Leftrightarrow \gamma(u)\in S_{1}(\varPhi) \: \Leftrightarrow \: u\in S_{1}(F).
\end{equation*}
Finally, if $(\widetilde{\varphi},\widetilde{\psi})\in \mathcal{P}(F,u_o)$ there exists $(\overline{\varphi},\overline{\psi}) \in \mathcal{P}(\varPhi,(0,0))$ such that $(\widetilde{\varphi},\widetilde{\psi})=\mathcal{T}[\gamma^{-1},\delta^{-1}](\overline{\varphi},\overline{\psi})$. Again, from (\ref{148}) we obtain that $J_{0}(\widetilde{\varphi},\widetilde{\psi})(u)=J_{0}(\overline{\varphi},\overline{\psi})(\gamma(u))$. On the other hand, from (\ref{137}) we know that $\gamma(u)\in S_{1}(\varPhi) \, \Rightarrow \,J_{0}(\overline{\varphi},\overline{\psi})(\gamma(u))=0$ and since $u \in S_{1}(F)\, \Leftrightarrow \, \gamma(u) \in S_{1}(\varPhi)$ we can conclude that $J_{0}(\widetilde{\varphi},\widetilde{\psi})(u)=0. \bracevert\\$\vspace{8pt}
\section{Classification of Singularities}\label{s2}
\subsection{Main Definitions and Differential Examples}\label{ss21}
\quad Here we introduce, with Definition \ref{D211}, our classification of singularities for a given Fredholm map. Subsequently we discuss some examples of differential problems exhibiting all of the considered singularities.\\
\indent The classification of simple singularities is made possible by the tools developed in the previous chapter. We point out that the adopted approach is quite different from those presented in other papers of Singularity Theory because of its local character, owing to the very definition of fibering pairs and fibering functionals. The previous approaches use pointwise, algebraic or analytical, conditions. By contrast, the analytic-geometric nature of our classification allows us to study singularities of all possible orders, at least when considering $C^{\infty}$ maps. The deep geometrical meaning of this kind of conditions will be made clear by the study of the stratification of singularities given in the following sections.\\
For the definition of a fibering pair $(\varphi,\psi)$ we refer to Definition \ref{d113}, while the related fibering functionals are introduced with Definition \ref{D116}. We also recall that in Theorem \ref{Teo146} we proved that when $u_o \in S_{1}(F)$, i.e. $u_o$ is a simple singularity for $F$, the set of all germs of $C^{d-1}$ fibering pairs for $F$ at $u_o$, denoted by $\mathcal{P}(F,u_o)$, is not empty.\\
\begin{definition}
\label{D211} Let  $U,V$  be open subsets in the $B$-spaces  $X,Y$, $F:U \subseteq X \rightarrow V \subseteq Y$ a  $C^d$ 0-Fredholm map $(d \geq 2)$ and $u_o \in S_{1}(F)$. Suppose that $(\varphi,\psi) \in \mathcal{P} (F,u_o)$ and let $k$ be an integer such that $1\leq k\leq d-1$. Then:\\
\begin{align*}
(\text{T}_k) \qquad & u_o \text{ is a } k\!-\!transverse  \: singularity  \text{ if}\\
& J_{0}(\varphi,\psi)(u_o)=\ldots=J_{k-1}(\varphi,\psi)(u_o)=0;\\
& I_{1}(\varphi,\psi)(u_o), \ldots, I_{k}(\varphi,\psi)(u_o) \text{ are } linearly \: independent \text{ (or l.i.).}\\
\\
(\text{S}_k) \qquad & u_o \text{ is a } k\!-\!singularity, \text{ or } k\!-\!ordinary \: singularity, \text{ if}\\
& J_{0}(\varphi,\psi)(u_o)=\ldots=J_{k-1}(\varphi,\psi)(u_o)=0, J_{k}(\varphi,\psi)(u_o)\neq 0;\\
& I_{1}(\varphi,\psi)(u_o),\ldots, I_{k-1}(\varphi,\psi)(u_o) \text{ are l.i.} \\
&\text{(the last condition is empty for } k = 1).\\
\\
(\text{M}_k) \qquad & u_o \text{ is a } maximal \: k\!-\!transverse \: singularity \text{ if}\\
& J_{0}(\varphi,\psi)(u_o)=\ldots=J_{k}(\varphi,\psi)(u_o)=0;\\
& I_{1}(\varphi,\psi)(u_o), \ldots, I_{k}(\varphi,\psi)(u_o)  \text{ are l.i.};\\
& I_{1}(\varphi,\psi)(u_o), \ldots, I_{k+1}(\varphi,\psi)(u_o) \text{ are not l.i.} \\
&\text{(the last condition is empty for } k =d-1).
\end{align*}
Finally, when $d=\infty \text{ or } d=\omega$:
\begin{align*}
\hspace{-10 cm}(\text{T}_\infty) \qquad &u_o \text{ is an } \infty\!-\!transverse \: singularity \text{ if} \qquad\qquad\qquad\\
&\text{for all integers } k\geq 1 \text{ one has}\\
& J_{0}(\varphi,\psi)(u_o)=\ldots=J_{k-1}(\varphi,\psi)(u_o)=0;\\
& I_{1}(\varphi,\psi)(u_o), \ldots, I_{k}(\varphi,\psi)(u_o) \text{ are l.i. }.\\
\end{align*}
\end{definition}
\vspace{4pt}
\indent Of course we have to check that the above definitions do not depend on the chosen fibering pair $(\varphi,\psi)\in \mathcal{P}(F,u_o)$: this is proved in Theorem \ref{Teo243}.\\
It is immediate to verify that if $u_o \in S_{1}(F)$ satisfies condition (T$_k$), i.e. it is a $k$-transverse singularity, then it is also an $h$-transverse singularity for $1\leq h \leq k$. It is also not difficult to prove that if $u_o \in S_{1}(F)$ satisfies condition (S$_k$) or (M$_k$) for an integer $k$ then this integer is unique (cf. Proposition \ref{Pro252}). Of course if $u_o$ satisfies condition (M$_k$) then it is also a $k$-transverse singularity. Moreover, we will prove in Lemma \ref{Lem244} that a $k$-singularity is a $k$-transverse singularity. Since an $\infty$-transverse singularity satisfies condition (T$_k$) for all integers $k \geq 1$, we can conclude that if $u_o \in S_{1}(F)$, where $F$ is a $C^\infty$ 0-Fredholm map, satisfies one of conditions (S$_k$),(M$_k$),(T$_\infty$) then $u_o$ satisfies condition (T$_1$), i.e. it is a 1-transverse singularity. Hence it is reasonable to raise the question whether the vice versa is true. As we will show in Proposition\ref{Pro252}, for a given $C^\infty$ 0-Fredholm map $F$ all 1-transverse singularities verify one and only one of conditions (S$_k$),(M$_k$),(T$_\infty$). In this sense, our classification of singularities can be considered as  \textit{complete}. Moreover, we will prove in the next section that the singular set near a 1-transverse singularity is a hypersurface; this important property justifies the interest for this kind of singular points.\\
Finally, a desirable property would be the invariance of this classification with respect to changes of coordinates. In the following two sections, by using the Algebraic and Geometric Lemmas (Section \ref{ss22}), we analyze the structure of the set of simple singularities and in Section \ref{ss24} we obtain a positive answer to all the above-listed issues. Specifically, simple singularities which are 1-transverse can be partitioned in: ordinary singularities, maximal transverse singularities and $\infty$-transverse singularities. Note that $k$-ordinary singularities, which are the analogues of Morin singularities in the finite-dimensional case (cf.\cite{Mo}), have been studied for small values of $k$ in several papers, quoted in the Introduction, and for all integers $k$ in \cite{Da}. On the other hand, to our knowledge, maximal transverse and $\infty$-transverse singularities have not been considered before. It is important to remark that $\infty$-transverse singularities can only occur in infinite dimensions.\\
\par In the rest of this section we briefly present or review some examples of differential problems which exhibit the singularities introduced above.\\
\indent
The first three examples show that there exist $\infty$-transverse singularities. For a first order periodic problem such as
\begin{equation*}
(\text{P}_{5})\begin{cases}
u'+a(t)(u-u^3)=h \;\; \text{ in\:\:} (0,1)\\
u(0) = u(1)\,, \end{cases}
\end{equation*}
where $h \in C^{0}([0,1]),u \in C^{1}([0,1])$ and the coefficient $a(t) \in C^{0}([0,1])\setminus \{0\}$ has mean value equal to zero, we can prove that $u \equiv +1$ and $u\equiv -1$ are $\infty$-transverse singularities for the associated map. The same conclusion can also be obtained for the point $u\equiv 0$ in the second order Neumann problem presented in the introduction\\
\begin{equation*}\label{}
(\text{P}_{3})\begin{cases} u'' + u u' = h \;\; \text{ in} \;\; (0,1) \\
u'(0) = 0\\
u'(1)= 0\,, \end{cases}
\end{equation*}\\
where $h \in C^{0}([0,1])$ and $u\in C^{2}([0,1]).$\\
However, even for these relatively simple ODEs, it seems difficult to prove the above statements by directly using condition $(T_\infty)$ of Definition \ref{D211}. In \cite{B-D 3} we prove these results by means of a consequence of the dual characterization of singularities.\\
\indent  A classical problem with a similar behaviour is given by the ``short'' pendulum without friction, i.e.  
\begin{equation*}
(\text{P}_{2})\begin{cases} u'' + A \sin u = h \;\; \text{ in} \;\; (0,T) \\
u(0) = u(T)\\
u'(0)= u'(T) \end{cases}
\end{equation*}
with $h \in C^{0}([0,T])$ and $u \in C^{2}([0,T])$. In the introduction we already recalled that, when $A > (\frac{2\pi}{T})^{2}$, problem $(\text{P}_{2})$ has a non-constant solution $\tilde{u}$  
to the right-hand side $h\equiv 0$. We also know that $\tilde{u}$ is an $\infty$-transverse singularity for the map associated with problem $(\text{P}_{2})$, \cite{B-D 4}. Note that the fact that $\tilde{u}$ is an $\infty$-transverse singularity for the associated map implies that, for any given $N \in \mathbb{N}$, there exists a right-hand side $\tilde{h}$ near 0 such that, for $h = \tilde{h}$, problem $(\text{P}_{2})$ admits at least $N$ distinct solutions near $\tilde{u}$ (cf. \cite{B-D 2}, Chapter 1, where the local behaviour near any kind of singular point is discussed).\\
\indent
In order to give examples of maximal $k$-transverse singularities for all integers $k$ we generalize the Riccati problem $(\text{P}_{1})$ stated in the introduction to the class of problems 
\begin{equation*}
(\text{P}_{\alpha,\beta})\begin{cases}
u'+a(t)u^{\alpha}e^{\beta u}=h  \qquad  \text{ in\:\:} (0,1)\\
u(0) = u(1), u>0 \;\; , \end{cases}
\end{equation*}
$\alpha, \beta \in \mathbb{R},\; h \in C^{0}([0,1]), u \in C^{1}([0,1])$ and $a(t) \in C^{0}([0,1]) \setminus \{0\}$ with mean value equal to zero. With regard to the above problem there exist singular points $u$ that are maximal 1-transverse singularities for the map associated with problem $(\text{P}_{\alpha,\beta})$ when $\alpha=0$ and $\beta \neq 0$, and singular points $u$ that are maximal $k$-transverse singularities, $k \geq 2$, when $\alpha=\frac{k}{k-1}$ and $\beta=0$, provided the function $a(t)$ is suitably chosen according to the values of $\alpha$ and $\beta$, \cite{B-D 4}. \\
\indent Finally, if one wishes to give explicit differential examples of all $k$-singularities it is possible to extend problem $(\text{P}_{3})$ to 
\begin{equation*}
(\text{P}_{6})\begin{cases}
u''+ uu' + u^{k+1}=h  \qquad  \text{ in\:\:} (0,1)\\
u'(0) = 0\\
u'(1) = 0 \;\; , \end{cases}
\end{equation*}
where $k\geq 1$ is an integer. It can be shown that $u\equiv 0$ is a $k$-singularity for the map associated with the above problem: this is proven in \cite{Ba4}, by using a tecnique completely different.\\
\par Importantly, in some cases it is possible to study the singularities of a given problem by direct use of the conditions in Definition \ref{D211}. For example let us consider the first order periodic problem  
\begin{equation*}
(\text{P}_{7})\begin{cases}
u'+a(t)u^{2}+p(t)u^{4}=h \qquad \text{ in\:\:} (0,1)\\
u(0)=u(1)  \end{cases}
\end{equation*}
with $a(t) \in C_{\#}^{1}([0,1])\setminus\{0\}$ and mean value zero, $p(t) \in C^{0}([0,1])$. For the map $F:C_{\#}^{1}([0,1])\rightarrow C^{0}([0,1])$ associated with $(\text{P}_{7})$ and defined as $F(u)=u'+g(t,u)$, where $g(t,u)=a(t)u^{2}+p(t)u^{4}$, we are interested in proving the following facts:\\
\indent i) when $p\equiv 0$ then $u\equiv 0$ is a maximal 2-transverse singularity for the map $F$;\\
\indent ii) when $\int_{0}^{1}p(t)\neq 0$ then $u\equiv 0$ is a 3-singularity for $F$.\\
\begin{remark} We note that, when $p\equiv 0$, $(\text{P}_{7})$ reduces to the Riccati problem $(\text{P}_{1})$ discussed in the introduction. Indeed, the result in i) recovers some of the properties for $(\text{P}_{1})$ initially shown in \cite{C-D 2}: by integration of the equation one can easily show that the solutions to $h\equiv 0$ are a continuum described by the unbounded real-analytic curve $\alpha \mapsto u_{\alpha}:=\dfrac{\alpha}{1+\alpha \int_{0}^{t}a(\tau)d\tau}$, defined for $
\alpha_{-}<\alpha<\alpha_{+}$ where $\alpha_{-}<0$ and $\alpha_{+} >0$ are suitable constants depending on the size of the function $a$. Furthermore, $u_{\alpha} \in C_{\#}^{1}([0,1])$ and the curve passes through 0, namely $u_{o}\equiv 0$.\\
\indent As for the result in ii), to our knowledge this is the first example of a 3-singularity for the class of first order differential problems. For an example of 3-singularity in the class of second order differential problems we refer to Theorem 2.5.1 in \cite{B-D 2}, where a different method of proof is adopted. Since $u\equiv 0$ is a 3-singularity for $F$ it is then possible to give an accurate description of the behaviour of the map $F$ between suitable neighbourhoods of $u\equiv 0$ and $h\equiv 0$: in particular, in these neighbourhoods one cannot have more than four solutions for the initial problem (cf. \cite{B-D 2}, Section 1.1). Instead, when $u\equiv 0$ is a maximal 2-transverse singularity for $F$, as in case i), we can only claim that there exist right-hand sides near $h\equiv 0$ which have either three solutions or a curve of solutions in a neighbourhood of $u\equiv 0$  (see again \cite{B-D 2}, Section 1.1). To the authors' knowledge, the first example of maximal 2-transverse singularities was shown in \cite{C-D 2}, even though those singularities were not classified as such yet. There a global description was given; in particular, it was shown that there are no right-hand sides with three solutions and also that a right-hand side has a curve of solutions \textit{iff} it is the image of a maximal 2-transverse singularity.
\end{remark}
\par We now briefly sketch the proof of properties i) and ii) for problem $(\text{P}_{7})$. As said before, in order to prove that $u\equiv 0$ is a maximal 2-transverse singularity or a 3-singularity for problem $(\text{P}_{7})$, we can directly use conditions (M$_2$) or (S$_3$) of Definition \ref{D211}. Let $(\varphi,\psi)$ be a fibering pair for $F$ defined near $u\equiv 0$. We saw in Example \ref{Ex115} that we can choose the following (global) fibering pair for $F$:
\begin{align*}
&(\varphi(u),\psi(u))=\\
&=\big(\text{exp}[t\cdotp \!\!\int_{0}^{t}\!g'(t,u)]\cdotp \text{exp}[-\! \!\int_{0}^{t}\!g'(\tau,u)], \text{exp}[-t \cdotp \!\! \int_{0}^{1}\!g'(t,u)] \cdotp \text{exp}[\int_{0}^{t}\!g'(\tau,u)]\big)\,.
\end{align*}
As discussed in Example \ref{Ex118}, this choice implies that
\begin{equation*}
J_{0}(u)= \int_{0}^{1}\!g'(t,u)dt,\,u\in C_{\#}^{1}([0,1]);
\end{equation*}
for the sake of simplicity, from now on we shall write $g^{(k)}(t,u)\equiv \dfrac{\partial^{k}g}{\partial u^{k}}(t,u)$ for $k \geq 1$ and, in the integrals, $u,v$ instead of $u(\cdotp), v(\cdotp)$. We shall also omit the dependence of the functionals on the fibering pair. From the last formula it follows that 
\begin{equation*}
I_{1}(u)v= \int_{0}^{1} \! g''(t,u)v\,dt, \;\; \text{for all} \;\; v \in C_{\#}^{1}([0,1])\; ,
\end{equation*}
and this implies that
\begin{equation*}
J_{1}(u)= \int_{0}^{1} \! g''(t,u)\varphi(u) dt. \medskip
\end{equation*}
Moreover, for all $v \in C_{\#}^{1}([0,1])$ we have that
\begin{equation*}
I_{2}(u)v= \int_{0}^{1} \! g'''(t,u)\varphi(u)v\,dt+\int_{0}^{1}\!g''(t,u)\varphi'(u)(v)dt \; ,
\end{equation*}
and this entails 
\begin{equation*}
J_{2}(u)= \int_{0}^{1}g'''(t,u)\varphi(u)^{2}dt+\int_{0}^{1}\!g''(t,u)\varphi_{1}(u)dt 
\end{equation*}
where $\varphi_{1}(u):=\varphi'(u)(\varphi(u))$, see also Definition 1.1.7. In a similar way we get that
\begin{align*}
I_{3}(u)v&= \int_{0}^{1}\!g^{(4)}(t,u)\varphi(u)^{2}v\,dt + \int_{0}^{1}g'''(t,u)2\varphi'(u)(v)dt \:+\\
&+ \int_{0}^{1}g'''(t,u)\varphi_{1}(u)v\,dt + \int_{0}^{1}g''(t,u)\varphi_{1}^{'}(u)(v) dt \;,
\end{align*}
and
\begin{equation*}
J_{3}(u)=\! \! \int_{0}^{1}\!g^{(4)}(t,u)\varphi(u)^{3}dt + 3 \!\!\int_{0}^{1}\! g'''(t,u)\varphi(u)\varphi_{1}(u) dt + \int_{0}^{1}g''(t,u)\varphi_{2}(u) dt 
\end{equation*}
where $\varphi_{2}(u):=\varphi^{'}_{1}(u)(\varphi(u))$ (cf. Definition \ref{D117}).\\
Since we have to evaluate the above functionals at $u = 0$ we only need to know the expressions of $\varphi(0), \varphi'(0)(v), \varphi_{1}(0),\varphi^{'}_{1} (0)(v)$ and $\varphi_{2}(0)$. To this end it is convenient to introduce the subspaces
\begin{align*}
&N_{1}:=\{v \in C_{\#}^{1}([0,1]): \int_{0}^{1}a(t)v\,dt = 0 \} \, ,\\
&N_{2}:=\{v \in C_{\#}^{1}([0,1]): \int_{0}^{1}(\int_{0}^{t}a(\tau) d\tau) a(t)v\,dt=0 \}
\end{align*}
and to observe that $u\equiv 1 \in N_{1}\cap N_{2}$, given that $a(t)$ has mean value zero.
By definition $\varphi(u)=\text{exp}[t\cdotp \int_{0}^{1}g'(t,u)- \int_{0}^{t}g'(\tau,u)]$, hence 
\begin{equation*}
\varphi(0)=1 .
\end{equation*}
By differentiating $\varphi(u)$ we get $\varphi'(u)(v)=\varphi(u)[t\cdotp \int_{0}^{1}g''(t,u)v - \int_{0}^{t} g''(\tau,u)v]$ and $\varphi_{1}(u) = \varphi(u)[t \cdotp \int_{0}^{1} g''(t,u)\varphi(u)- \! \int_{0}^{t} g''(\tau,u)\varphi(u)]$. It follows that 
\begin{equation*}
\varphi'(0)(v)=t \int_{0}^{1}2a(t)v\,dt - \int_{0}^{t} 2a(\tau)v\,d\tau = - \int_{0}^{t}2a(\tau)v \,d\tau, \;\; v \in N_{1} \;\; ,
\end{equation*}
\begin{equation*}
\varphi_{1}(0)= - \int_{0}^{t}2a(\tau) d\tau \;\;.
\end{equation*}
Again, by differentiating $\varphi_{1}(u)$, similar computations show that
\begin{align*}
\varphi^{'}_{1}(0)(v)=(\int_{0}^{t}\!2a(\tau) d\tau) (\int_{0}^{t} &2a(\tau)v\,d\tau)\!+\! \int_{0}^{t}(\int_{0}^{\tau}\!2a(\sigma)d\sigma)2a(\tau)v\,d\tau, v \in N_{1}\cap N_{2},\\
&\varphi_{2}(0)= \frac{3}{2}(\int_{0}^{t}\!2a(\tau) d\tau)^{2}\;. 
\end{align*}
By evaluating the functionals at $u=0$ and using the previous formulas we get: \medskip \\
$J_{0}(0)= \int_{0}^{1} \! g'(t,0)dt=0 \; ,$ \medskip\\
$I_{1}(0)v= \int_{0}^{1} \! 2a(t)v\,dt \; ,$ \medskip\\
$J_{1}(0)= \int_{0}^{1} \! 2a(t)dt=0 \; ,$ \medskip\\
$I_{2}(0)v=\! \int_{0}^{1}\! 2a(t)(t\int_{0}^{1}\! 2a(\tau)v\,d\tau -\! \int_{0}^{t}\! 2a(\tau)v \,d\tau)dt = \!- \int_{0}^{1}(\int_{0}^{t} 2a(\tau)d\tau)2a(t)v\,dt,$\medskip \\
(where the last equality holds for $v \in N_{1}$ and we integrated by parts to obtain it), \medskip \\
$J_{2}(0)= - \! \int_{0}^{1}(\int_{0}^{t}2a(\tau)d\tau)2a(t)dt=-\frac{1}{2}(\int_{0}^{1}2a(\tau)d\tau)^{2}=0 .$\medskip \\
Since $a(t)$ is nonzero it is not difficult to see that the functions $a(t)$ and $(\int_{0}^{t}a(\tau)d\tau)a(t)$  are linearly independent in $C_{\#}^{1}([0,1])$. Thus one can find a function $v \in C_{\#}^{1}([0,1])$ such that $I_{1}(0)v = 0$ and $I_{2}(0)v \neq 0$: this implies that $I_{1}(0)$ and $I_{2}(0)$ are linearly independent (cf. also the Algebraic Lemma in the next section, Lemma \ref{Lem221}).\\
\indent In order to prove property i) we now assume that $p\equiv 0$. Then, for $v \in N_{1}\cap N_{2}$ ,
\begin{align*}
I_{3}(0)v&=\int_{0}^{1}2a(t)[(\int_{0}^{t}2a(\tau)d\tau)(\int_{0}^{t}2a(\tau)v\,d\tau)+\int_{0}^{t}(\int_{0}^{t}2a(\sigma)d\sigma)2a(\tau)v\,d\tau]dt =\\
&=\int_{0}^{1}\frac{d}{dt}[2a(t)\int_{0}^{t}(\int_{0}^{\tau}2a(\sigma)d\sigma)2a(\tau)v\,d\tau]dt=0 .  
\end{align*}
Hence, when $I_{1}(0)v=0$ and $I_{2}(0)v=0$, one also has that $I_{3}(0)v=0$ and therefore $I_{1}(0)$, $I_{2}(0)$ and $I_{3}(0)$ must be linearly dependent (see also Algebraic Lemma). Summarizing, we proved that $J_{0}(0)=J_{1}(0)=J_{2}(0)=0$, the functionals $I_{1}(0)$ and $I_{2}(0)$ are l.i. but $I_{1}(0), I_{2}(0)$ and $I_{3}(0)$ are not l.i.. Thus condition $(\text{M}_2)$ of Definition 2.1.1 is satisfied and $u\equiv 0$ is a maximal 2-transverse singularity for $F$.\\
\indent
Finally, to prove property ii), under the hypothesis $\int_{0}^{1}p(t)dt\neq 0$ we also have that 
\begin{equation*}
J_{3}(0)= \int_{0}^{1}[24p(t)+2a(t)\frac{3}{2}(\int_{0}^{t}2a(\tau)d\tau)^2]dt=\int_{0}^{1}24p(t)dt \neq 0.
\end{equation*}
This allows us to conclude that $u\equiv 0$ is a 3-singularity for $F$ because $J_{0}(0)=J_{1}(0)=J_{2}(0)=0, J_{3}(0)\neq 0$ and $I_{1}(0), I_{2}(0)$ are l.i.: thus condition $(\text{S}_3)$ of Definition 2.1.1 is satisfied and the desired conclusion follows.\\
\begin{remark} Problem $(\text{P}_{7})$ can be generalized to include a linear term, i.e. it can be rewritten as 
\begin{equation*}
(\text{P}_{7}^{\,\prime})\begin{cases}
u'+q(t)u+a(t)u^{2}+p(t)u^{4}=h \qquad \text{ in\:\:} (0,1)\\
u(0)=u(1)\; ,  \end{cases}
\end{equation*}
where $a(t)$ and $p(t)$ are as above, while $q(t)\in C^{0}([0,1])$ has mean value zero and satisfies the relation $\int_{0}^{1}a(t)\cdot \text{exp}[-\int_{0}^{t}q(\tau)d\tau]dt = 0$. Then one can show that\\
\indent
j) when $p\equiv 0$ then $u\equiv 0$ is a maximal 2-transverse singularity for the map $F$ associated with the problem;\\
\indent
jj) when $p(t)\geq 0$ (or $p(t)\leq 0),p\neq 0\in C^{0}([0,1])$, then $u\equiv 0$ is a 3-singularity for $F$.\\
Properties j) and jj) of problem $(\text{P}_{7}^{\,\prime})$ can be proved in the same way of properties i) and ii) of $(\text{P}_{7})$: however, computations will be more complex because of $q(t)$ (for instance, in this case $\varphi(0) = \text{exp}[-\int_{0}^{t}q(\tau)d\tau]$).
\end{remark}
\vspace{2pt}
\subsection{Properties of Functionals}\label{ss22}
\quad Here we wish to study some algebraic and geometric properties of linear and nonlinear functionals which we will use extensively in the following. We point out that by applying the algebraic and geometric results to the fibering functionals (cf. Definition \ref{D116}) we derive, in the following sections, useful consequences for our classification. The first tool is an algebraic lemma which was partly suggested by lemma III.2 of \cite{Br}.

\begin{lemma}
\label{Lem221} (Algebraic Lemma). \textit{Let} $X$ \textit{be a vector space and} $I_{1}, \ldots, I_{k}$ \textit{linear functionals on} $X$. \textit{Then the following statements are equivalent}: 
\begin{enumerate}
\item[a)] $I_{1}, \ldots, I_{k}$ \textit{are linearly independent  (or }l.i.);
\item[b)] \textit{the linear map} $T:=(I_{1}, \ldots, I_{k}) : X \rightarrow \mathbb{R}^{k}$ \textit{is surjective};
\item[c)]  \textit{there exist vectors }$v_{j}\in X, j = 1,\ldots, k$,\textit{ such that the} $k \times k$  matrix
$\{I_{h}v_{j}\}_{h,j=1,\ldots,k}$  \textit{is non-singular};
\item[d)] \textit{there exist  vectors } $w_{j}\in X, j = 1,\ldots, k$ \textit{such that } $I_{h}w_{j}=\delta_{hj}, h, j = 1, \ldots, k$, \textit{where} $\delta_{hj}$ \textit{is the Kronecker delta};
\item[e)] $I_{h}\lvert_{\cap_{j=0}^{h-1}N(I_{j})}\;\neq 0,h = 1, \ldots, k$ (\textit{where} $N(I_{0}):=X$) ;
\item[f)] \textit{there exists}  $v_{k}\in \cap_{j=0}^{k-1}N(I_{j})$ \textit{such that }$I_{k}v_{k} \neq 0$  \textit{and, for} $k\geq 2, I_{1},\ldots, I_{k-1}$ \textit{are} l.i.;
\item[g)] \textit{for the map} $T$ \textit{we have that} $N(T)=\cap _{j=1}^{k}N(I_{j})$ \textit{with} codim $N(T)=k$.
\end{enumerate}
\end{lemma}
\vspace{8pt}
\par
\textbf{Proof.}\\
a) $\Rightarrow$ b): Let us suppose that the map $T$ is not surjective, i.e. $R(T)=\{Tu=(I_{1}u,\ldots,I_{k}u):u\in X\}$ is a proper vector subspace of $\mathbb{R}^{k}$. Then there exists a non-zero vector $\alpha=(\alpha_{1},\ldots,\alpha_{k})$ of $\mathbb{R}^{k}$ such that $\alpha \perp R(T)$, that is $\alpha_{1}I_{1}u + \ldots + \alpha_{k}I_{k}u = 0,\; \forall u \in X$. This gives $ \sum_{j=1}^{k} \alpha_{j}I_{j}=0$ with some $\alpha_{j}$ different from zero, hence $I_{1},\ldots,I_{k}$ are not  l.i..\\
b) $\Rightarrow$ c): Since $T$ is surjective, for every $j=1,\ldots,k$ there exists $v_{j}\in X$ such that $Tv_{j} =e_{j}$, where $e_{j}:=(0,\ldots,0,1,0,\ldots,0)$, with 1 at the j-th place. Hence $(I_{1}v_{j},\ldots,I_{k}v_{j}) = (0,\ldots,0,1,0 \ldots,0)$, i.e. $I_{h}v_{j}=\delta _{hj}$ with $h,j = 1,\ldots, k$, and so $\{I_{h}v_{j}\}_{h,j=1,\ldots,k}$ is the identity matrix. \\
c) $\Rightarrow$ d): Since $\{I_{h}v_{j}\}_{h,j=1,\ldots,k}$ is an invertible matrix, then the $k$ column vectors $ c_{j}:=(I_{1}v_{j},\ldots, I_{k}v_{j}),j=1,\ldots,k$, are a basis for $\mathbb{R}^{k}$. So there exist $\alpha _{ij}\in \mathbb{R}$, for $i,j=1,\ldots,k$, such that $e_{j}=\sum_{i=1}^{k}\alpha_{ij}c_{i}$. From this  
\begin{equation*}
e_{j}=\sum_{i=1}^{k}\alpha_{ij}(I_{1}v_{i}, \ldots, I_{k}v_{i})=(I_{1}(\sum_{i=1}^{k}\alpha_{ij}v_{i}),\ldots, I_{k}(\sum_{i=1}^{k}\alpha_{ij}v_{i})).
\end{equation*}
Then it suffices to take $ w_{j}:=\sum_{i=1}^{k}\alpha_{ij}v_{i}, j=1, \ldots, k .$ \\
d) $\Rightarrow$ e): Note that d) implies $I_{h}\lvert_{\cap_{j=1,j\neq h}^{k}N(I_{j})}\;\neq 0$, i.e. $I_{h}\lvert_{\cap_{j=0,j\neq h}^{k}N(I_{j})}\;\neq 0$, since by convention $ N(I_{0})=X$. In particular $ I_{h}\lvert_{\cap_{j=0}^{h-1}N(I_{j})}\;\neq 0 $.\\
e) $\Rightarrow$ f): For $h= k$ we have that $I_{k}\lvert_{\cap_{j=0}^{k-1}N(I_{j})}\;\neq 0$. This proves the existence of $v_{k}$. For  $ k\geq 2 $, e) implies the existence of vectors  $ v_{h}\in X, h = 1,\ldots, k-1 $, such that  
$ (I_{1},\ldots, I_{k-1})(v_{h})=(0,\ldots, 0, I_{h}(v_{h}),I_{h+1}(v_{h}),\ldots, I_{k-1}(v_{h}))$ with $h-1$  zeroes at the beginning and $ I_{h}(v_{h})\neq 0 $. If we let $w_{h}:=(I_{1},\ldots, I_{k-1})(v_{h})$  we obtain that $\{w_{1},\ldots, w_{k-1}\}$  is a basis of  $\mathbb{R}^{k-1} $. Let us choose  $ \alpha_{1},\ldots,\alpha_{ k-1}\in \mathbb{R} $  such that $\alpha_{1}I_{1} + \ldots + \alpha_{k-1}I_{k-1}=0 $. By evaluating this identity for $ v_{h}$ we have that $\alpha\perp w_{h}$, where  $\alpha:=(\alpha_{1},\ldots,\alpha_{k-1})\in \mathbb{R}^{k-1}$. Since $ \alpha \perp w_{h} $ for $ h=1,\ldots, k-1 $, we can say that $\alpha=0 \in \mathbb{R}^{k-1}$, i.e. $\alpha_{1}=\ldots=\alpha_{k-1}=0 $. Thus we conclude that  $ I_{1},\ldots, I_{k-1} $ are  l.i..\\
f) $\Rightarrow$ a): Let  $\alpha_{1}I_{1}+\ldots+\alpha_{k-1}I_{k-1}+\alpha_{k}I_{k}=0$ and  $(\alpha_{1},\ldots,\alpha_{k-1},\alpha_{k})\in \mathbb{R}^{k}$. We consider $v_{k}\in \cap_{j=0}^{k-1}N(I_{j}) $ such that $ I_{k}v_{k}\neq 0 $. Then $\alpha_{k}I_{k}v_{k}=0 $ and thus $\alpha_{k}=0 $. It follows that $\alpha_{1}I_{1}+\ldots+\alpha_{k-1}I_{k-1}=0$, and since $ I_{1},\ldots, I_{k-1} $ are l.i. we get that $ \alpha_{1}=\ldots= \alpha_{k-1}=0 $. This means that the functionals  $I_{1},\ldots, I_{k} $ are  l.i. since  $\alpha_{1}=\ldots=\alpha_{k-1}=\alpha_{k}=0 $.\\
g) $\Rightarrow$ b): By definition, $ Tu=(I_{1}u,\ldots,I_{k}u)=0\in \mathbb{R}^{k}$ \textit{ iff} $I_{j}u=0,j=1, \ldots,k$, thus  $ N(T)=\cap_{j=1}^{k}N(I_{j})$. Let us consider the map $ \tilde{T}:X/N(T)\rightarrow R(T) $ defined as $ \tilde{T}[u]:=T(u)$, where $ [u]\in X/N(T) $ is the equivalence class of $ u\in X $. Then $ \tilde{T} $ is a well-defined linear operator and, by construction, $ \tilde{T} $ is injective and surjective. Hence $ \tilde{T} $ is a linear isomorphism and so dim$\,X/N(T)=\text{dim}\,R(T)$ which, by definition of codimension, means that codim$\,N(T)=\text{dim}\,R(T) $. Thus $ T $ is surjective \textit{iff} codim$\,N(T)=\text{dim}\,R(T)=k. \bracevert$  \\
\par
A useful consequence of the Algebraic Lemma is the following
\begin{corollary}
\label{Co222} \textit{Let} $X$ \textit{be a normed vector space, }$\Omega$\textit{ a topological space,} $I_{1},\ldots,I_{k}:\Omega \rightarrow X^{\ast}$ \textit{continuous maps, and} $ \omega_{o} \in \Omega$  \textit{such that} $ I_{1}(\omega_{o}),\ldots,I_{k}(\omega_{o}) $ \textit{are} l.i. \textit{on} $X$.\! \textit{Then there exists a neighbourhood} $ U$\textit{ of} $\omega_{o}$ \textit{such that} $ I_{1}(\omega), \ldots,I_{k}(\omega) $ \textit{are} l.i. $ \forall \; \omega \in U $.
\end{corollary}
\vspace{4pt}
\textbf{Proof.} From the equivalence $ (a)\Leftrightarrow (c) $ of the previous lemma we have that there exist $ v_{j}\in X, j=1,\ldots,k,$ such that the matrix $ \{I_{h}(\omega_{o})v_{j}:h,j=1,\ldots,k\} $ is non-singular, that is $ \text{det}\{I_{h}(\omega_{o})v_{j}\}\neq 0$. By the continuity on $\Omega$ of this determinant there exists a neighbourhood $ U $ of $ \omega_{o} $ such that $ \text{det}\{I_{h}(\omega)v_{j}\}\neq 0, \; \forall \; \omega \in U $, that is $ \{I_{h}(\omega)v_{j}\} $ is a non-singular matrix. From the equivalence $ (a)\Leftrightarrow (c) $ we can conclude that $ I_{1}(\omega),\ldots, I_{k}(\omega) $ are  l.i.. $ \bracevert $ \\
\par
The above results about the independence of linear functionals turn out to be very useful when dealing with the zero-sets of real functionals defined on $B$-spaces.  More precisely, given $ C^{h} $ functionals $ f_{i} $ defined on a $ B $-space $ X, i=1,\ldots, k, $ with Fréchet derivatives satisfying a linear independence condition, we will see in the Geometric Lemma below that the intersection of the inverse images  $ f_{i}^{-1}(0) $ is a subset $ M $   which is actually a \textit{Banach manifold}. For the notion of $C^{h}$  Banach manifold we refer to \cite{A-M-R}, sections 1, 2 of chapter 3, where a thorough description of Banach manifolds of class $C^{h}$ is given as well as the notion of $C^{h}$ \textit{submanifold} of a manifold of class $C^{h}$. Moreover, it is interesting to consider \textit{maps of class} $C^{h}$ between $C^{h}$ manifolds and, from the notion of \textit{tangent space} of a manifold $ M $ at a point $ u $ (denoted by $ T_{u}M $), one can introduce the \textit{tangent} or \textit{derivative map} of a given $C^{h}$ map between $C^{h}$ manifolds (see e.g. \cite{A-M-R}, sections 2, 3 of chapter 3). We use the following definition of \textit{codimension} of a submanifold $S$ of a $C^{h}$ Banach manifold $ M $: 
\begin{equation*}
\text{codim}\, S:= \text{codim}_{T_{u}M}\,T_{u}S = \text{dim}\,T_{u}M/T_{u}S ,
\end{equation*}
which is well-defined if the right-hand side does not depend on $ u\in S $ (cf. \cite{Ze4}, section 11 of chapter 73).\\
\indent
We also recall a version of the Submersion Theorem suitable for our purposes (cf. \cite{Ze4}, theorem 73.C, for a more general statement and a detailed proof).\\
\\
 \textit{\textbf{Submersion Theorem.}} \textit{Let} $ M $ \textit{be  a} $C^{h}$ \textit{Banach manifold and} $ F : M \rightarrow \mathbb{R}^{k}$ \textit{ a} $C^{h}$ \textit{map,} $ h \geq 1 $. \textit{Let us assume that} $ h\in \mathbb{R}^{k}$ \textit{is a} regular value \textit{of} $ F $, i.e. $ F'(u) : T_{u}M \rightarrow \mathbb{R}^{k} $ \textit{is surjective for all} $ u \in F^{-1}(h)=\{w \in M : F(w) = h\} $. \textit{Then the set} $ S :=F^{-1}(h) $ \textit{is a} $C^{h}$ \textit{submanifold of}$  M  $\textit{ with tangent space given by} $ T_{u}S = N(F'(u)),\;\forall \; u \in S$.\\
 \par
Thanks to the above geometric notions and the Algebraic Lemma \ref{Lem221} we can prove the following

\begin{lemma} \label{Lem223} (Geometric Lemma). \textit{Let} $ X $ \textit{be a Banach space}, $ u_{o} \in X $, \textit{and} $ f_{1},\ldots, f_{k}\;C^{h}$ \textit{real functionals, defined on a neighbourhood} $ U $ \textit{of} $ u_{o} $, \textit{such that} $ f_{1}(u_{o}) =\ldots=f_{k}(u_{o}) = 0 $\textit{ with }$ f_{1}^{'}(u_{o}),\ldots, f_{k}^{'}(u_{o}) $ l.i \textit{on} $X^{\ast}$. \textit{If }$M:=\cap_{j=1}^{k} \{u\in U : f_{j}(u) = 0\}$ \textit{then} $ M $  \textit{is, near} $ u_{o} $,\textit{ a }$C^{h}$ \textit{k}-codimensional submanifold \textit{of }$X$ \textit{and for} $ u\in M $ \textit{we have that the tangent space to} $ M $ \textit{at} $ u $ \textit{is given by } $T_{u}M = \cap_{j=1}^{k} N(f_{j}^{'}(u))$.
\end{lemma}
\indent
\textbf{Proof.} We argue by induction. For $k = 1$ we have the usual statement about the zero-set of a differentiable functional having Fréchet derivative different from zero. This is proved either by a standard application of the Implicit Function Theorem, which allows one to represent $ M $ locally as a graph of a suitable $C^{h}$ map, or by using the Submersion Theorem.\\
\indent Now let us assume that the statement is true for $k\geq 1$. Then we shall prove it for the integer $k+1$. By hypothesis $f_{1}(u_{o})=\ldots=f_{k+1}(u_{o})=0$ and  $f_{1}^{'}(u_{o}),\ldots,f_{k+1}^{'}(u_{o})$ are l.i. and thus, by the inductive assumption for $k$, we have that $M:=\cap_{j=1}^{k} \{u\in U : f_{j}(u_{o})=0\}$ is, near $u_{o}$, a $C^{h}$ \textit{k}-codimensional submanifold of $X$ and $T_{u}M=\cap_{j=1}^{k} N(f_{j}^{'}(u)),u\in M$. Let us consider the restriction of $f_{k+1}$ to $M$, i.e. $f_{k+1}:M \rightarrow \mathbb{R}$, with  $f_{k+1}^{'}(u_{o}):T_{u_{o}}M=\cap _{j=1}^{k} N(f_{j}^{'}(u_{o}))\rightarrow \mathbb{R}$. This linear functional is not identically zero by the Algebraic Lemma \ref{Lem221}, a) $\Rightarrow$ e). Since $f_{k+1}(u_{o})=0$, again by the Submersion Theorem we obtain that $N:=\{u\in M:f_{k+1}(u)=0\}=\cap_{j=1}^{k+1}\{u\in U:f_{j}(u)=0\}$ is, near $u_{o}$, a $C^{h}$ submanifold of $M$ and, for $u\in N$, we have that $T_{u}N = \{v\in T_{u}M:v\in N(f_{k+1}^{'}(u))\} = \cap_{j=1}^{k+1}N(f_{j}^{'}(u))$. Hence $N$ is a 1-codimensional submanifold of $M$ and so $N$ is a submanifold of $X$. From the Algebraic Lemma it is quite easy to see that $N$ is of codimension $k+1$ in $X$. In fact, thanks to the equivalence g) $\Leftrightarrow$ b) the vector subspace $T_{u}N=\cap_{j=1}^{k+1}N(f_{j}^{'}(u))$ is of codimension $k+1$ in $X$ for all $u \in N$. $ \bracevert $\\
\par
It may be worthwhile to note that the Submersion Theorem can be directly used to prove the above lemma, too. In fact, let us consider the map $(f_{1},\ldots,f_{k}): U\rightarrow \mathbb{R}^{k}$. Its derivative at $u_{o}$ is $(f_{1}^{'}(u_{o}),\ldots, f_{k}^{'}(u_{o})): X \rightarrow \mathbb{R}^{k}$, which is surjective by the Algebraic Lemma. Thus $(f_{1}^{'}(u),\ldots, f_{k}^{'}(u)): X \rightarrow \mathbb{R}^{k}$ is surjective for $ u $ near $ u_{o} $ by the Algebraic Lemma and Corollary 2.2.2. Then, by the Submersion Theorem, $M= (f_{1},\ldots,f_{k})^{-1}(0)$  is a $C^{h}$ submanifold of  $X$ near $u_{o}$. Moreover, for $u$ near $u_{o},T_{u}M =\cap_{j=1}^{k} N(f_{j}^{'}(u))$ has codimension $k$ in $X$ because of the Algebraic Lemma, g) $\Leftrightarrow$ b). Hence $M$ is a \textit{k}-codimensional submanifold. \\
\par
The reason why we opted for an inductive proof of the Geometric Lemma \ref{Lem223} is illustrated in the following remark. This will be also useful when considering the stratification of singularities.

\begin{remark} \label{Rem224} Under the hypotheses of Lemma \ref{Lem223} let us define  $M_{0}:=X$  and, for $h=1,\ldots,k$, $M_{h} :=\cap_{j=1}^{h}\{u\in U:f_{j}(u_{o})=0\}$ (note that  $M_{k}\equiv M $). We claim that $M_{h}$ is, near $u_{o}$, a one-codimensional submanifold of $M_{h-1}$, for $h=1,\ldots,k$. In fact, the linear independence of $\{f_{j}^{'}(u_{o}), j=1,\ldots,k\}$ implies that $\{f_{j}^{'}(u_{o}),j=1,\ldots, h-1\}$ are l.i., for $h=2,\ldots,k$. Since the Geometric Lemma applies to the integer $h-1$ we get that $M_{h-1}$ is an ($h-1$)-codimensional submanifold of $X$. By considering $f_{h}:M_{h-1} \rightarrow \mathbb{R}$ and arguing as in the proof of the Geometric Lemma one can easily prove that  $M_{h}$  is a one-codimensional submanifold of  $M_{h-1},h =1,\ldots,k$.
\end{remark}
\par
Another way to look at the Geometric Lemma is through the notion of \textit{mutual transversality} (cf. \cite{ G-G} section 3, chapter III). We recall that the submanifolds $N_{j}$ of $X,j=1,\ldots,k$ where $k\geq 2$, are said to be in \textit{general position}, or \textit{mutually transversal}, at the point $u_{o}\in \cap_{j=1}^{k} N_{j}$ if
\begin{equation*}
\text{codim}(T_{u_{o}}N_{j_{1}}\cap \ldots \cap T_{u_{o}}N_{j_{s}})=\text{codim}\,T_{u_{o}}N_{j_{1}}+\ldots+\text{codim}\,T_{u_{o}}N_{j_{s}} ,
\end{equation*}
for every sequence of integers  $1\leq j_{1}<\ldots<j_{s}\leq k$, with $ 2\leq s\leq k$.\\
\indent Then, by using the Algebraic Lemma, we can see that, under the same hypotheses of the Geometric Lemma, the linear independence of $f_{1}^{'}(u_{o}),\ldots,f_{k}^{'}(u_{o})$ implies that the manifolds $N_{j}:=\{u \in U:f_{j}(u)=0\},j=1,\ldots,k$, are in general position at the point $u_{o}$.\\
To prove this fact we first remark that, since the functionals $f_{1}^{'}(u_{o}),\ldots,f_{k}^{'}(u_{o})$ are l.i., then each $f_{j}^{'}(u_{o})$ is not identically zero. Hence, by using the Geometric Lemma \ref{Lem223}, we get that each $N_{j}$ is, near $u_{o}$, a one-codimensional submanifold of $X$ and $T_{u_{o}}N_{j}=N(f_{j}^{'}(u_{o}))$. Thus codim$\,T_{u_{o}}N_{j}=\text{codim}\,N(f_{j}^{'}(u_{o}))=1,j=1,\ldots,k$. This implies that, if the integers $j_{1}, \ldots,j_{s}$ are chosen such that $1\leq j_{1}< \ldots < j_{s}\leq k$, it suffices to prove that
\begin{equation*}
\text{codim}(T_{u_{o}}N_{j_{1}}\cap \ldots \cap T_{u_{o}}N_{j_{s}})=s\,.
\end{equation*} 
In fact, since $f_{j_{1}}^{'}(u_{o}),\ldots,f_{j_{s}}^{'}(u_{o})$ are l.i. then, by the Algebraic Lemma \ref{Lem221}, a) $\Rightarrow$ g), we obtain that
\begin{equation*}
\text{codim}(T_{u_{o}}N_{j_{1}}\cap \ldots \cap T_{u_{o}}N_{j_{s}})=\text{codim}(N(f_{j_{1}}^{'}(u_{o})) \cap \ldots \cap N(f_{j_{s}}^{'}(u_{o})))=s ,
\end{equation*}
and thus we have that the manifolds $N_{j}$ are in general position.\\
\indent
Therefore the Geometric Lemma shows that the intersection of $k$ one-codimensional submanifolds $N_{j}$ of $X$ which are mutually transversal at $u_{o}$ is a \textit{k}-codimensional submanifold of $X$ near $u_{o}$; moreover, the tangent space at each point is the intersection of the tangent spaces. \vspace{2pt}
\subsection{1-transverse Singularities}\label{ss23}
\quad First of all, we want to prove that the definition of 1-transverse singularity given by condition (T$_1$) in Definition \ref{D211} is well-posed, i.e. it is independent from the chosen f-pair. Note that condition (T$_1$) is equivalent to $ J_{0}(\varphi,\psi)(u_{o})=0$ (which is always true from the very definition of $J_{0}$) and $I_{1}(\varphi,\psi)(u_{o})\neq 0$. To demonstrate the independence on the chosen f-pair we need the following technical result.\\
\begin{proposition}
\label{Pro231}\textit{ Let } $U,V$ \textit{be open subsets in the $B$-spaces} $X,Y$ \textit{and let}   $F:U\subseteq X \rightarrow V\subseteq Y$ \textit{be a } $C^{d}$ 0-\textit{Fredholm map} $(d\geqslant 2)$. \textit{If} $u_{o}\in S_{1}(F)$ \textit{then} $\forall u\in S_{1}(F)$, $u$ \textit{near} $u_{o}$, \textit{one has that for any} $(\varphi,\psi), (\widetilde{\varphi},\widetilde{\psi})\in \mathcal{P}(F,u_{o})$ 
\begin{equation*}
N(I_{1}(\varphi,\psi)(u))= N(I_{1}(\widetilde{\varphi},\widetilde{\psi})(u))=\{v \in X:F''(u)[v,N(F'(u))] \subseteq R(F'(u))\}.
\end{equation*}
\end{proposition}
\par
\textbf{Proof.} For a given $(\varphi,\psi)\in \mathcal{P}(F,u_{o})$ we have that every $(\widetilde{\varphi},\widetilde{\psi})\in \mathcal{P}(F,u_{o})$ can be represented, near $u_{o}$, as $\widetilde{\varphi}=\alpha\varphi+a,\widetilde{\psi}=\beta\psi+b$, where $\alpha, \beta \in C_{u_{o}}^{d-1}(X,\mathbb{R})$ are non-vanishing functions, and  $a \in C_{u_{o}}^{d-1}(X,X), b \in C_{u_{o}}^{d-1}(X,Y^{\ast})$ are equal to zero on $S_{1}(F)$ (as proved in Theorem \ref{Teo146}). By definition,  $J_{0}(\varphi,\psi)=\psi F'\varphi$ and $J_{0}(\widetilde{\varphi},\widetilde{\psi})=\widetilde{\psi}F'\widetilde{\varphi}$.\\
\indent We claim that: 
\begin{equation}\label{231}
\begin{split}
J_{0}(\widetilde{\varphi},\widetilde{\psi})&=\alpha\beta J_{0}(\varphi,\psi)+H ,\\
\text{with }H\in C_{u_{o}}^{d-1}(X,\mathbb{R})&,H=0\text{ and }H'=0\text{ on }S_{1}(F).
\end{split}
\end{equation}
In fact,  $J_{0}(\widetilde{\varphi},\widetilde{\psi})= \widetilde{\psi}F'\widetilde{\varphi}=(\beta\psi+b)F'(\alpha\varphi+a)=\alpha\beta J_{0}(\varphi,\psi)+H$, with $ H:= \beta\psi F'a+\alpha bF'\varphi+bF'a$. Note that every addendum of $H$ is obtained by composing maps such that two of them are equal to zero on $S_{1}(F)$. Precisely, by definition of fibering maps for $F$ and from the choice of $a,b$ we have that $\psi F',a,b,F'\varphi$ are zero on  $S_{1}(F)$. So, looking at the first addendum in $H$, one obtains $\beta\psi F'a=0$ on $S_{1}(F)$ and also $(\beta\psi F'a)'v=(\beta'v)\psi F'a+\beta((\psi F')'v)a +\beta \psi F'(a'v)=0$ on $S_{1}(F)$. To conclude it suffices to repeat this argument for the other two addenda.\\
\indent Since $I_{1}(\varphi,\psi)=J_{0}^{'}(\varphi,\psi)$ and  $I_{1}(\widetilde{\varphi},\widetilde{\psi})=J_{0}^{'}(\widetilde{\varphi},\widetilde{\psi})$ then from (\ref{231}) we get
\begin{equation*}
I_{1}(\widetilde{\varphi},\widetilde{\psi})v=((\alpha\beta)'v)J_{0}(\varphi,\psi)+\alpha\beta J_{0}^{'}(\varphi,\psi)v + H'v,\; v\in X.
\end{equation*}
Hence
\begin{equation}\label{232}
I_{1}(\widetilde{\varphi},\widetilde{\psi})v=\alpha\beta I_{1}(\varphi,\psi)v \;\; \text{on}\; S_{1}(F), \:\:\forall\; v\in X ,
\end{equation}
because $H'$ and $J_{0}(\varphi,\psi)$ vanish on $S_{1}(F)$.\\
\indent
Let us now compute $I_{1}(\varphi,\psi)$:\\
\begin{equation*}
I_{1}(\varphi,\psi)v=(\psi F'\varphi)'v=(\psi'v)F'\varphi+\psi F''[v,\varphi]+\psi F'(\varphi'v).
\end{equation*} 
Thus\\
\begin{equation}\label{233}
I_{1}(\varphi,\psi)v= \psi F''[v,\varphi]   \text{ on }  S_{1}(F), \;\; \forall \; v \in X ,
\end{equation}
because $\psi F',F'\varphi$ are zero on  $S_{1}(F)$. \\
Now let us recall that $\alpha(u_{o})\,\neq 0, \beta (u_{o})\,\neq 0, \varphi(u_{o}) $ spans $N(F'(u_{o}))$ and $\psi(u_{o})$ spans $R(F'(u_{o}))^{\perp}$, i.e. $\forall \, w \in Y$ one has $\psi(u_{o})w=0 \textit{ iff } w \in  R(F'(u_{o}))$. From this and (\ref{232}), (\ref{233}) we draw the desired conclusion. $\bracevert$\\
\par
An immediate consequence of the above result is the following\\
\begin{corollary}
\label{Co232}\textit{ Let }$U,V$ \textit{be open subsets in the $B$-spaces} $X,Y$ \textit{and let} $F: U\subseteq X\rightarrow V\subseteq Y$ \textit{be a} $C^{d}$ 0-\textit{Fredholm map} ($d\geq 2$). \textit{For the point} $u_{o}\in S_{1}(F)$ \textit{the following conditions are equivalent:}
\begin{enumerate}
\item[i)]\textit{there exists a} $C^{d-1}$ \textit{ f-pair for} $F$ \textit{near} $u_{o}$,\textit{ say} $(\varphi,\psi)\in \mathcal{P}(F,u_{o})$, \textit{such that } $I_{1}(\varphi,\psi)(u_{o})=\psi(u_{o})F''(u_{o})[\, \cdotp,\varphi(u_{o})]\,\neq 0$;
\item[ii)] \textit{for each } $(\widetilde{\varphi},\widetilde{\psi})\in \mathcal{P}(F,u_{o})$ \textit{ one has }\\
$I_{1}(\widetilde{\varphi},\widetilde{\psi})(u_{o})=\widetilde{\psi}(u_{o})F''(u_{o})[\cdotp,\widetilde{\varphi}
(u_{o})]\,\neq 0$. \vspace{5pt}\end{enumerate}
\end{corollary}
\noindent Note that, from what we said about $J_{0}$ at the beginning of the section, this means that the definition of 1-transverse singularity is well-posed.\\
\par
Now we wish to prove two properties of 1-transverse singularities. We shall see that, when $u_{o}$ is a 1-transverse singularity, the singular set $S_{1}(F)$ is locally a \textit{hypersurface} of $X$, i.e. $S_{1}(F)$ is a one-codimensional manifold near $u_{o}$. Moreover $S_{1}(F)$ coincides, near $u_{o}$, with the zero-set of the functional $J_{0}(\widetilde{\varphi},\widetilde{\psi})$ for any f-pair $(\widetilde{\varphi},\widetilde{\psi})$ defined on a neighbourhood of  $u_{o}$.\\
\begin{theorem}\label{Teo233} \textit{Let} $U,V$ \textit{be open subsets in the $B$-spaces} $X,Y, F: U\subseteq X\rightarrow V\subseteq Y$\textit{ a } $C^{d}$  0-\textit{Fredholm map} ($d\geq 2$) \textit{and let} $u_{o}\in  S_{1}(F)$ \textit{be} 1-transverse \textit{for} $F$. \textit{Then } $S_{1}(F)$ \textit{is, near} $u_{o}$,\textit{ a }$C^{d-1}$ \textit{one-codimensional submanifold of} $X$ \textit{and, for each} $u\in S_{1}(F)$ \textit{ and }$(\varphi,\psi)\in \mathcal{P}(F,u_{o})$ \textit{it follows that}
\begin{center}
$T_{u}S_{1}(F)=N(I_{1}(\varphi,\psi)(u))=\{v\in X:F''(u)[v,N(F'(u))]\subseteq R(F'(u))\}$.
\end{center}
\textit{Moreover, for every} $(\widetilde{\varphi},\widetilde{\psi})\in \mathcal{P}(F,u_{o})$ \textit{we have that, near} $u_{o}$, \textit{the following equivalence holds:}
\begin{equation*}
J_{0}(\widetilde{\varphi},\widetilde{\psi})(u)=0 \Leftrightarrow u\in S_{1}(F).
\end{equation*}
\end{theorem}
\indent
\textbf{Proof.} As we saw in  Theorem \ref{Teo146}  there exists $(\varphi,\psi)\in \mathcal{P}(F,u_{o})$ such that  
\begin{center}
$J_{0}(\varphi,\psi)(u)=0 \Leftrightarrow u\in S_{1}(F).$
\end{center}
Let us consider the $C^{d-1}$ functional $J_{0}(\varphi,\psi):U_{o}\subseteq X\rightarrow \mathbb{R}$, where $U_{o}\subseteq U$ is a suitable neighbourhood of $u_{o}$. By definition $J_{0}^{'}(\varphi,\psi)(u_{o})=I_{1}(\varphi,\psi)(u_{o})\,\neq 0$, given that $u_{o}$ is 1-transverse for $F$. Since  $J_{0}(\varphi,\psi)(u_{o})=0$ the Geometric Lemma \ref{Lem223} implies that, up to shrinking $U_{o}$, the subset $S_{1}(F)\cap U_{o}=\{u\in U_{o}:J_{0}(\varphi,\psi)(u)=0\}$ is a $C^{d-1}$ one-codimensional submanifold of $X$. Also, for  $u\in S_{1}(F)$ we have that $ T_{u}S_{1}(F)=N(J_{0}^{'}(\varphi,\psi))(u)=N(I_{1}(\varphi,\psi)(u))=\{v\in X:F''(u)[v,N(F'(u))] \subseteq R(F'(u))\}$ (the last equality follows from  Proposition \ref{Pro231}).\\
Now let $(\widetilde{\varphi},\widetilde{\psi})\in \mathcal{P}(F,u_{o})$ be another f-pair, and consider $J_{0}(\widetilde{\varphi},\widetilde{\psi}):\widetilde{U}_{o}\subseteq X\rightarrow \mathbb{R}$ with $\widetilde{U}_{o}$ a suitable neighbourhood of $u_{o}$ in $U$. We set $\Sigma_{1}(F):=\{u\in \widetilde{U}_{o}: J_{0}(\widetilde{\varphi},\widetilde{\psi})(u)=0\}$; $\Sigma_{1}(F)$ is not empty because $u_{o} \in S_{1}(F)\subseteq \Sigma_{1}(F)$ (the inclusion following from point 3) of Theorem \ref{Teo146}). By Corollary \ref{Co232} we get that $J_{0}^{'}(\widetilde{\varphi},\widetilde{\psi})(u_{o})=I_{1}(\widetilde{\varphi},\widetilde{\psi})(u_{o})\,\neq 0$. Then we can prove, in the same way as above, that $\Sigma_{1}(F)$ is, near $u_{o}$, a $C^{d-1}$ one-codimensional submanifold of $X$. Therefore $S_{1}(F)$ and $\Sigma_{1}(F)$ are, near $u_{o},\; C^{d-1}$ one-codimensional submanifolds of $X$ which contain $u_{o}$ and $S_{1}(F)\subseteq\Sigma_{1}(F)$. Hence it must be $S_{1}(F)=\Sigma_{1}(F)$. This quite plausible fact is an easy consequence of the Implicit Function Theorem and is shown below.$\bracevert$\\
\begin{lemma} \textit{Let} $X$ \textit{be a $B$-space and} $S,\Sigma$ \textit{two} $C^{1}$ \textit{one-codimensional submanifolds of} $X$ \textit{such that} $S\subseteq \Sigma$. \textit{Then, for any }$u_{o} \in S$,\textit{ it follows that} $S=\Sigma$\textit{ near }$u_{o}$.
\end{lemma}
\par
\textbf{Proof.} We can assume that $u_{o}=0 \in S$; hence, by hypothesis, we get that $0\in \Sigma$. Now we shall prove that for any $u\in \Sigma$ which is close enough to $0$ one has that $u\in S$.\\
By definition of submanifold we know that there exist a neighbourhood $U$ of $0$ and $C^{1}$ functionals $J,\vartheta:U\subseteq X\rightarrow \mathbb{R}$ such that $S\cap U=J^{-1}(0),\Sigma\cap U=\vartheta^{-1}(0)$ and  $T_{0}S=N(J'(0)),T_{0}\Sigma = N(\vartheta'(0))$.\\
We contend that $N(J'(0))=N(\vartheta'(0))$. Since  $N(J'(0)),N(\vartheta'(0))$ are 1-codimensional closed subspaces it suffices to show that $N(J'(0))\subseteq N(\vartheta'(0))$. Given $v\in N(J'(0))=T_{0}S$, by definition of tangent space there exists a $C^{1}$ curve $c:\mathbb{R}\rightarrow X$ such that $c(\mathbb{R})\subseteq S,c(0)=0$ and $c'(0)=v$. Since $S\subseteq \Sigma$ then, for small $t\in \mathbb{R}$, one has that $\vartheta(c(t))=0$. By differentiation we obtain that $\vartheta'(0)c'(0)=0$, hence $c'(0)=v\in N(\vartheta'(0))$.\\
Let us now set $N:=N(J'(0))=N(\vartheta'(0))$ and let $v_{o}\in X$ be such that $v_{o}\notin N$. In this way, $X=N\oplus \mathbb{R}v_{o}$. On a suitable neighbourhood $U'$ of the origin of $N\times \mathbb{R}$ we can consider the $C^{1}$ functional $f_{J}:U'\subseteq N \times \mathbb{R}\rightarrow \mathbb{R},f_{J}(n,r):=J(n+rv_{o})$. Then $f_{J}(0,0)=0$ and $\frac{\partial f_{j}}{\partial r}(0,0)=J'(0)v_{o}\,\neq 0$. From the Implicit Function Theorem we obtain that there exist neighbourhoods $U_{N},U_{\mathbb{R}}$ of $0\in N,0\in \mathbb{R}$ respectively, such that $U_{N}\times U_{\mathbb{R}}\subseteq U'$, and a $C^{1}$ function $j:U_{N}\subseteq N\rightarrow U_{\mathbb{R}}\subseteq \mathbb{R}$ such that $j(0)=0$ and $f_{J}(n,r)=0  \textit{ iff } r=j(n)$, for all $n\in U_{N}, r\in U_{\mathbb{R}}$. Hence $J(n+rv_{o})=0$ \textit{iff } $r=j(n), n\in U_{N}, r \in U_{\mathbb{R}}$.\\
By considering the $C^{1}$ functional $f_{\vartheta}:U'\subseteq N\times \mathbb{R}\rightarrow \mathbb{R},f_{\vartheta}(n,r):=\vartheta(n+rv_{o})$, and arguing in the same way as before, we conclude the existence of neighbourhoods $V_{N},V_{\mathbb{R}}$ of $0\in N,0\in \mathbb{R}$ respectively, such that $V_{N}\times V_{\mathbb{R}}\subseteq U'$, and a $C^{1}$ function $\Theta:V_{N}\subseteq N\rightarrow V_{\mathbb{R}}\subseteq  \mathbb{R}$ such that $\Theta(0)=0$ and  $\vartheta(n+rv_{o})=0 $ \textit{iff} $r=\Theta(n)$, for all $n\in V_{N}, r\in V_{\mathbb{R}}$.\\ 
\indent
Let us consider $u\in \Sigma \cap U :$ then $u=n+rv_{o}$, for suitable $n\in N$ and $r\in \mathbb{R}$. If $u$ is taken near $0$ we may suppose that $n\in U_{N}\cap V_{N}$ and $r\in U_{\mathbb{R}}\cap V_{\mathbb{R}}$. By construction $\vartheta(u)=0$, that is $\vartheta(n+rv_{o})=0$. The last equality implies that $r=\Theta(n)$, i.e. $u=n+\Theta(n)v_{o}$. Moreover $J(n+j(n)v_{o})=0$. Hence, by construction, $n+j(n)v_{o}\in S$. We assumed that $S \subseteq \Sigma$, therefore $n+j(n)v_{o}\in \Sigma$. Thus we obtain that $\vartheta(n+j(n)v_{o})=0$  and, by uniqueness, it follows that $j(n)=\Theta(n)$. Hence $u=n+\Theta(n)v_{o}=n+j(n)v_{o}\in S. \bracevert $\\ \vspace{2pt}
\subsection{Transversality and Singular Strata}\label{ss24}
\quad In the previous section we studied, in Theorem \ref{Teo233}, some useful properties deriving from 1-transversality. Moreover, as we noted in Remark \ref{Rem224}, the Geometric Lemma (i.e. Lemma \ref{Lem223}) allows us to study what is reasonable to call a \textit{stratification of submanifolds}. Actually the hypotheses of the Geometric Lemma are \textit{transversality conditions}, generalizing the condition used to define 1-transverse singularities, which are related to the study of suitable geometric objects naturally associated with the given Fredholm map. As we shall see below, when the assumptions of the Geometric Lemma are satisfied it is possible to conclude that these geometric objects are nested submanifolds which exhibit properties similar to those proved in Theorem \ref{Teo233} for $S_{1}(F)$ near a 1-transverse singularity. The geometric objects we have just referred to are introduced formally in the following
\begin{definition} \label{D241} Let $U,V$ be open subsets in the $B$-spaces $X,Y$ and $F:U\subseteq X\rightarrow V\subseteq Y$ a $C^{d}$ 0-Fredholm map $(d\geq 2)$. For $u_{o}\in S_{1}(F),(\varphi,\psi)\in \mathcal{P}(F,u_{o})$ and $1\leq k\leq d$ we say that the set 
\begin{center}
$S_{1_{k}}(F)(\varphi,\psi):=\{u \text{ near }u_{o}:J_{0}(\varphi,\psi)(u)=\ldots=J_{k-1}(\varphi,\psi)(u)=0\} $
\end{center}
is the \textit{k-th singular stratum} for $F$ (in fact we are again dealing with germs of singular strata).
\end{definition}
\par
To simplify the notations we will sometimes write $S_{1_{k}}(\varphi,\psi)$ instead of $S_{1_{k}}(F)(\varphi,\psi)$ for a fixed map $F$. Moreover, as we will prove in Proposition \ref{Pro245}, if $u_{o}$ is a \textit{k}-transverse singularity for $F$ then the \textit{k}-th singular stratum is independent of the f-pair $(\varphi,\psi)$.\\
\indent We now prove, under the transversality conditions given in (T$_{k}$) of Definition \ref{D211}, a few geometric properties of the \textit{k}-th singular stratum for $F$.\\
\begin{proposition} \label{Pro242} \textit{Let }$U,V$ \textit{be open subsets in the $B$-spaces} $X,Y$\textit{ and }$F:U\subseteq X\rightarrow V\subseteq Y$ \textit{a} $C^{d}$ 0-\textit{Fredholm map} $(d\geq 2)$.\textit{ Let }$u_{o}\in S_{1}(F)$\textit{ and }$(\varphi,\psi)\in \mathcal{P}(F,u_{o})$. \textit{Suppose that} $u_{o}$ \textit{is a} \textit{k}-transverse \textit{singularity  for} $1 \leq k \leq d-1$, i.e.
\begin{center}
$J_{0}(\varphi,\psi)(u_{o})=\ldots=J_{k-1}(\varphi,\psi)(u_{o})=0$ ;\\
$I_{1}(\varphi,\psi)(u_{o}),\ldots,I_{k}(\varphi,\psi)(u_{o})$ \textit{are} l.i. . 
\end{center}
\textit{Then} $S_{1_{k}}(F)(\varphi,\psi)$\textit{ is, near }$u_{o}$,\textit{ a }$C^{d-k}$  \textit{k}-codimensional \textit{submanifold of }$X$\textit{ and,} $\forall\;u\in S_{1_{k}}(F)(\varphi,\psi)$,\textit{ one has that } $T_{u}S_{1_{k}}(F)(\varphi,\psi)=\cap_{j=1}^{k}N(I_{j}(\varphi,\psi)(u))$.
\end{proposition}
\par
\textbf{Proof.} $J_{h}(\varphi,\psi)$ is of class $C^{d-h-1}$, hence $J_{0}(\varphi,\psi),\ldots,J_{k-1}(\varphi,\psi)$ are of class  $C^{d-k}$, with $d-k\geq 1$. By Definition \ref{D116} we know that $J_{h}^{'}(\varphi,\psi)(u_{o})=I_{h+1}(\varphi,\psi)(u_{o})$,\\
$h=0,1,\ldots, k-1$, and the conclusion follows from the Geometric Lemma \ref{Lem223}. $\bracevert$  \\
\par
We actually have some more information: the stratum $S_{1_{h}}(F)(\varphi,\psi)$ is, near $u_{o}$, a one-codimensional submanifold of $S_{1_{h-1}}(F)(\varphi,\psi)$ for $h=1,\ldots,k$ (where $S_{1_{0}}(F)(\varphi,\psi):=X)$. Namely, $\{J_{\eta}^{'}(\varphi,\psi)(u_{o})=I_{\eta+1}(\varphi,\psi)(u_{o}), \eta=0, 1,\ldots,h-1\}$ is  l.i., for  $h=1,\ldots,k$, and thus this is just a special case of the situation described in Remark \ref{Rem224}.\\
We therefore have the following \textit{one-codimensional stratification of manifolds}: 
\begin{center}
$S_{1_{k}}(F)(\varphi,\psi)\subseteq S_{1_{k-1}}(F)(\varphi,\psi) \subseteq\ldots\subseteq S_{1_{2}}(F)(\varphi,\psi)\subseteq S_{1_{1}}(F)(\varphi,\psi)=S_{1}(F)\subseteq X $,
\end{center}
i.e. each manifold has codimension $1$ in the next one. Note that  $S_{1_{1}}(F)(\varphi,\psi)=S_{1}(F)$  because $J_{0}^{-1}(0)=S_{1}(F)$ near $u_{o}$, as proved in Theorem \ref{Teo233}. We can use this result since the linear independence of $I_{1}(\varphi,\psi)(u_{o}),\ldots,I_{k}(\varphi,\psi)(u_{o})$ implies (in particular) that $I_{1}(\varphi,\psi)(u_{o})\,\neq 0$, i.e. $u_{o}$ is a 1-transverse singularity. The above chain of inclusions gives rise to a \textit{stratification of singularities} for $F$ near $u_{o}$. It is important to note that in the following we will also consider other singular strata $S_{1_{h}}(F)(\varphi,\psi)$ for $h=k+1,\ldots,d$ which, however, can possibly be empty sets; in this case the corresponding statements are to be considered void.\\
\par
We now prove that conditions $(\text{T}_k),(\text{S}_k),(\text{M}_k)\text{ and }(\text{T}_\infty)$ of Definition \ref{D211} are independent of the chosen f-pair.\\
\begin{theorem}\label{Teo243}\textit{ Let }$U,V$\textit{ be open subsets in the $B$-spaces} $X,Y$ \textit{and } $F:U\subseteq X\rightarrow V\subseteq Y$\textit{ a }$C^{d}$ 0-\textit{Fredholm map }$(d\geq 2)$.\textit{ Let }$u_{o}\in S_{1}(F)$\textit{ and } $(\varphi,\psi),(\widetilde{\varphi},\widetilde{\psi})\in \mathcal{P}(F,u_{o})$. \textit{Suppose that} $1\leq k\leq d-1$ \textit{and consider the following conditions for the f-pair} $(\varphi,\psi)$:
\begin{enumerate}
\item[i)]$J_{0}(\varphi,\psi)(u_{o})=\ldots=J_{k-1}(\varphi,\psi)(u_{o})=0$ ;\\
$I_{1}(\varphi,\psi)(u_{o}),\ldots,I_{k}(\varphi,\psi)(u_{o})$\textit{ are }l.i. .
\item[ii)]$J_{0}(\varphi,\psi)(u_{o})=\ldots=J_{k-1}(\varphi,\psi)(u_{o})=0, J_{k}(\varphi,\psi)(u_{o})\,\neq 0$ ;\\
$I_{1}(\varphi,\psi)(u_{o}),\ldots,I_{k-1}(\varphi,\psi)(u_{o})$\textit{ are }l.i. .
\item[iii)] $J_{0}(\varphi,\psi)(u_{o})=\ldots=J_{k}(\varphi,\psi)(u_{o})=0$ ;\\
$I_{1}(\varphi,\psi)(u_{o}),\ldots, I_{k}(\varphi,\psi)(u_{o})$ \textit{ are } l.i.;\\
$I_{1}(\varphi,\psi)(u_{o}),\ldots,I_{k+1}(\varphi,\psi)(u_{o})$\textit{  are not  }l.i. .
\item[iv)] \textit{when } $d=\infty$\textit{ or }$d=\omega$ \textit{for all integers} $k\geq 1$ \textit{one has}\\
$J_{0}(\varphi,\psi)(u_{o})=\ldots=J_{k-1}(\varphi,\psi)(u_{o})=0$ ;\\
$I_{1}(\varphi,\psi)(u_{o}),\ldots,I_{k}(\varphi,\psi)(u_{o})$ \textit{ are } l.i..
\end{enumerate}
\textit{Then these conditions are true if and only if the analogous statements relative to the f-pair} $(\widetilde{\varphi}, \widetilde{\psi})$ \textit{are true.}
\end{theorem}
\par For the proof we need two results: the first one is a lemma which says that a \textit{k}-singularity is a \textit{k}-transverse singularity, while the second tool is a technical result whose proof is postponed after the proof of Theorem \ref{Teo243}. \\
\begin{lemma} \label{Lem244} \textit{Under the hypotheses of  Theorem} 2.4.3, \textit{if }
\begin{equation*}
\begin{split}
J_{1}(\varphi,\psi)&(u_{o})=\ldots=J_{k-1}(\varphi,\psi)(u_{o})=0 \; ,\; J_{k}(\varphi,\psi)(u_{o})\,\neq 0 \\
\textit{and}\hspace{110pt}&\\
I_{1}(\varphi,\psi)&(u_{o}),\ldots,I_{k-1}(\varphi,\psi)(u_{o})\textit{ are }\text{ l.i.}\, ,\\
\textit{then}\hspace{110pt}& \\
I_{1}(\varphi,\psi)&(u_{o}),\ldots,I_{k}(\varphi,\psi)(u_{o})\textit{ are }\text{  l.i.} .
\end{split}
\end{equation*}
\end{lemma}
\par \textbf{Proof.}  By definition,  
\begin{equation*}
\begin{split}
J_{h}(\varphi,\psi)&(u_{o})=I_{h}(\varphi,\psi)(u_{o})\varphi(u_{o})=0 , \;\; h=1,\ldots,k-1, \\
J_{k}(\varphi,\psi)&(u_{o})=I_{k}(\varphi,\psi)(u_{o})\varphi(u_{o})\, \neq 0 . 
\end{split}
\end{equation*}
Hence $\varphi(u_{o})\in \cap_{h=1}^{k-1}N(I_{h}(\varphi,\psi)(u_{o}))$ and $\varphi(u_{o})\notin N(I_{k}(\varphi,\psi)(u_{o}))$. Since $I_{1}(\varphi,\psi)(u_{o})$,\\
$\ldots,I_{k-1}(\varphi,\psi)(u_{o})$ are l.i. then, by the Algebraic Lemma \ref{Lem221}, $f) \Rightarrow a)$, we have that $I_{1}(\varphi,\psi)(u_{o}),\ldots,I_{k}(\varphi,\psi)(u_{o})$ are l.i. . $\bracevert$  \\
\par
The above-mentioned technical result is given by the following\\
\begin{proposition} \label{Pro245}\textit{ Let }$U,V$ \textit{be open subsets in the $B$-spaces} $X,Y$\textit{ and }$F:U \subseteq X\rightarrow V\subseteq Y$\textit{ a }$C^{d}$ 0\textit{-Fredholm map} $(d\geq 2)$. \textit{Let }$u_{o}\in S_{1}(F),(\varphi,\psi)\in \mathcal{P}(F,u_{o})$ \textit{and let us assume that} $u_{o}$\textit{ is a }\textit{k}-transverse \textit{singularity for} $1\leq k\leq d-1$, \textit{i.e.} 
\begin{align*}
J_{0}(\varphi,\psi)(u_{o})=\ldots=J_{k-1}(\varphi,\psi)(u_{o})=0\, ; \\
I_{1}(\varphi,\psi)(u_{o}),\ldots,I_{k}(\varphi,\psi)(u_{o})\textit{ are} \text{ l.i.}.
\end{align*} 
\textit{If} $(\widetilde{\varphi},\widetilde{\psi})\in \mathcal{P}(F,u_{o})$ \textit{is another fibering pair then, near }$u_{o}$, \textit{there exist two non-vanishing real-valued functions }$\alpha,\beta$ \textit{of class} $C^{d-1}$ \textit{such that the following relations hold}:
\begin{align*}
\text{(j)}_{k}\qquad & \begin{cases}
I_{h+1}(\widetilde{\varphi},\widetilde{\psi})(u)v=\alpha^{h+1}(u)\beta(u)I_{h+1}(\varphi,\psi)(u)v\\
\textit{for  } u\in S_{1_{h+1}}(F)(\varphi,\psi),\; v\in T_{u}S_{1_{h}}(F)(\varphi,\psi), h=0, 1 ,\ldots,k
\end{cases}\\
\\
\text{(jj)}_{k}\qquad & \begin{cases}
\cap_{j=1}^{h+1}N(I_{j}(\widetilde{\varphi},\widetilde{\psi})(u))=\cap_{j=1}^{h+1}N(I_{j}(\varphi,\psi)(u))\\
\textit{for } u\in S_{1_{h+1}}(F)(\varphi,\psi),h=0,1,\ldots,k
\end{cases}\\
\\
\text{(jjj)}_{k}\qquad & \begin{cases}
J_{h+1}(\widetilde{\varphi},\widetilde{\psi})(u)=\alpha^{h+2}(u)\beta(u)J_{h+1}(\varphi,\psi)(u)\\
\textit{for } u\in S_{1_{h+1}}(F)(\varphi,\psi), h=0,1,\ldots,k
\end{cases}\\
\\
\text{(jv)}_{k}\qquad &  S_{1_{h+1}}(F)(\widetilde{\varphi},\widetilde{\psi})=S_{1_{h+1}}(F)(\varphi,\psi), h=0,1,\ldots, k+1 
\end{align*}
\\
\textit{(for }$k=d-1$,\textit{ conditions }\text{(j)}$_{k}$,  \text{(jj})$_{k}$, \text{(jjj)}$_{k}$\textit{ are valid only for} $h=0,1,\ldots,k-1$, \textit{and condition }\text{(jv)}$_{k}$ \textit{ for } $h=0,1,\ldots,k$\textit{)}.
\end{proposition}
\par
\textbf{Proof of Theorem 2.4.3.} By interchanging the roles of the f-pairs $(\varphi,\psi)$ and $(\widetilde{\varphi},\widetilde{\psi})$, it suffices to show that if the statements i),...,iv) for $(\varphi,\psi)$ are true then the analogous statements for $(\widetilde{\varphi},\widetilde{\psi})$ are true.\\
\indent Let us suppose that conditions i), ii) and iii) relative to the pair $(\varphi,\psi)$ are satisfied. Thanks to the above lemma, in all cases  i), ii) and iii) we can conclude that
\begin{center}
$J_{1}(\varphi,\psi)(u_{o})=\ldots=J_{k-1}(\varphi,\psi)(u_{o})=0$ ;\\
$I_{1}(\varphi,\psi)(u_{o}),\ldots,I_{k}(\varphi,\psi)(u_{o})$ are l.i. ,
\end{center}  
i.e. $u_{o}$ is a \textit{k}-transverse singularity for $F$. Hence relations (j)$_{k}$, (jj)$_{k}$, (jjj)$_{k}$, (jv)$_{k}$ of Proposition \ref{Pro245} are true for all cases i), ii) and iii). In particular, if we consider (jv)$_{k}$ for $h=k-1$ and $h=k$ it follows that  $S_{1_{k}}(F)(\widetilde{\varphi},\widetilde{\psi})=S_{1_{k}}(F)(\varphi,\psi)$ and $S_{1_{k+1}}(F)(\widetilde{\varphi},\widetilde{\psi})=S_{1_{k+1}}(F)(\varphi,\psi)$ near $u_{o}$. Hence, from Definition \ref{D241}, we obtain the following two facts: since $J_{0}(\varphi,\psi)(u_{o})=\ldots=J_{k-1}(\varphi,\psi)(u_{o})=0$  then \vspace{-10 pt}\\

\hspace{30 pt}a)$\qquad J_{0}(\widetilde{\varphi},\widetilde{\psi})(u_{o})=\ldots=J_{k-1}(\widetilde{\varphi},\widetilde{\psi})(u_{o})=0$ \\
and  \vspace{-10 pt} \\ 

\hspace{30 pt}b)$\qquad J_{k}(\varphi,\psi)(u_{o})=0$ if and only if $J_{k}(\widetilde{\varphi},\widetilde{\psi}(u_{o})=0$ .\vspace{-10 pt} \\

\noindent On the other hand, we can also prove that \vspace{-10 pt} \\

\hspace{30pt}c)$\qquad I_{1}(\widetilde{\varphi},\widetilde{\psi})(u_{o}),\ldots,I_{k}(\widetilde{\varphi},\widetilde{\psi})(u_{o})$ are  l.i. . \vspace{-10 pt}\\

\noindent Thus we obtain that \vspace{-10 pt}\\

\hspace{30pt}d)$\qquad I_{1}(\widetilde{\varphi},\widetilde{\psi})(u_{o}),\ldots,I_{k-1}(\widetilde{\varphi},\widetilde{\psi})(u_{o})$  are  l.i. .\vspace{-10 pt}\\

\noindent It is clear how the points  a), b), c), d)  prove the statements  i) and ii) relative to the pair $(\widetilde{\varphi},\widetilde{\psi})$.\\
\indent Let us prove claim c). As usual we set $S_{1_{0}}(F)(\varphi,\psi):=X,N(I_{0}(\varphi,\psi)(u_{o})):= X,N(I_{0}(\widetilde{\varphi},\widetilde{\psi})(u_{o})):=X$. Since $J_{0}(\varphi,\psi)(u_{o})=\ldots= J_{k-1}(\varphi,\psi)(u_{o})=0$ then $u_{o}\in S_{1_{\eta}}(F)(\varphi,\psi)$, for $\eta=1,\ldots, k$. Moreover,  $I_{1}(\varphi,\psi)(u_{o}),\ldots,I_{\eta}(\varphi,\psi)(u_{o})$ are l.i. for $\eta=1,\ldots,k.$ 
Then, from Proposition \ref{Pro242}  for the integers $\eta=1,\ldots,k$, one has that  \vspace{-18 pt} \\
\begin{center}
$T_{u_{o}}S_{1_{\eta}}(F)(\varphi,\psi)=\cap_{j=1}^{\eta}N(I_{j}(\varphi,\psi)(u_{o})),\;\,\eta=1,\ldots,k$.
\end{center}
Given that $u_{o}\in S_{1_{\eta}}(F)(\varphi,\psi),\eta=1,\ldots,k$, if we set $h+1:=\eta$ in (jj)$_{k}$ of Proposition \ref{Pro245} we find \vspace{-18 pt} \\
\begin{center}
$\cap_{j=1}^{\eta}N(I_{j}(\widetilde{\varphi},\widetilde{\psi})(u_{o}))=\cap_{j=1}^{\eta}N(I_{j}(\varphi,\psi)(u_{o})),\eta=1,\ldots,k$.
\end{center}
Hence, for  $h=0,\ldots,k-1 $, we can write that \vspace{-18 pt} \\
\begin{center}
$T_{u_{o}}S_{1_{h}}(F)(\varphi,\psi)=\cap_{j=0}^{h}N(I_{j}(\varphi,\psi)(u_{o}))=\cap_{j=0}^{h}N(I_{j}(\widetilde{\varphi},\widetilde{\psi})(u_{o}))$.
\end{center}
Since  $u_{o}\in S_{1_{h+1}}(F)(\varphi,\psi),h=0,\ldots,k-1$, we may use (j)$_{k}$ of  Proposition 2.4.5  to obtain
\begin{equation}\label{241}
\begin{split}
\qquad &\quad \; I_{h+1}(\widetilde{\varphi},\widetilde{\psi})(u_{o}) \lvert\cap_{j=0}^{h}N(I_{j}(\widetilde{\varphi},\widetilde{\psi})(u_{o})) =\\
&= I_{h+1}(\widetilde{\varphi},\widetilde{\psi})(u_{o})\lvert T_{u_{o}}S_{1_{h}}(F)(\varphi,\psi)=\\
&= \alpha^{h+1}(u_{o})\beta(u_{o})I_{h+1}(\varphi,\psi)(u_{o})\lvert T_{u_{o}}S_{1_{h}}(F)(\varphi,\psi)=\\
&= \alpha^{h+1}(u_{o})\beta(u_{o})I_{h+1}(\varphi,\psi)(u_{o})\lvert \cap_{j=0}^{h}N(I_{j}(\varphi,\psi)(u_{o})), \;h=0,\ldots,k-1,
\end{split}
\end{equation}
where the vertical bar $\lvert$ means \textquotedblleft restricted to\textquotedblright. We recall that $\alpha,\beta$ are $C^{d-1}$ real functions defined near $u_{o}$ such that $\alpha(u_{o})\neq 0$ and $\beta(u_{o})\neq 0$.\\
Since $I_{1}(\varphi,\psi)(u_{o}),\ldots,I_{k}(\varphi,\psi)(u_{o})$ are l.i. then, by the Algebraic Lemma \ref{Lem221}, a) $\Rightarrow$ e), we have that
\begin{center}
$I_{h+1}(\varphi,\psi)(u_{o})\lvert  \cap_{j=0}^{h}N(I_{j}(\varphi,\psi)(u_{o}))\,\neq 0, \;h=0,\ldots,k-1$. 
\end{center} 
Hence (\ref{241}) implies that 
\begin{center}
$I_{h+1}(\widetilde{\varphi},\widetilde{\psi})(u_{o})\lvert \cap_{j=0}^{h}N(I_{j}(\widetilde{\varphi},\widetilde{\psi})(u_{o}))\,\neq 0, \, h=0,\ldots,k-1$.
\end{center}
By using again the Algebraic Lemma, e) $\Rightarrow$ a), we see that the functionals  $I_{1}(\widetilde{\varphi},\widetilde{\psi})(u_{o}),$\\
$\ldots,I_{k}(\widetilde{\varphi},\widetilde{\psi})(u_{o})$ are l.i.. Hence point c) is true.\\
\indent
Now we continue by proving condition iii) for the f-pair $(\widetilde{\varphi},\widetilde{\psi})$. Let us suppose that
\begin{equation*}
\begin{split}
\qquad \qquad \quad &J_{0}(\varphi,\psi)(u_{o})=\ldots,J_{k}(\varphi,\psi)(u_{o})=0\,;\\
&I_{1}(\varphi,\psi)(u_{o}),\ldots,I_{k}(\varphi,\psi)(u_{o}) \text{ are l.i.},\\
&I_{1}(\varphi,\psi)(u_{o}),\ldots,I_{k+1}(\varphi,\psi)(u_{o}) \text{ are not l.i.}.
\end{split}
\end{equation*}
By virtue of the above points  a), b)  and  c), we only have to show that 
\begin{center}
$I_{1}(\widetilde{\varphi},\widetilde{\psi})(u_{o}),\ldots,I_{k+1}(\widetilde{\varphi},\widetilde{\psi})(u_{o}) \text{ are not l.i.}$.
\end{center}
Since $u_{o}\in S_{1_{h+1}}(F)(\varphi,\psi),h=0,\ldots,k$, one can show, in a similar way as for (\ref{241}), the validity of the following formula:

\begin{equation}\label{242}
\begin{split}
\qquad \quad I_{k+1}&(\widetilde{\varphi},\widetilde{\psi})(u_{o}) \lvert\cap_{j=0}^{k}N(I_{j}(\widetilde{\varphi},\widetilde{\psi})(u_{o})) =\\
&=\alpha^{k+1}(u_{o})\beta(u_{o})I_{k+1}(\varphi,\psi)(u_{o})\lvert \cap_{j=0}^{k}N(I_{j}(\varphi,\psi)(u_{o})).
\end{split}
\end{equation}

We contend that $I_{k+1}(\varphi,\psi)(u_{o})\lvert \cap_{j=0}^{k}N(I_{j}(\varphi,\psi)(u_{o}))= 0.$ \\
In fact, by contradiction let us suppose that $I_{k+1}(\varphi,\psi)(u_{o})\lvert \cap_{j=0}^{k}N(I_{j}(\varphi,\psi)(u_{o}))\, \neq 0$. Since  $I_{1}(\varphi,\psi)(u_{o}),\ldots,I_{k}(\varphi,\psi)(u_{o})$ are l.i. then $I_{1}(\varphi,\psi)(u_{o}),\ldots, I_{k+1}(\varphi,\psi)(u_{o})$ are l.i. by the Algebraic Lemma \ref{Lem221}, f) $\Rightarrow$ a). Yet the linear independence of $I_{1}(\varphi,\psi)(u_{o})$,\\
$\ldots, I_{k+1}(\varphi,\psi)(u_{o})$ contradicts the initial assumption. \\
By inspection of the formula (\ref{242}) and recalling that $\alpha(u_{o})\,\neq 0,\beta(u_{o})\,\neq 0$, we get that $I_{k+1}(\widetilde{\varphi},\widetilde{\psi})(u_{o}) \lvert\cap_{j=0}^{k}N(I_{j}(\widetilde{\varphi},\widetilde{\psi})(u_{o}))=0$. This implies (by the Algebraic Lemma, a) $\Rightarrow$ f)) that $I_{1}(\widetilde{\varphi},\widetilde{\psi})(u_{o}),\ldots, I_{k+1}(\widetilde{\varphi},\widetilde{\psi})(u_{o})$ are not l.i. .\\
\indent
Finally, let us assume that statement iv) holds for the pair $(\varphi,\psi)$. Thus i) is true for $(\varphi,\psi)$ for any integer $k$: thanks to what we have just shown, this implies that i) holds for the pair $(\widetilde{\varphi},\widetilde{\psi})$ for any $k$. Therefore statement iv) is also true for the pair $(\widetilde{\varphi},\widetilde{\psi}). \; \bracevert$\\
\par
\textbf{Proof of Proposition 2.4.5.} We argue by induction on $k$. \\
By convention we set  $S_{1_{0}}(F)(\varphi,\psi)\equiv S_{1_{0}}(\varphi,\psi)\!:=\!X$ and  $T_{u}S_{1_{0}}(\varphi,\psi)\!:=\!X,\;\forall \,u\in X$. In order to prove the Proposition we first establish some formulas valid on a suitable neighbourhood of $u_{o}$. We recall that, near $u_{o}$, we proved the following equality (see (\ref{232}) in the proof of  Proposition \ref{Pro231}):
\begin{equation}\label{243}
I_{1}(\widetilde{\varphi},\widetilde{\psi})(u)v=\alpha(u)\beta(u)I_{1}(\varphi,\psi)(u)v\; ,\; u\in S_{1}(F)\, , \, \forall \, v\in X 
\end{equation}
for suitable $C^{d-1}$ functions $\alpha,\beta$, defined near $u_{o}$, which are always different from zero and such that $\widetilde{\varphi}=\alpha\varphi,\widetilde{\psi}=\beta\psi$ on $S_{1}(F)\equiv S_{1}$.\\
Hence, given that $\alpha,\beta$ are non-vanishing functions, it follows that
\begin{equation}\label{244}
N(I_{1}(\widetilde{\varphi},\widetilde{\psi})(u))= N(I_{1}(\varphi,\psi)(u)),\; u\in S_{1}.
\end{equation}
Moreover, from (\ref{243}) we obtain 

\begin{equation}\label{245}
\begin{split}
J_{1}(\widetilde{\varphi},\widetilde{\psi})(u)&= I_{1}(\widetilde{\varphi},\widetilde{\psi})(u))\widetilde{\varphi}(u)=\alpha(u)\beta(u)I_{1}(\varphi,\psi)(u)\widetilde{\varphi}(u)=\\
&=\alpha^{2}(u)\beta(u)I_{1}(\varphi,\psi)(u)\varphi(u)=\alpha^{2}(u)\beta(u)J_{1}(\varphi,\psi)(u)\, ,\, u\in S_{1}.
\end{split}
\end{equation}
\indent
Now let us assume that $u_{o}$ is a 1-transverse singularity for F, i.e. $k=1$ and $I_{1}(\varphi,\psi)(u_{o})$ is l.i. (while $J_{0}(\varphi,\psi)(u_{o})=0$ is always satisfied). By using Theorem \ref{Teo233} we have that the singular set $S_{1}$ is, near $u_{o}$, a $C^{d-1}$ one-codimensional submanifold of $X$ and
\begin{equation}\label{246}
T_{u}S_{1}=N(I_{1}(\widetilde{\varphi},\widetilde{\psi})(u))=N(I_{1}(\varphi,\psi)(u)),\, u\in S_{1}.
\end{equation}
Moreover, again from Theorem \ref{Teo233}, it also follows that
\begin{equation}\label{247}
S_{1_{1}}(\widetilde{\varphi},\widetilde{\psi})=S_{1_{1}}(\varphi,\psi)=S_{1}.
\end{equation}
By definition $ S_{1_{2}}(\widetilde{\varphi},\widetilde{\psi})=\{u\in S_{1_{1}}(\widetilde{\varphi},\widetilde{\psi})=S_{1}:J_{1}(\widetilde{\varphi},\widetilde{\psi})(u)=0\}$; a similar equality defines $S_{1_{2}}(\varphi,\psi) $. From (\ref{245}) and by taking into account that $\alpha \,\neq 0$ and $\beta\, \neq 0$  we obtain the equality
\begin{equation}\label{248}
S_{1_{2}}(\widetilde{\varphi},\widetilde{\psi})=S_{1_{2}}(\varphi,\psi).
\end{equation}
By differentiating the identity (\ref{245}) on $S_{1}$ we obtain that, for $u\in S_{1}, v\in T_{u}S_{1}$,
\begin{equation*}
\begin{split}
I_{2}(\widetilde{\varphi},\widetilde{\psi})(u)v&=J_{1}^{'}(\widetilde{\varphi},\widetilde{\psi})(u)v=\\
&=(\alpha^{2}\beta)'(u)v\cdotp J_{1}(\varphi,\psi)(u)+\alpha^{2}(u)\beta(u)J_{1}^{'}(\varphi,\psi)(u)v=\\
&=(\alpha^{2}\beta)'(u)v\cdotp J_{1}(\varphi,\psi)(u)+\alpha^{2}(u)\beta(u)I_{2}(\varphi,\psi)(u)v.
\end{split}
\end{equation*}
By definition we have that $S_{1_{2}}(\varphi,\psi) \subseteq S_{1_{1}}(\varphi,\psi)=S_{1}$ and that $J_{1}(\varphi,\psi)(u)=0$ on $S_{1_{2}}(\varphi,\psi)$. Hence the last formula gives
\begin{equation}\label{249}
I_{2}(\widetilde{\varphi},\widetilde{\psi})(u)v=\alpha^{2}(u)\beta(u)I_{2}(\varphi,\psi)(u)v, \, u\in S_{1_{2}}(\varphi,\psi), \, v\in T_{u}S_{1}.
\end{equation}
Note that formulas (\ref{243}) and (\ref{249}) give (j)$_{1}$.\\
From  (\ref{246}), (\ref{249}) and since $\alpha\neq 0, \beta\neq 0$ it follows that
\begin{equation}\label{2410}
\begin{split}
&N(I_{1}(\widetilde{\varphi},\widetilde{\psi})(u))\cap N(I_{2}(\widetilde{\varphi},\widetilde{\psi})(u))=\\
&=N(I_{1}(\varphi,\psi)(u))\cap N(I_{2}(\varphi,\psi)(u)), u\in S_{1_{2}}(\varphi,\psi).
\end{split}
\end{equation}
Thanks to equalities (\ref{244}), (\ref{2410}) we obtain (jj)$_{1}$.\\
Let $u\in S_{1_{2}}(\widetilde{\varphi},\widetilde{\psi})=S_{1_{2}}(\varphi,\psi)\subseteq S_{1}$. In particular one has
\begin{center}
$J_{1}(\widetilde{\varphi},\widetilde{\psi})(u)=0=J_{1}(\varphi,\psi)(u)$\,,
\end{center} 
i.e.
\begin{center}
$I_{1}(\widetilde{\varphi},\widetilde{\psi})(u)\widetilde{\varphi}(u)=0=I_{1}(\varphi,\psi)(u)\varphi(u)$.
\end{center}
Hence
\begin{center}
$\widetilde{\varphi}(u)\in N(I_{1}(\widetilde{\varphi},\widetilde{\psi})(u))\;\, , \;\, \varphi(u)\in N(I_{1}(\varphi,\psi)(u))$.
\end{center}
This means, from (\ref{246}), that
\begin{center}
$\widetilde{\varphi}(u),\varphi(u)\in T_{u}S_{1}$.
\end{center}
By virtue of (\ref{249}) and by definition of $J_{2}$ one has that
\begin{equation}\label{2411}
\begin{split}
J_{2}(\widetilde{\varphi},\widetilde{\psi})(u)&=I_{2}(\widetilde{\varphi},\widetilde{\psi})(u)\widetilde{\varphi}(u)=\alpha^{2}(u)\beta(u)I_{2}(\varphi,\psi)(u)\widetilde{\varphi}(u)=\\
&=\alpha^{3}(u)\beta(u)I_{2}(\varphi,\psi)(u)\varphi(u)=\alpha^{3}(u)\beta(u)J_{2}(\varphi,\psi)(u)
\end{split}
\end{equation}
for $u\in S_{1_{2}}(\widetilde{\varphi},\widetilde{\psi})=S_{1_{2}}(\varphi,\psi)$.\\
Note that (\ref{245}) and (\ref{2411}) imply (jjj)$_{1}$.\\
\indent Since $\alpha,\beta$ are different from zero, we get 
\begin{equation}\label{2412}
\begin{split}
S_{1_{3}}(\widetilde{\varphi},\widetilde{\psi})&= \{u \in S_{1_{2}}(\widetilde{\varphi},\widetilde{\psi}):J_{2}(\widetilde{\varphi},\widetilde{\psi})(u)=0\}=\\
&=\{u \in S_{1_{2}}(\varphi,\psi):J_{2}(\widetilde{\varphi},\widetilde{\psi})(u)=0\}=\\
&=\{u \in S_{1_{2}}(\varphi,\psi):J_{2}(\varphi,\psi)(u)=0\}=S_{1_{3}}(\varphi,\psi).
\end{split}
\end{equation}
Finally, by recalling that $S_{1_{2}}(\varphi,\psi)\subseteq S_{1_{1}}(\varphi,\psi)=S_{1}$ we have that  equalities (\ref{247}), (\ref{248}) and (\ref{2412}) prove (jv)$_{1}$. Hence we have proved the thesis for $k=1$. \\
\indent Now we assume the result to be true for $k=1,2,\ldots,d-2$. We shall show that the thesis holds for the integer $k+1$. The proof is similar to that given for $k=1$.\\
\indent Let us suppose that $u_{o}$ is a $(k+1)$-transverse singularity for $F$, i.e.\\
\begin{equation*}
\begin{split}
J_{0}(\varphi,\psi)(u_{o})=&\ldots=J_{k}(\varphi,\psi)(u_{o})= 0\,\, ;\\
I_{1}(\varphi,\psi)(u_{o}),\ldots&,I_{k+1}(\varphi,\psi)(u_{o})\text{  are l.i. .}
\end{split}
\end{equation*}
In particular $J_{0}(\varphi,\psi)(u_{o})=\ldots =J_{k-1}(\varphi,\psi)(u_{o})=0$ and moreover $I_{1}(\varphi,\psi)(u_{o}),\ldots$,\\
$I_{k}(\varphi,\psi)(u_{o})$ are l.i.. Then, by the inductive hypothesis, it follows that statements (j)$_{k}$, (jj)$_{k}$, (jjj)$_{k}$, (jv)$_{k}$  are true.\\
\indent Since $u_{o}$ is $(k+1)$-transverse then, thanks to Proposition \ref{Pro242}, we have that $S_{1_{k+1}}(\varphi,\psi)$ is, near $u_{o}$, a submanifold such that $T_{u}S_{1_{k+1}}(\varphi,\psi)=\cap_{j=1}^{k+1} N(I_{j}(\varphi,\psi)(u))$, for $u \in S_{1_{k+1}}(\varphi,\psi) $.\\
\indent By using relation (jjj)$_{k}$ for $h=k$ one has  
\begin{center}
$J_{k+1}(\widetilde{\varphi},\widetilde{\psi})(u)=\alpha^{k+2}(u)\beta(u)J_{k+1}(\varphi,\psi)(u) \,, \; u \in S_{1_{k+1}}(\varphi,\psi).$
\end{center}
By differentiation of this equality on $S_{1_{k+1}}(\varphi,\psi)$ we obtain that
\begin{center}
$I_{k+2}(\widetilde{\varphi},\widetilde{\psi})(u)v=(\alpha^{k+2}\beta)'(u)v\cdotp J_{k+1}(\varphi,\psi)(u)+\alpha^{k+2}(u)\beta(u)I_{k+2}(\varphi,\psi)(u)v$ 
\end{center}
for $u\in S_{1_{k+1}}(\varphi,\psi)\, , v\in T_{u}S_{1_{k+1}}(\varphi,\psi)$.\\
Since, by definition, $S_{1_{k+2}}(\varphi,\psi)=\{u\in S_{1_{k+1}}(\varphi,\psi):J_{k+1}(\varphi,\psi)(u)=0\} $  we get the identity
\begin{equation}\label{2413}
I_{k+2}(\widetilde{\varphi},\widetilde{\psi})(u)v = \alpha^{k+2}(u)\beta(u)I_{k+2}(\varphi,\psi)(u)v \,,
\end{equation}
for $u\in S_{1_{k+2}}(\varphi,\psi)$ and $v\in T_{u}S_{1_{k+1}}(\varphi,\psi)$.\\
\indent From relation (jj)$_{k}$, used for $h=k$, one has that
\begin{equation}\label{2414}
\begin{split}
\cap_{j=1}^{k+1} N(I_{j}(\widetilde{\varphi},\widetilde{\psi})(u))&=\cap_{j=1}^{k+1}N(I_{j}(\varphi,\psi)(u))=\\
&=T_{u}S_{1_{k+1}}(\varphi,\psi) \,, u\in S_{1_{k+1}}(\varphi,\psi).
\end{split}
\end{equation}
Since $S_{1_{k+2}}(\varphi,\psi)\subseteq S_{1_{k+1}}(\varphi,\psi)$ and $\alpha,\beta$ are not zero relations (\ref{2413}), (\ref{2414}) yield
\begin{equation}\label{2415}
\cap_{j=1}^{k+2} N(I_{j}(\widetilde{\varphi},\widetilde{\psi})(u))=\cap_{j=1}^{k+2} N(I_{j}(\varphi,\psi)(u))\, , u\in S_{1_{k+2}}(\varphi,\psi) .
\end{equation}
By using relation (jv)$_{k}$ for  $h=k+1$ one has that $S_{1_{k+2}}(\widetilde{\varphi},\widetilde{\psi})= S_{1_{k+2}}(\varphi,\psi)$.
Hence, for $u \in S_{1_{k+2}}(\widetilde{\varphi},\widetilde{\psi})=S_{1_{k+2}}(\varphi,\psi)$, we have by the very definition of the strata $S_{1_{h}}$ that  
\begin{center}
$J_{h}(\widetilde{\varphi},\widetilde{\psi})(u)=0=J_{h}(\varphi,\psi)(u)$, for $h=0,1,\ldots,k+1$.
\end{center}
Specifically, for $h=1,\ldots,k+1$, we have from the definition of $J_{h}$ that  
\begin{center}
$I_{h}(\widetilde{\varphi},\widetilde{\psi})(u)\widetilde{\varphi}(u)=0=I_{h}(\varphi,\psi)(u)\varphi(u)$, i.e \vspace{5 pt}\\
$\widetilde{\varphi}(u)\in \cap_{j=1}^{k+1}N(I_{j}(\widetilde{\varphi},\widetilde{\psi})(u)), \, \varphi(u)\in \cap_{j=1}^{k+1}N(I_{j}(\varphi,\psi)(u))$.
\end{center}
From (\ref{2414}) we then get that $\widetilde{\varphi}(u), \varphi(u)\in T_{u}S_{1_{k+1}}(\varphi,\psi)$ for all $u\in S_{1_{k+2}}(\widetilde{\varphi},\widetilde{\psi})=S_{1_{k+2}}\varphi,\psi) $.\\
From (\ref{2413}) it follows that
\begin{equation}\label{2416}
\begin{split}
J_{k+2}(\widetilde{\varphi},\widetilde{\psi})(u)&=I_{k+2}(\widetilde{\varphi},\widetilde{\psi})(u)\widetilde{\varphi}(u)=\alpha^{k+2}(u)\beta(u)I_{k+2}(\varphi,\psi)(u)\widetilde{\varphi}(u)=\\
&=\alpha^{k+2}(u)\beta(u)I_{k+2}(\varphi,\psi)(u)\alpha(u)\varphi(u) =\\
&=\alpha^{k+3}(u)\beta(u)I_{k+2}(\varphi,\psi)(u)\varphi(u) = \alpha^{k+3}(u)\beta(u)J_{k+2}(\varphi,\psi)(u) ,\\
&\text{ for }u \in  S_{1_{k+2}}(\widetilde{\varphi},\widetilde{\psi})= S_{1_{k+2}}(\varphi,\psi).
\end{split}
\end{equation}
\indent
By virtue of the above relation, and since $\alpha \neq 0,\beta \neq 0 $, we obtain that
\begin{equation}\label{2417}
\begin{split}
S_{1_{k+3}}(\widetilde{\varphi},\widetilde{\psi})&=\{u \in S_{1_{k+2}}(\widetilde{\varphi},\widetilde{\psi}): J_{k+2}(\widetilde{\varphi},\widetilde{\psi})(u)=0\} =\\
&= \{u \in S_{1_{k+2}}(\varphi,\psi): J_{k+2}(\widetilde{\varphi},\widetilde{\psi})(u)=0\} =\\
&= \{u \in S_{1_{k+2}}(\varphi,\psi): J_{k+2}(\varphi,\psi)(u)=0\} = S_{1_{k+3}}(\varphi,\psi).
\end{split}
\end{equation}
\indent We can now conclude the proof. Namely, statement (j)$_{k}$  and formula (\ref{2413})
yield (j)$_{k+1}$, relation (jj)$_{k}$ and equality (\ref{2415}) give (jj)$_{k+1}$, condition (jjj)$_{k}$ and formula (\ref{2416}) prove (jjj)$_{k+1}$ and, finally, statement (jv)$_{k}$ and equality (\ref{2417}) yield (jv)$_{k+1}. \bracevert$\\
\begin{remark} \label{Rem246} It may be useful to summarize the main results obtained in this section about the stratification of singularities near a simple singularity $u_{o}$ for a $C^{d}$ map $F$.\\
\indent a)  Since, by Definition \ref{D241}, the singular strata for $F$ are the sets
\begin{center}
$S_{1_{h}}(F)(\varphi,\psi):=\{u \text{ near }u_{o}:J_{0}(\varphi,\psi)(u)=\ldots=J_{h-1}(\varphi,\psi)(u)=0\}$,
\end{center}
where $1\leq h\leq d$ and $(\varphi,\psi)\in \mathcal{P}(F,u_{o})$, by construction we always have the  inclusions of sets 
\begin{center}
$S_{1_{d}}(F)(\varphi,\psi)\subseteq S_{1_{d-1}}(F)(\varphi,\psi) \subseteq ...\subseteq S_{1_{2}}(F)(\varphi,\psi)\subseteq S_{1_{1}}(F)(\varphi,\psi)\subseteq X$.
\end{center}
\qquad b)  The ``zero-set'' condition
\begin{center}
$(\text{b}_{1}) \qquad  J_{0}(\varphi,\psi)(u_{o})=\ldots =J_{k-1}(\varphi,\psi)(u_{o})=0$
\end{center}
and the ``independence'' condition
\begin{center}
$(\text{b}_{2}) \qquad \;\;\;\;  I_{1}(\varphi,\psi)(u_{o})=\ldots =I_{k}(\varphi,\psi)(u_{o}) \text{ l.i.}$
\end{center}
do not depend on the given f-pair $(\varphi,\psi)$, for $1\leq k\leq d-1$ (cf. statement  i)  of Theorem \ref{Teo243}). Hence the definition of $k$-transverse singularity, given in (T$_{k}$) of Definition \ref{D211}, is well-posed.\\
\indent c)  When $u_{o}$ is a $k$-transverse singularity, i.e. conditions (b$_{1}$) and (b$_{2}$) hold, we know that $S_{1_{h}}(F)(\varphi,\psi)=S_{1_{h}}(F)(\widetilde{\varphi},\widetilde{\psi})$ near $u_{o}$, for $h=1, 2,\ldots, k+2$ and $(\varphi,\psi),(\widetilde{\varphi},\widetilde{\psi})\in \mathcal{P}(F,u_{o})$, as proved in (jv)$_{k}$ of Proposition \ref{Pro245} (note that, for $k=d-1$, the equality is true for $h=1, 2,\ldots,k+1)$. On the other hand (b$_{2}$) guarantees that $u_{o}$ is 1-transverse for $F$ and so, from Theorem \ref{Teo233}, we have that $S_{1_{1}}(F)(\varphi,\psi)=S_{1}(F)$ near $u_{o}$, for $(\varphi,\psi) \in \mathcal{P}(F,u_{o})$. Hence $S_{1_{1}}(F)(\varphi,\psi) = S_{1_{1}}(F)(\widetilde{\varphi},\widetilde{\psi})= S_{1}(F)$ on a suitable neighbourhood of $u_{o}$. Since $S_{1}(F)$ is a globally defined object, the previous equalities can be used to ease the notation, though with a slight abuse, and to write $ S_{1_{1}}(F)$ instead of $S_{1_{1}}(F)(\varphi,\psi)$. As we shall prove in the next section (cf. Theorem \ref{Teo254}), under conditions (b$_{1}$) and (b$_{2}$) similar identities between local and global objects are also valid for $S_{1_{h}}(F)(\varphi,\psi)$. Hence, to simplify notations, from now on we will simply write $S_{1_{h}}(F)$ instead of $S_{1_{h}}(F)(\varphi,\psi)$, for  $h=1,2,\ldots,k+2$ (or, when $k=d-1$, for $h=1, 2,\ldots,k+1)$. We thus have the following inclusions of sets:
\begin{center}
$S_{1_{k+2}}(F)\subseteq S_{1_{k+1}}(F)\subseteq \ldots \subseteq S_{1_{2}}(F)\subseteq S_{1_{1}}(F)=S_{1}(F)\subseteq X $
\end{center} 
or
\begin{center}
$S_{1_{k+1}}(F)\subseteq \ldots \subseteq S_{1_{2}}(F)\subseteq S_{1_{1}}(F)=S_{1}(F)\subseteq X$ for $k=d-1$.
\end{center}
We point out that $S_{1_{k+2}}(F)$ and $S_{1_{k+1}}(F)$ can possibly be empty sets. Precisely, it is not difficult to show that, near $u_{o}, S_{1_{k+1}}(F)$ is empty \textit{iff} $u_{o}$ is a $k$-singularity; an analogous statement holds for $S_{1_{k+2}}(F)$.\\
\indent d)  Again under conditions (b$_{1}$) and (b$_{2}$), for $h=1,\ldots,k, \text{ the stratum } S_{1_{h}}(F)$ is a one-codimensional submanifold of $S_{1_{h-1}}(F)$ near $u_{o}$, as noted in the comment following Proposition \ref{Pro242}.  Thus we have the one-codimensional stratification of manifolds
\begin{center}
$S_{1_{k}}(F) \subseteq S_{1_{k-1}}(F) \subseteq\ldots\subseteq S_{1_{2}}(F) \subseteq S_{1_{1}}(F)=S_{1}(F)\subseteq X$.
\end{center}
We recall that $S_{1_{h-1}}(F)$ is, near $u_{o}$, a $h$-codimensional submanifold of $X$ for $h=1,2,\ldots,k$ (cf. Proposition \ref{Pro242}). We also showed that, for $u$ near $u_{o}$ and $h=1,2,\ldots,k$,
\begin{center}
$T_{u}S_{1_{h}}(F) =\cap_{j=1}^{h} N(I_{j}(\varphi,\psi)(u)) ,\;\;\forall (\varphi,\psi)\in \mathcal{P}(F,u_{o})$.
\end{center}
This equality is obtained by combining the description of the tangent spaces of $S_{1_{h}}(F)$, given in Proposition \ref{Pro242}, with part (jj)$_{k}$  of  Proposition \ref{Pro245}.\\
\indent  e)  The strata $S_{1_{2}}(F),\ldots,S_{1_{k+1}}(F)$ near a singularity $u_{o}$ for which (b$_{1}$) and (b$_{2}$) are satisfied have a simple geometric meaning. Indeed, it is possible to show that, for  $h=1,2,\ldots,k$,
\begin{center}
$S_{1_{h+1}}(F)=\{u\in S_{1_{h}}(F):$ the straight line $N(F'(u))$ of $X$ is contained in the\\
$h$-codimensional subspace  $T_{u}S_{1_{h}}(F)$ of $X\}$.
\end{center}
In fact, $\forall \;(\varphi,\psi)\in \mathcal{P}(F,u_{o})$,
\begin{center}
$S_{1_{h}}(F)=\{u \text{ near }u_{o}:J_{0}(\varphi,\psi)(u)=\ldots=J_{h-1}(\varphi,\psi)(u)=0\}$.
\end{center}
Hence for all $u\in S_{1_{h}}(F)$ we have that $I_{\eta}(\varphi,\psi)(u)\varphi(u)=J_{\eta}(\varphi,\psi)(u)=0,\eta=1,2,\ldots,h-1$. \\
Thus $\varphi(u)\in T_{u}S_{1_{h}}(F)=\cap_{j=1}^{h}N(I_{j}(\varphi,\psi)(u))$ \textit{iff} $I_{h}(\varphi,\psi)(u)\varphi(u)=J_{h}(\varphi,\psi)(u)=0$ \textit{iff} $u\in S_{1_{h+1}}(F)$. Moreover, $S_{1_{h}}(F)\subseteq S_{1}(F)$ and so, for $u\in S_{1_{h}}(F),\varphi(u)$ spans $N(F'(u))$ and this allows us to conclude.
\end{remark}

\begin{remark} In the above remark we saw that, under conditions of $k$-transversality (b$_{1}$) and (b$_{2}$),
\begin{equation*}
\begin{split}
S_{1_{h+1}}(F)=\{&u\in S_{1_{h}}(F):\text{the straight line }N(F'(u))\text{ of }X\text{ is contained in the }\\
&h\text{-codimensional subspace }T_{u}S_{1_{h}}(F)\text{ of }X \}.
\end{split}
\end{equation*}
Hence we can say that $S_{1_{h+1}}(F), h=0,1,\ldots,k$, coincides with the subset of all points $u$ of $S_{1_{h}}(F)$ (setting $S_{1_{0}}(F):=X$) where the restriction of the map $F$ to the submanifold $S_{1_{h}}(F), F:S_{1_{h}}(F)\subseteq X \rightarrow V\subseteq Y$, has derivative  $F'(u): T_{u}S_{1_{h}}(F) \subseteq X \rightarrow Y$ with one-dimensional kernel $N(F'(u))$.\\
\indent We shall prove below the following two facts, for any $u\in S_{1_{h}}(F)$ near $u_{o}$: 
\begin{enumerate}
\item[1)] if $F'(u):T_{u}S_{1_{h}}(F)\subseteq X\rightarrow Y$ is injective, then\\
codim$_{Y} F'(u)(T_{u}S_{1_{h}}(F)) = \text{codim}_{X}T_{u}S_{1_{h}}(F)=h$;
\item[2)] if $F'(u):T_{u}S_{1_{h}}(F)\subseteq X\rightarrow Y$ has one-dimensional kernel $N(F'(u))$, then  codim$_{Y} F'(u)(T_{u}S_{1_{h}}(F)) = \text{codim}_{X}T_{u}S_{1_{h}}(F)+1=h+1$.
\end{enumerate}
In the finite-dimensional case, i.e. dim$\,X\,$=\,dim$\,Y<+\infty$, 1) and 2) imply that $S_{1_{h+1}}(F)$ is, near $u_{o}$, the subset of all points $u$ of  $S_{1_{h}}(F)$ where  $F'(u):T_{u}S_{1_{h}}(F)\rightarrow Y$ drops rank by 1.\\
\indent
This agrees with the usual \textit{Thom-Boardman stratification} of singularities, cf.\cite{Bo} and \cite{G-G}, i.e. it is immediate to check that
\begin{equation*}
S_{1_{h+1}}(F)\equiv S_{\underbrace{\text{\scriptsize 1,...,1,1}}_{h+1}}(F) \;\; (\text{here denoted by }  S_{1,\ldots,1,1}(F)).
\end{equation*}
We recall that the \textit{Thom-Boardman stratum }$S_{1,\ldots,1,1}(F)$ is inductively defined in the following way: as usual, the stratum $S_{1}(F)$ is the subset of all simple singularities for $F$ (cf. Definition \ref{d112}) and, provided that 
\begin{equation*}
S_{\underbrace{\text{\scriptsize 1,...,1}}_{h}}(F) \;\; (\text{here denoted by }  S_{1,\ldots,1}(F))
\end{equation*}
is a submanifold, then the next stratum is given by\\
$S_{1,\ldots,1,1}(F):=\{u\in S_{1,\ldots,1}(F):$ the map $F:S_{1,\ldots,1}(F)\subseteq X\rightarrow V\subseteq Y$  drops rank by 1\}.\\
\indent Let us prove claims 1) and 2).\\
We recall that  dim\,$N(F'(u))=\text{codim}\,R(F'(u))=1$, with $u\in S_{1}(F)$ and $F'(u)$ a 0-Fredholm operator. Moreover, codim$_{X}T_{u}S_{1_{h}}(F)=h$.\\
If $F'(u):T_{u}S_{1_{h}}(F)\subseteq X\rightarrow Y$ is injective, i.e. $N(F'(u))\nsubseteq T_{u}S_{1_{h}}(F)$, then there exists a $(h-1)$-dimensional subspace $X_{h-1}$ of $X$ such that  
\begin{equation*}
X=T_{u}S_{1_{h}}(F) \oplus N(F'(u))\oplus X_{h-1} .
\end{equation*}
Hence $R(F'(u))=F'(u)(X)=F'(u)(T_{u}S_{1_{h}}(F))\oplus F'(u)(X_{h-1})$. By construction $X_{h-1}\cong F'(u)(X_{h-1})$, thus $F'(u)(T_{u}S_{1_{h}}(F))$ has codimension $h-1$ in $R(F'(u))$ which has codimension 1 in $Y$. Therefore $F'(u)(T_{u}S_{1_{h}}(F))$ has codimension $h$ in $Y$.\\
\indent
When $F'(u):T_{u}S_{1_{h}}(F)\subseteq X\rightarrow Y$ has one-dimensional kernel $N(F'(u))$, i.e. $N(F'(u))\subseteq T_{u}S_{1_{h}}(F)$, we consider an h-dimensional subspace $X_{h}$ of $X$ such that $X=T_{u}S_{1_{h}}(F)\oplus X_{h}$. Since $N(F'(u))\nsubseteq X_{h}$, arguing as before one proves that $F'(u)(T_{u}S_{1_{h}}(F))$  has codimension $h+1$ in  $Y$.
\end{remark}
\vspace{2pt}
\subsection{Results on the Classification}\label{ss25}
\quad Here we give a positive answer to all the issues stated after Definition \ref{D211}. First, we discuss the independence of the definition from the chosen f-pair and the completeness of the classification. Then we provide some information about the singular points near a 1-transverse singularity. Additionally, in Theorem \ref{Teo254} we present a useful characterization of the singular points near a $k$-transverse singularity. Finally, we prove the invariance of the definition with respect to changes of coordinates. \vspace{8pt}
\begin{remark} \label{Rem251} (Well-posedness of Definition \ref{D211}). As already noted in Remark \ref{Rem246}, statement i) of Theorem \ref{Teo243} allows us to conclude that the definition of $k$-transverse singularity, given in (T$_{k}$) of Definition \ref{D211}, is well-posed. The same result is true for the definitions of $k$-ordinary and maximal $k$-transverse singularity, given in (S$_{k}$) and (M$_{k}$), and this derives from statements ii) and iii) of Theorem \ref{Teo243}, respectively. Also when the map $F$ is of class $C^{\infty}$ we can consider condition (T$_{\infty}$) in Definition \ref{D211}. Then, by using statement iv) of Theorem \ref{Teo243} we can conclude that the definition of $\infty$-transverse singularity is well-posed too.
\end{remark}
\vspace{4pt}
\par We are now able to prove that a 1-transverse singularity for a $C^{\infty}$ 0-Fredholm map $F$ verifies one and only one of conditions (S$_{k}$), (M$_{k}$), (T$_{\infty}$). In fact, if $u_{o}\in S_{1}(F)$ satisfies one of conditions (S$_{k}$), (M$_{k}$), (T$_{\infty}$) for some $(\varphi,\psi)\in \mathcal{P}(F,u_{o})$ then, from Lemma \ref{Lem244}, we get that $u_{o}$ is a $k$-transverse singularity for $F$. In particular, $I_{1}(\varphi,\psi)(u_{o})\neq 0$ and thus $u_{o}$ is 1-transverse for $F$. Hence a $k$-ordinary [maximal $k$-transverse, $\infty$-transverse] singularity is a 1-transverse singularity. Vice versa, each 1-transverse simple singularity for a $C^{\infty}$ 0-Fredholm map $F$ has to satisfy one and only one condition among (S$_{k}$), (M$_{k}$), (T$_{\infty}$) for a suitable $k$: this is better stated and proved in the following result. For the sake of simplicity, we only consider the most interesting situation where all kind of singularities can possibly occur, i.e. smooth maps and infinite dimensions, even though this result can be easily generalized to all cases. \\ 
\begin{proposition} \label{Pro252}
(Classification of Singularities). \textit{ Let} $U, V$ \textit{be open subsets in the }\textit{B-spaces} $X,Y$ \textit{and} $F:U \subseteq X\rightarrow V\subseteq Y$\textit{ a } $C^{d}$ 0\textit{-Fredholm map } $(d=\infty\textit{ or } d=\omega)$.\textit{ Let us choose }$u_{o}\in S_{1}(F)$\textit{; then }$u_{o}$\textit{ is a }1-transverse singularity \textit{if and only if one and only one of the following conditions is satisfied}:
\begin{enumerate}
\item[i)]\textit{there exists a unique integer } $k\geq 1$ \textit{such that }$u_{o}$\textit{ is a} $ k $-singularity;
\item[ii)]\textit{there exists a unique integer } $k\geq 1$\textit{ such that }$u_{o}$\textit{ is a }maximal k-transverse
singularity;
\item[iii)]$u_{o}$\textit{ is an }$\infty$-transverse singularity.
\end{enumerate}
\end{proposition}
\par
\textbf{Proof.} From the above discussion, we just need to prove the \textit{only if} part of the statement. To this end we define 
\begin{center}
$H:=\{h\in \mathbb{N},h\geq 1:u_{o}$ is an $h$-transverse singularity for $F\}$.
\end{center}
Since such a subset is non-empty (because $1\!\in\!H$) we can consider $\text{sup}\,H$. If $H$ is bounded then $k$ := sup$\,H$ coincides with max$\,H$. Thus $u_{o}$ is $k$-transverse and, for a given $(\varphi,\psi)\in \mathcal{P}(F,u_{o})$, one has either $J_{k}(\varphi,\psi)(u_{o}) \neq 0$ or $J_{k}(\varphi,\psi)(u_{o})=0$. From Definition 2.1.1 it is not difficult to see that in the first case $u_{o}$ is a  $k$-singularity, while in the second one $u_{o}$ is a maximal $k$-transverse singularity. Finally, when $H$ is not bounded it is easily proved that $u_{o}$ is an $\infty$-transverse singularity for $F$. \\
\indent Now let us suppose that $u_{o}$ is a $k$-singularity for some $k\geq 1$ and let us show that $u_{o}$  is not an $m$-singularity for any integer $m\neq k$ and that neither ii) nor iii) can occur.\\
If  $u_{o}$ were also an $m$-singularity, $m>k$, then $u_{o}$ would be $m$-transverse and this would contradict $k=\text{max }H$; if $u_{o}$ were an $m$-singularity, $m\leq k-1$, then $J_{m}(\varphi,\psi)(u_{o})\neq 0$ but this is absurd because $u_{o}$ is $k$-transverse and this implies that $J_{0}(\varphi,\psi)(u_{o})=\ldots=J_{k-1}(\varphi,\psi)(u_{o})=0$.\\
Now let us assume that ii) holds, i.e. $u_{o}$ is also a maximal $m$-transverse singularity for some integer $m$: if $m\geq k$, then by definition $J_{0}(\varphi,\psi)(u_{o})=\ldots=J_{m}(\varphi,\psi)(u_{o})=0$ but, on the other hand, $J_{k}(\varphi,\psi)(u_{o}) \neq 0$ because $u_{o}$ is a $k$-singularity. If $\;m+1\leq k\;$ then, again by definition, $I_{1}(\varphi,\psi)(u_{o}),\ldots,I_{m+1}(\varphi,\psi)(u_{o})$ are not l.i. and this cannot be because $k=\text{max} \,H$, hence $I_{1}(\varphi,\psi)(u_{o}),\ldots,I_{k}(\varphi,\psi)(u_{o})$ are l.i.. Finally, $u_{o}$ is not an $\infty$-transverse singularity because for any integer $m$ one should have $J_{m}(\varphi,\psi)(u_{o})=0$ while $J_{k}(\varphi,\psi)(u_{o})\neq 0$.\\
\indent When $u_{o}$ is a maximal $k$-transverse singularity or an $\infty$-transverse singularity one can argue in a similar way and thus conclude the proof of the proposition.$\bracevert$ \\
\par The next result provides some information about the singular points near a 1-transverse singularity that is of some interest in itself and will be also used in the following.\\
\begin{proposition} \label{Pro253}
\textit{ Let} $U, V$ \textit{be open subsets in the }\textit{B-spaces} $X,Y$ and $F:U \subseteq X\rightarrow V\subseteq Y$\textit{ a } $C^{d}$ 0\textit{-Fredholm map }$(d\geq 2)$.\textit{ Let }$k$ \textit{be an integer such that }$1\leq k\leq d-1 $\textit{; if }$ u_{o}\in S_{1}(F)$\textit{ is a  }1-\textit{transverse singularity for }$F$\textit{ then, for a given }$(\varphi,\psi)\in \mathcal{P}(F,u_{o})$,\textit{ there exists a neighbourhood }$W$\textit{ of }$u_{o}$\textit{ in }$X$\textit{ such that for any }$u\in W$\textit{ the following equivalence holds:} \\
\indent $u\in S_{1}(F)$\textit{ and }$u$\textit{ is a }$k$-transverse [$k$-ordinary, maximal $k$-transverse, $\infty$-transverse] singularity\textit{ for }$F \Leftrightarrow$\textit{ the conditions }(T$_{k}$) [(S$_{k}$), (M$_{k}$), (T$_{\infty}$)]\textit{ in }Definition \ref{D211}\textit{ are satisfied by the functionals }$J_{h}(\varphi,\psi),I_{h}(\varphi,\psi)$\textit{ evaluated at }$u$.
\end{proposition}
\par\textbf{ Proof.} To prove the equivalence in the case of $k$-transverse singularities we note that, given that $u_{o}$ a 1-transverse singularity, one has that $I_{1}(\varphi,\psi)(u_{o})\neq 0$ and, from Theorem 2.3.3, $J_{0}(\varphi,\psi)^{-1}(0)=S_{1}(F)$ near $u_{o}$. Thus it suffices to take a neighbourhood W of $u_{o}$ such that\\
i)\; the fibering pair $(\varphi,\psi)$  is defined on $W$;\\
ii) for $u \in W$ it is true that $J_{0}(\varphi,\psi)(u)=0 \textit{ iff }u\in S_{1}(F)$.\\
With this choice the above equivalence is easily proved. As a preliminary remark, note that $W$ is also a neighbourhood of any $u\in W$ and then, for a given $u\in S_{1}(F)\cap W$, one has that $(\varphi,\psi)\in \mathcal{P}(F,u)$: this is true because of i) and the very definition of fibering pair (see Definition \ref{d113}).\\
\textit{Proof of }$\Rightarrow$: if $u\in S_{1}(F)\cap W$ is a $k$-transverse singularity for $F$ then, since Definition \ref{D211} is independent from the given $(\varphi,\psi)\in \mathcal{P}(F,u)$ (as recalled in Remark \ref{Rem251}), we have that $J_{0}(\varphi,\psi)(u)=\ldots=J_{k-1}(\varphi,\psi)(u)=0$ and $I_{1}(\varphi,\psi)(u),\ldots, I_{k}(\varphi,\psi)(u)$ are l.i.. \\
\textit{Proof of }$\Leftarrow$: if $u\in W$, $J_{0}(\varphi,\psi)(u)=\ldots=J_{k-1}(\varphi,\psi)(u)=0$ and $I_{1}(\varphi,\psi)(u),\ldots$,\\
$I_{k}(\varphi,\psi)(u)$ are l.i. then, in particular, $u\in S_{1}(F)$ because of ii). Thus $(\varphi,\psi)\in  \mathcal{P}(F,u)$ and, from Definition \ref{D211} applied to the pair $(\varphi,\psi)$, we get that $u$ is a $k$-transverse singularity for $F$.\\
\indent The other equivalences are proved in a similar way. $\bracevert$\\
\par
The following characterizations provide a deep insight into the close relationship between the classification of singularities and the singular strata $S_{1_{h}}(F)$. Thanks to these characterizations, it is possible to justify the slight abuse of notation when writing $S_{1_{h}}(F)$ instead of $S_{1_{h}}(F)(\varphi,\psi)$, as already announced in c) of Remark \ref{Rem246}. For instance, given $(\varphi,\psi), (\widetilde{\varphi},\widetilde{\psi})\in \mathcal{P}(F,u_{o})$ where $u_{o}$ is a $k$-transverse singularity for $F$, we obtain from equivalence a) below that $S_{1_{k}}(F)(\varphi,\psi)=S_{1_{k}}(F)(\widetilde{\varphi},\widetilde{\psi})= \{u\in S_{1}(F):u$ is $k$-transverse for F\} $\cap\,W$ for a suitable neighbourhood $W$ of $u_{o}$. Moreover, these characterizations are used in  \cite{B-D 2} to describe the local behaviour of the map $F$ near a $k$-transverse singularity $u_{o}$ and in \cite{B-D 3} to find an operative way to determine the kind of a given singularity.\\
\begin{theorem} \label{Teo254}
\textit{ Let }$U,V$ \textit{ be open subsets in the $B$-spaces }$X,Y$ \textit{and }$F:U\subseteq X\rightarrow V\subseteq Y$\textit{ a }$C^{d}$ 0-\textit{Fredholm map }$(d\geq 2)$. \textit{Let us assume that }$u_{o}\in S_{1}(F)$\textit{ is a }$k$-transverse singularity\textit{ for }$F$\textit{ where }$1\leq k \leq d-1$.\textit{ Then, for a given }$(\varphi,\psi)\in \mathcal{P}(F,u_{o})$,\textit{ there exists a neighbourhood }$W$\textit{ of }$u_{o}$\textit{ in }$X$\textit{ such that for }$u\in W$\textit{ the following equivalences hold}:
\begin{align*}
\text{a)}\; u\,\in \, S_{1_{k}}(F)\,\Leftrightarrow \, u\in S_{1}&(F)\textit{ and it is }k\text{-transverse};\\
\text{b)}\;u\,\in \, S_{1_{k}}(F)\setminus S_{1_{k+1}}(F)&\Leftrightarrow u\in S_{1}(F)\textit{ and it is a } k\text{-singularity} \\
&\Leftrightarrow u\in S_{1_{k}}(F)\textit{ and }N(F'(u))\nsubseteq T_{u}S_{1_{k}}(F) \\
&\Leftrightarrow u\in S_{1_{k}}(F)\textit{ and }\varphi(u)\notin T_{u}S_{1_{k}}(F);\\
\text{c)}\; u \in S_{1_{k+1}}(F)\Leftrightarrow u\in S_{1}&(F)\textit{ and it is }k\text{-transverse},\textit{ not a } k\text{-singularity}\\
\Leftrightarrow u\in S_{1_{k}}&(F)\textit{ and }N(F'(u))\subseteq T_{u}S_{1_{k}}(F)\\
\Leftrightarrow u\in S_{1_{k}}&(F)\textit{ and }\varphi(u)\in T_{u}S_{1_{k}}(F);\\
\text{d)}\; \textit{for }d\geq 3\textit{ and }1\leq k\leq &\,d-2,\\
u\in S_{1_{k+1}}(F) \setminus S_{1_{k+2}}&(F)\Leftrightarrow u\in S_{1}(F)\textit{ and it is a } (k+1)\text{-singularity};\\
\text{e)}\;\textit{for }d\geq 3\textit{ and }1\leq k\leq &\,d-2,\\
u\in S_{1_{k+2}}(F) \Leftrightarrow u\in &S_{1}(F)\textit{ and it is }k\text{-transverse},\textit{ not a }k\text{-singularity},\\
\textit{not a }&(k+1)\text{-singularity}.
\end{align*}
\end{theorem}
\par
\textbf{Proof.} Since $u_{o}$ is $k$-transverse it suffices to take $W$ such that \\
j)    the equivalences of  Proposition \ref{Pro253}  are valid on $W$;\\
jj)   for $u\in W$ the functionals $I_{1}(\varphi,\psi)(u),\ldots,I_{k}(\varphi,\psi)(u)$  are  l.i..\\
\indent In fact the equivalences mentioned in j) can be used because a $k$-transverse singularity is 1-transverse.
Moreover, near $u_{o}$, the linear independence of $I_{1}(\varphi,\psi)(u),\ldots$,\\
$I_{k}(\varphi,\psi)(u)$ is assured by the continuity of the functionals $I_{h}(\varphi,\psi)$ as proved in Corollary \ref{Co222}, hence jj) is true on a neighbourhood of $u_{o}$.\\
\indent From conditions j) and jj) we can prove the equivalences a), ..., e),  for a given $u\in W$.\\
\indent \textit{Proof of} a). From Proposition \ref{Pro253} we have that $u\in S_{1}(F)$ and it is $k$-transverse \textit{iff} $J_{0}(\varphi,\psi)(u)=\ldots=J_{k-1}(\varphi,\psi)(u)=0$ \;\;and\;\; $I_{1}(\varphi,\psi)(u),\ldots,I_{k}(\varphi,\psi)(u)$  are l.i.. Because of jj) this is true \textit{iff} $J_{0}(\varphi,\psi)(u)=\ldots=J_{k-1}(\varphi,\psi)(u)=0$. Then one obtains that $u\in S_{1}(F)$ and it is  $k$-transverse \textit{iff} $u\in S_{1_{k}}(F)$.\\
\indent \textit{Proof of} b). From Proposition \ref{Pro253} we know that $u\in S_{1}(F)$ and it is a $k$-singularity \textit{iff} $J_{0}(\varphi,\psi)(u)=\ldots= J_{k-1}(\varphi,\psi)(u)=0, \; J_{k}(\varphi,\psi)(u)\neq 0$ and $I_{1}(\varphi,\psi)(u)$,\\
$\ldots,I_{k-1}(\varphi,\psi)(u)$ are l.i.. From jj) this is true \textit{iff} $J_{0}(\varphi,\psi)(u)=\ldots=J_{k-1}(\varphi,\psi)(u)=0, \; J_{k}(\varphi,\psi)(u)\neq 0$. Then one has that $u \in S_{1}(F)$ and it is a $k$-singularity \textit{iff} $u\in S_{1_{k}}(F)$ but $u\notin S_{1_{k+1}}(F)$ and this proves the first equivalence of b).\\
\indent  In order to prove the other equivalences we recall that we have the equality $S_{1_{k+1}}(F)=\{u\in S_{1_{k}}(F):N(F'(u))\subseteq T_{u}S_{1_{k}}(F)\}$ (see Remark \ref{Rem246}, e)). Therefore $u\in S_{1_{k}}\setminus (F) \ S_{1_{k+1}}(F)$ \textit{iff} $u\in S_{1_{k}}(F)$ and $N(F'(u))\nsubseteq T_{u}S_{1_{k}}(F)$. The last statement is equivalent to saying that $\varphi(u)\notin T_{u}S_{1_{k}}(F)$. In fact, $u\in S_{1_{k}}(F)$ implies that $u\in S_{1}(F)$ (thanks to a), for instance) and thus $\varphi(u)$ spans $N(F'(u))$ (cf. Definition \ref{d113}). Hence $N(F'(u))\nsubseteq T_{u}S_{1_{k}}(F)$ \textit{iff} $\varphi(u)\notin T_{u}S_{1_{k}}(F)$.\\
\indent \textit{Proof of} c). We recall that $S_{1_{k+1}}(F)\subseteq S_{1_{k}}(F)\subseteq S_{1}(F)$ (see Remark \ref{Rem246}, c)); then $u\in S_{1_{k+1}}(F)$  \textit{iff} $u\in S_{1_{k}}(F)$ but $u\notin S_{1_{k}}(F)\setminus S_{1_{k+1}}(F)$. By combining a) and b) this amounts to saying that $u\in S_{1}(F)$, $u$ is $k$-transverse but it is not a $k$-singularity. The other equivalences are also an immediate consequence of the analogues in b).\\
\indent  \textit{Proof of} d). From Proposition \ref{Pro253} one has that $u\in S_{1}(F)$ and it is a $(k+1)$-singularity \textit{iff}  $J_{0}(\varphi,\psi)(u)=\ldots=J_{k}(\varphi,\psi)(u)=0,\; J_{k+1}(\varphi,\psi)(u)\neq 0$ and $I_{1}(\varphi,\psi)(u), \ldots,I_{k}(\varphi,\psi)(u)$ are l.i., that is \textit{iff} $J_{0}(\varphi,\psi)(u)=\ldots=J_{k}(\varphi,\psi)(u)=0, J_{k+1}(\varphi,\psi)(u)\neq 0$, from jj). Hence we get that $u\in S_{1}(F)$ and it is a $(k+1)$-singularity \textit{iff} $u\in S_{1_{k+1}}(F)$ and $u\notin S_{1_{k+2}}(F)$.\\
\indent  \textit{Proof of} e). Since  $S_{1_{k+2}}(F)\subseteq S_{1_{k+1}}(F)\subseteq S_{1}(F)$, $u\in S_{1_{k+2}}(F)$  \textit{iff} $u\in S_{1_{k+1}}(F)$  but $u\notin S_{1_{k+1}}(F) \setminus S_{1_{k+2}}(F)$. By combining c) and d), this is equivalent to saying that $u\in S_{1}(F)$, $u$ is $k$-transverse, $u$ is not a $k$-singularity and it is not a $(k+1)$-singularity. $\bracevert$\\
\par With the following theorem we prove that Definition \ref{D211} is\textit{ invariant }with respect to changes of coordinates.\vspace{6pt}
\begin{theorem} \label{Teo255}
(Invariance Theorem). \textit{Let us suppose that there exists a} l.c.d. \textit{of class}  $C^{d}\,(d \geq 2)$\\
\begin{center}
\begin{tabular}{ccc}
&$F \quad $&  \medskip \\
$u_o \in U \subseteq X \quad $&$ \rightarrow \quad $&$ V \subseteq Y $ \medskip  \\
$\gamma \downarrow $&&$ \downarrow \delta$ \medskip \\ 
$\tilde{u}_{o} \in \widetilde{U} \subseteq \tilde{X} \quad $&$ \rightarrow  \quad $&$ \widetilde{V} \subseteq \widetilde{Y}$ , \medskip \\
&$\widetilde{F} \quad $&\\
\end{tabular}
\end{center}
\textit{where }$U,\,V,\, \widetilde{U}$\textit{ and }$\widetilde{V}$\textit{ are open subsets in the $B$-spaces }$X,\,Y,\,\widetilde{X}$\textit{ and }$\widetilde{Y}$\textit{ respectively. Moreover, let }$F$\textit{ be a }0\textit{-Fredholm map and }$u_{o}\in S_{1}(F)$. \textit{Then, for a fixed integer }$k, 1\leq k\leq d-1$,\textit{ we have that:}\\
\indent $u_{o}$\textit{ is a }$k$-transverse [$k$-ordinary, maximal $k$-transverse, $\infty$-transverse] singularity\textit{ for }$F$\textit{ if and only if }$\tilde{u}_{o}$\textit{ is a }$k$-transverse [$k$-ordinary, maximal $k$-transverse, $\infty$-transverse] singularity \textit{for} $\widetilde{F}$.
\end{theorem}
\par
\textbf{Proof.} Let us note that, up to shrinking the neighbourhoods $U,\,V,\,\widetilde{U},\,\widetilde{V}$, we already considered a similar commutative diagram in Subsection 1.4.2, where we studied how f-pairs and fibering functionals are modified by changes of coordinates. There we proved that $\widetilde{F}$ is a 0-Fredholm map and $\tilde{u}_{o}\in S_{1}(\widetilde{F})$. Indeed, by definition of l.c.d. $\tilde{u}_{o}=\gamma(u_{o})$ and, from (\ref{144}), $S_{1}(\widetilde{F})=\gamma(S_{1}(F))$. Therefore in order to prove the statement it suffices to show that if $u_{o}$ is a $k$-transverse [$k$-ordinary, maximal $k$-transverse, $\infty$-transverse] singularity for $F$ then $\tilde{u}_{o}$ is a $k$-transverse [$k$-ordinary, maximal $k$-transverse, $\infty$-transverse] singularity for $\widetilde{F}$.\\
\indent Let $(\varphi,\psi)$ be a fibering pair for $F$ near $u_{o}$. Then $(\widetilde{\varphi},\widetilde{\psi}):= \mathcal{T}[\gamma,\delta](\varphi,\psi)$, the pair-transform of $(\varphi,\psi)$, is a fibering pair for $\widetilde{F}$ near  $\tilde{u}_{o}$ thanks to Proposition \ref{Pro144}. Moreover, formulas (\ref{146}) give the fibering functionals associated with $(\widetilde{\varphi},\widetilde{\psi})$. Then we have, for all $\tilde{u} \in \widetilde{U}$,
\begin{equation}\label{251}
\begin{split}
J_{h}(\widetilde{\varphi},\widetilde{\psi})(\tilde{u})&=J_{h}(\varphi,\psi)(\gamma^{-1}(\tilde{u}))\, , \, 0\leq h\leq d-1 ;\\
I_{h}(\widetilde{\varphi},\widetilde{\psi})(\tilde{u})&=I_{h}(\varphi,\psi)(\gamma^{-1}(\tilde{u})) (\gamma^{-1})'(\tilde{u})\, , \, 1\leq h \leq d-1.
\end{split}
\end{equation}
From the first equalities we easily get that
\begin{equation}\label{252}
J_{h}(\widetilde{\varphi},\widetilde{\psi})(\tilde{u}_{o})= 0 \, \Leftrightarrow \,J_{h}(\varphi,\psi)(u_{o})=0.
\end{equation}
The second ones yield
\begin{equation}\label{253}
I_{1}(\widetilde{\varphi},\widetilde{\psi})(\tilde{u}_{o}),\ldots,I_{h}(\widetilde{\varphi},\widetilde{\psi})(\tilde{u}_{o})\text{ are l.i.}\Leftrightarrow I_{1}(\varphi,\psi)(u_{o}),\ldots,I_{h}(\varphi,\psi)(u_{o})\text{ are l.i..}\vspace{8 pt}
\end{equation}
\indent Let us prove the equivalence (\ref{253}). Let $I_{1}(\varphi,\psi)(u_{o}),\ldots,I_{h}(\varphi,\psi)(u_{o})$ be l.i.. From the Algebraic Lemma \ref{Lem221}, a) $\Rightarrow$ d), there exist $v_{s}\in X, s=1,\ldots,h$, such that $I_{t}(\varphi,\psi)(u_{o})v_{s}=\delta_{ts},t,s=1,\ldots,h$. It is now convenient to define vectors $w_{s}\in \widetilde{X}$ as $w_{s}:=\gamma\hspace{1pt}'(u_{o})v_{s},s=1,\ldots,h$. Since  $((\gamma^{-1})'(\tilde{u}_{o}))^{-1}=\gamma\hspace{1pt}'((\gamma^{-1}(\tilde{u}_{o}))=\gamma\hspace{1pt}'(u_{o})$, from (\ref{251}) we obtain that $I_{t}(\widetilde{\varphi},\widetilde{\psi})(\tilde{u}_{o})w_{s}=I_{t}(\varphi,\psi)(\gamma^{-1}(\tilde{u}_{o}))(\gamma^{-1})'(\tilde{u}_{o})w_{s}=I_{t}(\varphi,\psi)(u_{o})(\gamma\hspace{1pt}'(u_{o}))^{-1}w_{s}=I_{t}(\varphi,\psi)(u_{o})v_{s}=\delta_{ts}$, for  $t,s=1,\ldots,h $. Then the Algebraic lemma, d) $\Rightarrow$ a), implies that $I_{1}(\widetilde{\varphi},\widetilde{\psi})(\tilde{u}_{o}),\ldots,I_{h}(\widetilde{\varphi},\widetilde{\psi})(\tilde{u}_{o})$ are l.i.. \\
We have shown the sufficient condition of (\ref{253}). The converse follows in the same manner by writing $I_{t}(\varphi,\psi)(u_{o})=I_{t}(\widetilde{\varphi},\widetilde{\psi})(\tilde{u}_{o})\gamma\hspace{1pt}'(u_{o}) $.\\
\indent Since Definition \ref{D211} does not depend on the f-pairs $(\varphi,\psi)$ and $(\widetilde{\varphi},\widetilde{\psi})$, defined near $u_{o}$ and $\tilde{u}_{o}$ respectively, and the equivalences (\ref{252}) and (\ref{253}) hold, then the proof is complete.$\bracevert$\\
\par
We could use the previous result to give an analogous definition of $k$-transverse [$k$-ordinary, maximal $k$-transverse, $\infty$-transverse] singularities for 0-Fredholm maps $F$ between \textit{Banach manifolds}. One would only have to work by means of local charts and apply Definition 2.1.1 to the related local representative map of $F$, which is of course between Banach spaces. Then the Invariance Theorem would imply that this definition is independent from the chosen local charts.\vspace{6pt}
\begin{remark} From formulas (\ref{251}) we easily get
\begin{equation*}
S_{1_{h}}(\widetilde{F})(\widetilde{\varphi},\widetilde{\psi})=\gamma(S_{1_{h}}(F)(\varphi,\psi)), \;1\leq h\leq d ,
\end{equation*} 
because $S_{1_{h}}(F)(\varphi,\psi)=\{u\text{ near }u_{o}:J_{\eta}(\varphi,\psi)(u)=0, \eta=0,1,\ldots,h-1\}$ and  $S_{1_{h}}(\widetilde{F})(\widetilde{\varphi},\widetilde{\psi})=\{\tilde{u}\text{ near }\tilde{u}_{o}:J_{\eta}(\widetilde{\varphi},\widetilde{\psi})(\tilde{u})=0, \eta=0,1,\ldots,h-1\}$.\\
\indent Moreover, when $u_{o}$ satisfies one of the conditions in Definition \ref{D211} $u_{o}$ is a $k$-transverse singularity for $F$. Thanks to the Invariance Theorem, we can say that $\tilde{u}_{o}$ is a $k$-transverse singularity for $\widetilde{F}$. Therefore, also as a consequence of the characterizations proved in Theorem \ref{Teo254}, we obtain that if $u_{o}$ is a $k$-transverse [$k$-ordinary, maximal $k$-transverse, $\infty$-transverse] singularity  then:\\
\begin{equation*}
S_{1_{h}}(\widetilde{F})=\gamma(S_{1_{h}}(F)),\text{ for }1\leq h\leq k+2\text{  where }1\leq k\leq d-2
\end{equation*}
or
\begin{equation*}
S_{1_{h}}(\widetilde{F})=\gamma(S_{1_{h}}(F)),\text{ for }1\leq h\leq k+1\text{  where }k=d-1.
\end{equation*}
\end{remark}
\subsection{Singularities for LS-maps and ``Polynomial'' Examples}\label{ss26}
\quad In order to give examples of all the singularities introduced in Definition \ref{D211} it is useful to see how conditions (T$_{k}$), (S$_{k}$), (M$_{k}$), (T$_{\infty}$) can be restated for the LS-maps. For this class of maps we find the equivalent conditions $(\text{T}_{k}^{\,\prime})$, $(\text{S}_{k}^{\,\prime})$, $(\text{M}_{k}^{\,\prime})$, $(\text{T}_{\infty}^{\,\prime})$  which are also a basic tool in the proof of the dual characterization of singularities we give in Chapter 1 of \cite{B-D 3}.
\subsubsection{} Let $\varXi$ be a $B$-space, $U,V$ open subsets in $\mathbb{R}\times \varXi$ and $F:U\subseteq \mathbb{R} \times \varXi \rightarrow V \subseteq \mathbb{R} \times \varXi$ a map having the form $F(t,\xi)=(f(t,\xi),\xi), \, \forall \,(t,\xi)\in U$. Here $f:U\subseteq \mathbb{R} \times \varXi \rightarrow \mathbb{R}$  is a $C^{d}$ function, $d\geq 2$, and so $F$ is a $C^{d}$ LS-map. We recall that this kind of maps were studied in Subsection 1.3.1. There we saw that any LS-map is a 0-Fredholm map and $S_{1}(F)=\{(t,\xi)\in U:\dfrac{\partial f}{\partial t}(t,\xi)=0\}$. Moreover, the formulas \\
\begin{equation*}
\varphi_{C}(t,\xi)=(1,0)\in \mathbb{R}\times \varXi \; , \; \psi_{C}(t,\xi)=(1,-\frac{\partial f}{\partial t}(t,\xi))\in \mathbb{R}\times\varXi^{\ast}\cong(\mathbb{R}\times \varXi)^{\ast} ,
\end{equation*}
define the canonical f-pair $(\varphi_{C},\psi_{C})$ for $F$. We obtained in (\ref{135}) that\\
$J_{h}(\varphi_{C},\psi_{C})(t,\xi)=\dfrac{\partial^{h+1}f}{\partial t^{h+1}}(t,\xi)\in \mathbb{R}, \, h\geq 0\,,$\\
$I_{h}(\varphi_{C},\psi_{C})(t,\xi)=(\dfrac{\partial^{h+1}f}{\partial t^{h+1}}(t,\xi),\dfrac{\partial^{h+1}f}{\partial t^{h} \partial \xi}(t,\xi))\in \mathbb{R} \times \varXi^{\ast}\cong(\mathbb{R}\times \varXi)^{\ast} , (t,\xi)\in U, h\geq 1$.\\
The equivalent conditions $(\text{T}_{k}^{\,\prime})$, $(\text{S}_{k}^{\,\prime})$, $(\text{M}_{k}^{\,\prime})$, $(\text{T}_{\infty}^{\,\prime})$ for LS-maps we talked above are stated in the following four propositions.\\
\begin{proposition}\label{Pro262T}\textbf{(T).}\textit{ Let }$F:U\subseteq \mathbb{R} \times \varXi \rightarrow V\subseteq \mathbb{R} \times \varXi$\textit{ be a }$C^{d}$\textit{ LS-map }$(d\geq 2)$.\textit{ Then, for a given integer }$k$\textit{ such that }$1\leq k\leq d-1$,\textit{ a point }$(t,\xi)\in U$\textit{ is a }$k$-transverse singularity\textit{ for }$F$ \textit{if and only if the following condition is satisfied:}\\
\begin{align*}
(\text{T}_{k}^{\,\prime}) \qquad \qquad & \frac{\partial f}{\partial t}(t,\xi)=\ldots=\frac{\partial^{k}f}{\partial t^{k}}(t,\xi)=0 \,\, ;\\
& \textit{there exist } r_{\eta} \in \mathbb{R}, w_{\eta} \in \varXi, \eta=1,\ldots,k \textit{ such that}\\
& \frac{\partial ^{h+1}f}{\partial t^{h}\partial \xi }(t,\xi) w_{\eta}=\delta_{h\eta}, \, h=1,\ldots,k-1, \, \eta=1,\ldots,k \; ,\\
&r_{\eta}\frac{\partial^{k+1}f}{\partial t^{k+1}}(t,\xi)+ \frac{\partial^{k+1}f}{\partial t^{k}\partial \xi}(t,\xi) w_{\eta}=\delta_{k\eta}\; , \eta=1,\ldots,k. 
\end{align*}
\end{proposition}
\vspace{8pt}\par
\textbf{Proof.} From Definition \ref{D211}, $ (t,\xi)\in U$ is a $k$-transverse singularity for $F$ \textit{iff }$(t,\xi) \in S_{1}(F)$ and the conditions (T$_{k}$) hold for the canonical f-pair, i.e. \\
\begin{align*}
J_{h}(\varphi_{C},\psi_{C})(t,\xi)&=\frac{\partial^{h+1}f}{\partial t^{h+1}}(t,\xi)=0 \; ,\; h = 0,\ldots,k-1\,, \\ 
I_{h}(\varphi_{C},\psi_{C})(t,\xi)&=(\frac{\partial^{h+1}f}{\partial t^{h+1}}(t,\xi)+ \frac{\partial^{h+1}f}{\partial t^{h}\partial \xi}(t,\xi))\text{  are  l.i. , } h = 1, ..., k .
\end{align*}
Hence $(t,\xi)\in S_{1}(F)$  because  $\dfrac{\partial f}{\partial t}(t,\xi)=0$. \bigskip \\
By the Algebraic Lemma \ref{Lem221}, a) $\Leftrightarrow$ b), the above functionals  $I_{h}(\varphi_{C},\psi_{C})(t,\xi)$ are l.i.  \textit{iff} there exist $r_{\eta}\in \mathbb{R}, w_{\eta}\in \varXi, \eta = 1,\ldots,k$, such that  \\
\begin{center}
$(\dfrac{\partial^{h+1}f}{\partial t^{h+1}}(t,\xi)+ \dfrac{\partial^{h+1}f}{\partial t^{h}\partial \xi}(t,\xi))(r_{\eta}, w_{\eta})=r_{\eta}\dfrac{\partial^{h+1}f}{\partial t^{h+1}}(t,\xi) + \dfrac{\partial^{h+1}f}{\partial t^{h}\partial \xi}(t,\xi)w_{\eta}=\delta_{h\eta}$
\end{center}
for $h= 1,\ldots,k$.\bigskip \\
Since  $\dfrac{\partial^{h+1}f}{\partial t^{h+1}}(t,\xi)=0$ for $h=0,\ldots,k-1$ we easily conclude the proof.$\bracevert$\\
\begin{proposition}\label{Pro263S}
\textbf{(S).} \textit{Let} $F:U\subseteq \mathbb{R} \times \varXi \rightarrow V\subseteq \mathbf{R} \times \varXi$\textit{ be a }$C^{d}$\textit{ LS-map }$(d\geq 2)$.\textit{ Then, for a given integer }$k$\textit{ such that }$1\leq k\leq d-1$\textit{, a point }$(t,\xi)\in U$\textit{ is a }$k$-singularity\textit{ for }$F$ \textit{if and only if the following condition is satisfied:}
\begin{align*}
(\text{S}_{k}^{\,\prime}) \qquad \qquad & \frac{\partial f}{\partial t}(t,\xi)=\ldots=\frac{\partial^{k}f}{\partial t^{k}}(t,\xi)=0 \; , \; \frac{\partial^{k+1}f}{\partial t^{k+1}}(t,\xi) \neq 0 \,\, ;\\
&\textit{there exist } w_{\eta}\in \varXi , \eta=1,\ldots, k-1\textit{ such that }\\
&\frac{\partial ^{h+1}f}{\partial t^{h}\partial \xi}(t,\xi) w_{\eta}=\delta_{h\eta} \, , \, h,\eta=1,\ldots,k-1.
\end{align*}
\end{proposition}
\par
\textbf{Proof. }From Definition \ref{D211}, $(t,\xi)\in U$ is a $k$-singularity \textit{iff} $(t,\xi)\in S_{1}(F)$ and the conditions (S$_{k}$) hold for the canonical f-pair, that is 
\begin{align*} 
J_{h}(\varphi_{C},\psi_{C})(t,\xi) &=\frac{\partial^{h+1}f}{\partial t^{h+1}} (t,\xi)= 0, h=0\ldots,k-1 \, ,\\
J_{k}(\varphi_{C},\psi_{C})(t,\xi)&=\frac{\partial^{k+1}}{\partial t^{k+1}}(t,\xi) \neq 0 \,; \\
I_{h}(\varphi_{C},\psi_{C})(t,\xi)&=(\frac{\partial^{h+1}f}{\partial t^{h+1}} (t,\xi), \frac{\partial^{h+1}f}{\partial t^{h}\partial \xi}(t,\xi))  \text{ are l.i. }, \;\text{for } h=1,\ldots,k-1.
\end{align*}
In the previous proof we already showed that  $(t,\xi)\in S_{1}(F)$ and that the functionals $I_{h}(\varphi_{C},\psi_{C})(t,\xi),h=1,\ldots,k-1$, are  l.i. \textit{iff} there exist $r_{\eta}\in \mathbb{R}, w_{\eta}\in \varXi,\eta=1,\ldots,k-1$, such that
\begin{center}
$(\dfrac{\partial^{h+1}f}{\partial t^{h+1}} (t,\xi), \dfrac{\partial^{h+1}f}{\partial t^{h}\partial \xi}(t,\xi))(r_{\eta},w_{\eta})= r_{\eta}\,\dfrac{\partial^{h+1}f}{\partial t^{h+1}} (t,\xi) +  \dfrac{\partial^{h+1}f}{\partial t^{h}\partial \xi}(t,\xi)w_{\eta}=\delta_{h\eta}$\, ,
\end{center}
for $h=1,\ldots,k-1$.\medskip \\
Since $\dfrac{\partial^{h+1}f}{\partial t^{h+1}} (t,\xi)=0, \; h=1,\ldots,k-1$, the Proposition is  proved.$\bracevert$\\

\begin{proposition}\label{Pro264M}\textbf{(M).}\textit{ Let }$F:U\subseteq \mathbb{R}\times \varXi\rightarrow V \subseteq \mathbb{R}\times \varXi$\textit{ be a }$C^{d}$\textit{ LS-map }$(d\geq 2)$.\textit{ Then, for a given integer }$k$\textit{ such that }$1\leq k\leq d-1$,\textit{ a point }$(t,\xi)\in U$\textit{ is a }maximal $k$-transverse singularity\textit{ for }$F$ iff\textit{ the following condition is satisfied:}
\begin{align*}
(\text{M}_{k}^{\,\prime}) \qquad & \frac{\partial f}{\partial t}(t,\xi)=\ldots=\frac{\partial^{k+2}f}{\partial t^{k+2}}(t,\xi)=0\, ;\\
&\exists \, w_{\eta}\in \varXi , \eta=1,\ldots,k:\frac{\partial^{h+1}f}{\partial t^{h}\partial \xi}(t,\xi)w_{\eta}=\delta_{h\eta}, \; h,\eta=1.\ldots,k \, ; \\
&\forall \;  w\in \varXi: \frac{\partial ^{h+1}f}{\partial t^{h}\partial \xi}(t,\xi) w=0,\; h=1,\ldots,k \Rightarrow \frac{\partial^{k+2}f}{\partial t^{k+1}\partial \xi}(t,\xi) w=0,
\end{align*}
\textit{where the last condition is empty for } $k=d-1$.
\end{proposition}
\textbf{Proof. }From Definition \ref{D211},  $(t,\xi)\in U$ is a maximal $k$-transverse singularity \textit{iff} $(t,\xi) \in S_{1}(F)$ and conditions (M$_{k}$) hold for the canonical f-pair, i.e.
\begin{align*}
J_{h}(\varphi_{C},\psi_{C})(t,\xi) &=\frac{\partial^{h+1}f}{\partial t^{h+1}} (t,\xi)= 0, h=0\ldots,k;\\
I_{h}(\varphi_{C},\psi_{C})(t,\xi) &\text{ are l.i.}, \;  h=1,\ldots,k ; \\
I_{h}(\varphi_{C},\psi_{C})(t,\xi) &\text{ are not l.i.}, \; h=1,\ldots,k + 1.
\end{align*}
Note that we can add the requirement $ J_{k+1}(\varphi_{C},\psi_{C})(t,\xi)=\dfrac{\partial^{k+2}f}{\partial t^{k+2}}(t,\xi)=0$ to the above conditions. Otherwise, if $J_{k+1}(\varphi_{C},\psi_{C})(t,\xi)\neq 0$ then, by virtue of Proposition \ref{Pro245}, we should have that $I_{h}(\varphi_{C},\psi_{C})(t,\xi)$ are l.i., for $h=1,\ldots,k+1$. \\ 
By the Algebraic Lemma \ref{Lem221}, a) $\Leftrightarrow$ f), the conditions
\begin{center}
$I_{h}(\varphi_{C},\psi_{C})(t,\xi)\text{ are l.i., }h=1,\ldots,k\, ,$\\
$I_{h}(\varphi_{C},\psi_{C})(t,\xi)\text{ are not l.i.}, h=1,\ldots,k+1\,,$
\end{center}
are equivalent to:
\begin{align*}
&\exists \; v_{\eta}=(r_{\eta},w_{\eta}) \in \mathbb{R}\times \varXi,\; \eta=1,\ldots,k:I_{h}(\varphi_{C},\psi_{C})v_{\eta}=\delta_{h\eta}, \; h,\eta=1,\ldots,k \, ,\\
&\forall \, v=(r,w)\in \mathbb{R} \times \varXi: I_{h}(\varphi_{C},\psi_{C})(t,\xi)v=0,\, h=1,\ldots,k  \Rightarrow \\
&\qquad \qquad \qquad \qquad \qquad I_{k+1}(\varphi_{C},\psi_{C})(t,\xi)v=0.
\end{align*}
Since $I_{h}(\varphi_{C},\psi_{C})(t,\xi)=(\dfrac{\partial^{h+1}f}{\partial t^{h+1}} (t,\xi),\dfrac{\partial^{h+1}f}{\partial t^{h}\partial \xi}(t,\xi))= (0,\dfrac{\partial^{h+1}f}{\partial t^{h}\partial \xi}(t,\xi)),\;  h=1,\ldots,k+1$, by arguing as in the previous propositions we can conclude the proof.$\bracevert$ \\
\begin{proposition}\label{Pro265Inf}\textbf{($\infty $).}\textit{ Let } $F:U\subseteq\mathbb{R}\times \varXi \rightarrow V \subseteq \mathbb{R} \times \varXi $\textit{ be a }$C^{d}$\textit{ LS-map }$(d=\infty$\textit{ or }$d=\omega)$.\textit{ Then a point }$(t,\xi)\in U$\textit{ is an }$\infty$-transverse singularity\textit{ for }$F$ if\textit{ the following condition is satisfied:}
\begin{align*}
(\text{T}_{\infty}^{\,\prime}) \qquad \qquad & \frac{\partial^{h+1} f}{\partial t^{h+1}}(t,\xi)=0 \, , \, h\geq 0\, ;\\
&\exists w_{\eta}\in \varXi , \eta\geq 1 : \frac{\partial^{h+1} f}{\partial t^{h}\partial\xi}(t.\xi)w_{\eta}=\delta_{h\eta}, \; h,\eta\geq 1.
\end{align*}
\end{proposition}
\par \textbf{Proof. }From Definition \ref{D211}, $(t,\xi)\in U$ is an $\infty$-transverse singularity for $F$\textit{ iff} $(t,\xi) \in S_{1}(F)$ and the conditions  (T$_{\infty}$) hold for the canonical f-pair, that is
\begin{align*}
&J_{h}(\varphi_{C},\psi_{C})(t,\xi)=\frac{\partial^{h+1}f}{\partial t^{h+1}} (t,\xi)= 0 , \; h\geq 0\, ;\\
&\{I_{h}(\varphi_{C},\psi_{C})(t,\xi), \; h\geq 1\} \text{ are l.i. .}
\end{align*}
The Algebraic Lemma \ref{Lem221}, d) $\Rightarrow$ a), implies that the last condition is satisfied if there exist vectors  $v_{\eta}\in \mathbb{R} \times \varXi,\, \eta\geq 1$, such that $I_{h}(\varphi_{C},\psi_{C})(t,\xi)v_{\eta}=\delta_{h\eta}, \;h,\eta \geq 1$.\\
Since $I_{h}(\varphi_{C},\psi_{C})(t,\xi)=(0,\dfrac{\partial^{h+1}f}{\partial t^{h}\partial \xi}(t,\xi)), \; h\geq 1$, this proves the thesis. $\bracevert$ \\
\par
Actually, one could show that $(t,\xi)\in U$ is an $\infty$-transverse singularity for $F$ \textit{iff}
\begin{align*}
&\dfrac{\partial^{h+1}f}{\partial t^{h+1}}(t,\xi)=0,\; h\geq 0\,;\\
&\exists \, w_{\eta}\in \varXi,\;\eta\geq 1:  \dfrac{\partial^{h+1}f}{\partial t^{h}\partial \xi}(t,\xi)w_{\eta}=\delta_{h\eta}, 1\leq h\leq \eta .
\end{align*}
However, for the sake of simplicity, in stating the above Proposition we only gave a sufficient condition for $(t,\xi)\in U $ to be an $\infty$-transverse singularity. We also note that the vice versa of the proposition could be not true: the linear independence of the functionals  $(0,\dfrac{\partial^{h+1}f}{\partial t^{h}\partial \xi}(t,\xi))$, for  $h\geq 1$, i.e. the fact that $\dfrac{\partial^{h+1}f}{\partial t^{h}\partial \xi}(t,\xi)$ are l.i., $h\geq 1$, does not ensure the existence of vectors $w_{\eta}\in \varXi, \eta\geq 1$, such that $\dfrac{\partial^{h+1}f}{\partial t^{h}\partial \xi}(t,\xi)w_{\eta}=\delta_{h\eta}$, for $h,\eta \geq 1$. The following example better illustrates this behaviour.\\
\indent Consider the Hilbert space  $\textit{l}\,^{2}(\mathbb{N}):=\{\text{ real sequences }(t_{h})_{h\geq 1}=(t_{1},t_{2},t_{3},\ldots):\, \lVert(t_{h})_{h\geq 1}\rVert \,:=(\sum_{h=1}^{\infty}\!\mid\! t_{h}\!\mid^{2})^{1/2}<+\infty\}$ and the canonical orthonormal basis $\{e_{h},h\geq 1\}$. Given the vector $e_{0} :=(1,1/2,1/3,1/4,\ldots)\in \textit{l}\,^{2}(\mathbb{N})$, one can easily see that the family of vectors $\{e_{h},h\geq 0\}$ is l.i. in $\textit{l}\,^{2}(\mathbb{N})$. Define the functionals $I_{h}\in \textit{l}\,^{2}(\mathbb{N})^{\ast}, h\geq 0$, as $I_{h}(\cdotp):=\,<\!e_{h},\cdotp>$, where $<\cdotp,\cdotp>$ is the usual $\textit{l}\,^{2}(\mathbb{N})$ scalar product. Then the family $\{I_{h}, h\geq 0\}$ is  l.i. in $\textit{l}\,^{2}(\mathbb{N})^{\ast}$, but it cannot exist $v\in \textit{l}\,^{2}(\mathbb{N})$ such that $I_{0}(v)=1$ and $I_{h}(v)=0$ for $h\geq 1$ because the last condition is equivalent to $v=0$.\\
\par
We are now able to add, for the class of LS-maps, examples of all the singularities considered in Definition \ref{D211}. Similar examples, but related to nonlinear differential problems, are studied in \cite{B-D 2}, \cite{B-D 3} and \cite{B-D 4}.\\
\begin{example}
\label{Exe266} For $k\geq 1$, let $Z$ be a $B$-space and let us consider the map
\begin{align*}
F:\mathbb{R}^{k+1}\times Z &\rightarrow \mathbb{R}^{k+1}\times Z\\
(t,t_{1},\ldots,t_{k},z) &\mapsto (\sum_{h=1}^{k}t_{h}t^{h},t_{1},\ldots,t_{k},z).
\end{align*}
If we set $\varXi :=\mathbb{R}^{k}\times Z$ and $\xi:=(t_{1},\ldots,t_{k},z)\in \varXi$ we are in the situation described in Subsection 2.6.1, where  $f(t,\xi)=f(t,t_{1},\ldots,t_{k},z):=t_{k}t^{k}+t_{k-1}t^{k-1}+\ldots+ t_{2}t^{2}+t_{1}t$.\\
We want to show that the points  $(0,\ldots,0,z), z\in Z$, are $k$-transverse singularities for $F$. It is clear that $\dfrac{\partial^{h+1}f}{\partial t^{h+1}}(0,\ldots,0,z)=0$ for every $h\geq 0, z\in Z$. It is now convenient to define vectors $e_{\eta}\in \mathbb{R}^{k+1 }\times Z, \eta =1,\ldots,k$, by $e_{\eta}:=(0,\ldots,0,1,0,\ldots,0,0_{Z})$, with $1$ at the $(\eta+1)-$th place. Hence we easily get that $\dfrac{\partial^{h+1}f}{\partial t^{h} \partial \xi}(t,\xi)e_{\eta}=\dfrac{\partial^{h+1}f}{\partial t^{h} \partial t_{\eta}}(t,\xi)$, for $h\geq 1$, and in particular $\dfrac{\partial^{h+1}f}{\partial t^{h} \partial \xi}(0,\ldots,0,z)e_{\eta}=\dfrac{\partial^{h+1}f}{\partial t^{h} \partial t_{\eta}}(0,\ldots,0,z)=\eta!\delta_{h\eta} ,\text{ for }h,\eta=1,\ldots,k,z\in Z$.
Since $\dfrac{\partial^{k+1}f}{\partial t^{k+1}}(t,\xi)=0$ we can take $r_{\eta}$ in an arbitrary way, and by choosing $w_{\eta}:=e_{\eta}/\eta!, \eta=1,\ldots,k$, we get that conditions $(\text{T}_{k}^{\,\prime})$ in Proposition 2.6.1(T) are satisfied. Thus the points $(0,\ldots,0,z), z\in Z$, are  $k$-transverse singularities  for $F$.
\end{example}
\begin{example} \label{Exe267} Given $k,n\in \mathbb{N}$ let us  study the map \medskip
\begin{align*}
F:\mathbb{R}^{k+1}\times Z&\rightarrow \mathbb{R}^{k+1}\times Z\\
(t,t_{1},\ldots,t_{k},z)&\mapsto((1-\delta_{0n})t^{n}+\sum_{h=1}^{k} t_{h}t^{h},t_{1},\ldots,t_{k},z).
\end{align*}
\indent We note that, for any given $k,n\in \mathbb{N}$ and $z\in Z$, we have $F(0,\ldots,0,z)=(0,\ldots,0,z)$. Now we shall classify the points $(0,\ldots,0,z)$ with respect to $k,n\in \mathbb{N}$. We can prove that:\\

\noindent Case $k=0$:\\
- for $n=0$ or for $ n\geq 3 $ the points $(0,z),z\in Z $, are not 1-transverse 
singularities\\
\indent for $F$;\\
- for $n=1$ the points $(0,z)$ are \textit{regular} points for $F$, i.e. $F'(0,z)$ is an 
isomorphism;\\
- for $n=2$ the points $(0,z)$ are 1-singularities for $F$.\\

\noindent Case $k\geq 1$:\\
- for $n=0$ or for $n\geq k+3$ the points $(0,\ldots,0,z),z\in Z$, are maximal $k$-transverse\\
\indent singularities for $F$;\\
- for $n=1$ the points $(0,\ldots,0,z)$ are \textit{regular} points for $F$, i.e. $F'(0,\ldots,0,z)$ is an \\
\indent isomorphism;\\
- for $n=2,\ldots,k+2$ the points $(0,\ldots,0,z)$ are $(n-1)$-singularities for $F$.
\end{example}
\noindent \textit{Proof of case }$k=0$.\\
We can write  $F:\mathbb{R} \times Z\rightarrow \mathbb{R}\times Z$ as the map $F(t,z) =(f(t,z),z)$, with $f(t,z):=(1-\delta_{0n})t^{n}$. We also recall that 
\begin{center}
$S_{1}(F)=\{(t,z)\in U:\dfrac{\partial f}{\partial t}(t,z)=0\}$.
\end{center}
- For $n=1$ one has $\dfrac{\partial f}{\partial t}(t,z)=1$ and so the points $(0,z)$ are not simple singularities  for F. Indeed,
\begin{equation*}
F'(0,z)= \left[ \begin{array}{cc} 
\dfrac{\partial f}{\partial t}(0,z) & \dfrac{\partial f}{\partial \xi}(0,z) \\
\\
0 & \textbf{1}_{\textbf{Z}} \end{array} \right]
\end{equation*}
is an isomorphism (see also Subsection \ref{sss131}), and thus the points $(0,z)$ are regular for $F$.\\
When $n\neq 1$ it follows that $\dfrac{\partial f}{\partial t}(t,z)=0$ and so the points $(0,z)$ are singular for $F$. Let us study these singularities.\\
- For $n=0$ or $n\geq 3$, it is easily checked that $\dfrac{\partial^{2}f}{\partial t^{2}}(0,z)=\dfrac{\partial^{2}f}{\partial t \partial z}(0,z)=0$ and so, by applying condition $(\text{T}_{1}^{\,\prime})$ in Proposition \ref{Pro262T}(T), we deduce that the points $(0,z),z\in Z$, are not 1-transverse.\\
- Finally, for $n=2$ one has that $\dfrac{\partial^{2} f}{\partial t^{2}}(0,z)\neq 0$ and, by condition $(\text{S}_{1}^{\,\prime})$ in Proposition \ref{Pro263S}(S), we get that the points $(0,z)$ are 1-singularities for $F$.\\

\noindent \textit{Proof of case }$k\geq 1$.\\
As seen in Example \ref{Exe266}, we consider $F:\mathbb{R}\times \varXi \rightarrow \mathbb{R}\times \varXi$ with $\varXi:=\mathbb{R}^{k}\times Z$ and $F(t,\xi)=(f(t,\xi),\xi)$. Here $\xi:=(t_{1},\ldots,t_{k},z),f(t,\xi):=(1-\delta_{0n})t^{n}+t_{k}t^{k}+t_{k-1}t^{k-1}+\ldots+t_{2}t^{2}+t_{1}t$. It is useful to define $g(t,\xi):=t_{k}t^{k}+t_{k-1}t^{k-1}+\ldots+t_{2}t^{2}+t_{1}t$, thus $f(t,\xi)=(1-\delta_{0n})t^{n}+g(t,\xi)$. Moreover, as recalled in the case $k=0$, one has that the points $(0,\ldots,0,z), z\in Z$, are simple singularities \textit{iff} $\dfrac{\partial f}{\partial t}(0,...,0,z)=0$.\\
- For $n=1$ one has $\dfrac{\partial f}{\partial t}(0,\ldots,0,z)=1$ and the points $(0,\ldots,0,z),z\in Z$, are regular for $F$.\\
When $n\neq 1$ one has that $\dfrac{\partial f}{\partial t}(0,\ldots,0,z)=0$ and we now proceed to study the nature of the singular points $(0,\ldots,0,z), z\in Z$.\\
- For $n=0$ or $n\geq k+3$ one has $\dfrac{\partial^{h+1}f}{\partial t^{h+1}}(0,\ldots,0,z)=0,h=0,\ldots,k+1$. Moreover, if we argue as in Example \ref{Exe266}, by choosing  $w_{\eta}:=e_{\eta}/\eta!,\eta= 1,\ldots,k$,  we obtain that 
\begin{center}
$\dfrac{\partial^{h+1}f}{\partial t^{h}\partial \xi}(0,\ldots,0,z)w_{\eta}=1/\eta!\dfrac{\partial^{h+1}f}{\partial t^{h}\partial t_{\eta}}(0,\ldots,0,z)=1/\eta!\dfrac{\partial^{h+1}g}{\partial t^{h}\partial t_{\eta}}(0,\ldots,0,z)=\delta_{h\eta},$
\end{center}
$h,\eta=1,\ldots, k, z\in Z$. Finally, $\dfrac{\partial^{k+2}f}{\partial t^{k+1}\partial t_{\eta}}\equiv 0, \eta=1,\ldots,k$, and $\dfrac{\partial^{k+2}f}{\partial t^{k+1}\partial z}\equiv 0$. Hence $\dfrac{\partial^{k+2}f}{\partial t^{k+1}\partial \xi}(t,\xi)w\equiv 0, (t,\xi)\in \mathbb{R}\times \varXi, w \in \varXi$. Since conditions $(\text{M}_{k}^{\,\prime})$ in Proposition \ref{Pro264M}(M) are satisfied, this yields that the points $(0,\ldots,0,z),z\in Z$, are maximal $k$-transverse singularities.\\
- For $n=2,\ldots,k+2$, we have that $\dfrac{\partial f}{\partial t}(0,\ldots,0,z)=\ldots=\dfrac{\partial^{n-1}f}{\partial t^{n-1}}(0,\ldots,0,z)=0$ and $\dfrac{\partial^{n}f}{\partial t^{n}}(0,\ldots,0,z)=n!\neq 0$. Furthermore, as seen above, for $w_{\eta}= e_{\eta}/\eta!, \eta=1,\ldots,n-2\leq k$, it follows that
\begin{center}
$\dfrac{\partial ^{h+1}f}{\partial t^{h}\partial \xi}(0,\ldots,0,z)w_{\eta}=\delta_{hn} ,$
\end{center}
$h,\eta=1,\ldots,n-2, z\in Z$. Thus, from condition $(\text{S}_{n-1}^{\,\prime})$ in Proposition \ref{Pro263S}(S), the points $ (0,\ldots,0,z)$ are $(n-1)$-singularities.\\

\begin{example}
\label{Exe268} In the above example we proved that, given  $k\geq 0$ and a $B$-space $Z$, the points $(0,\ldots,0,z), z\in Z$, are $(k+1)$-singularities for the map $F$ defined as
\begin{align*}
F:\mathbb{R}^{k+1}\times Z&\rightarrow \mathbb{R}^{k+1}\times Z\\
(t,t_{1},\ldots,t_{k},z)&\mapsto(t^{k+2}+\sum_{h=1}^{k} t_{h}t^{h},t_{1},\ldots,t_{k},z).
\end{align*}
When the above map is rewritten in the form
\begin{align*}
w_{k,Z}:\mathbb{R}^{k}\times Z&\rightarrow \mathbb{R}^{k}\times Z\\
(t,t_{1},\ldots,t_{k-1},z)&\mapsto(t^{k+1}+\sum_{h=1}^{k-1} t_{h}t^{h},t_{1},\ldots,t_{k-1},z),
\end{align*}
for an integer $k\geq 1$, it is usually called \textit{generalized Whitney map} $w_{k,Z}$ or, more simply, \textit{Whitney map }(when $Z=\{0\}$ we just write $w_{k}$ instead of $w_{k,\{0\}})$. Thus, by definition, 
\begin{center}
$w_{k,Z}(t,t_{1},\ldots,t_{k-1},z)=(t^{k+1}+t_{k-1}t^{k-1}+\ldots+t_{2}t^{2}+t_{1}t,t_{1},\ldots,t_{k-1},z).$
\end{center}
Hence, for $k\geq 1$, the points $(0,\ldots,0,z), z\in Z$, are $k$-singularities for $w_{k,Z}.$
\end{example}
\begin{example}\label{Exe269} Now we shall give an example of an $\infty$-transverse singularity. Consider $\mathbb{R}_{1}:=\{t\in \mathbb{R}:\lvert t \rvert < 1\}$, the Hilbert space $\textit{l\,}^{2}(\mathbb{N})=\{\text{real sequences } (t_{h})_{h\geq 1}=(t_{1},t_{2},t_{3},\ldots):\lVert(t_{h})_{h\geq 1}\rVert:=(\sum_{h=1}^{\infty}\lvert t_{h}\lvert^{2})^{1/2}<+\infty \}$ and let $Z$ be a $B$-space. Define the map
\begin{align*}
F:\mathbb{R}_{1}\times \textit{l}\,^{2}(\mathbb{N}) \times Z\subseteq \mathbb{R}\times \textit{l}\,^{2}(\mathbb{N}) \times Z &\rightarrow \mathbb{R}\times \textit{l}\,^{2}(\mathbb{N}) \times Z\\
(t,t_{1},t_{2},t_{3},\ldots,z) &\mapsto (\sum_{h=1}^{\infty}t_{h}t^{h},t_{1},t_{2},t_{3},\ldots,z).
\end{align*}
As usual, we are considering $F$ of the form  $F:\mathbb{R}_{1}\times \varXi \subseteq \mathbb{R}\times \varXi \rightarrow \mathbb{R}\times \varXi$, where $\varXi:=\textit{l}\,^{2}(\mathbb{N}) \times Z$ and $F(t,\xi)=(f(t,\xi),\xi)$, for  $\xi:=(t_{1},t_{2},t_{3},...,z)$ and $f(t,\xi):=\sum_{h=1}^{\infty}t_{h}t^{h}$.\\
The real function $f$ is well-defined. In fact, from the Cauchy-Schwarz inequality, 
\begin{center}
$\lvert\sum_{h=1}^{\infty}t_{h}t^{h}\lvert\leq(\sum_{h=1}^{\infty}(t_{h})^{2})^{1/2}(\sum_{h=1}^{\infty}(t^{h})^{2})^{1/2}=\lVert(t_{h})_{h\geq 1}\lVert(t^{2}/1-t^{2})^{1/2}< +\infty.$
\end{center}
One can also show that $f$ is a $C^{\omega}$ function and can be termwise differentiated. This implies that $F$ is a $C^{\omega}$ map. Hence, for $k\geq 1, \dfrac{\partial^{k}f}{\partial t^{k}}(t,\xi)=\sum_{h=k}^{\infty}h(h-1)\ldots(h-k+1)t_{h}t^{h-k}$ and so $\dfrac{\partial^{k}f}{\partial t^{k}}(0,0,0,0,\ldots,z)=0$. In a similar way, if for $\eta \geq 1$ we define $e_{\eta}:=(0,0,0,0,\ldots,0,1,0,\ldots,z)\in \mathbb{R}\times\textit{l}\,^{2}(\mathbb{N}) \times Z$, with $1$ at the $(\eta+1)$-th place, we obtain that
$\dfrac{\partial^{k+1}f}{\partial t^{k}\partial \xi}(0,0,0,0,\ldots,z)(e_{\eta}/\eta!)=1/\eta!\cdotp \dfrac{\partial^{k+1}f}{\partial t^{k}\partial t_{\eta}}(0,0,0,0,\ldots,z)=1/\eta!\cdotp k!\cdotp \delta_{k\eta}=\delta_{k\eta}$.\\
By choosing $w_{\eta}=e_{\eta}/\eta!$ we get that condition $(\text{T}_{\infty}^{\,\prime})$ in Proposition \ref{Pro265Inf}($\infty$) is satisfied. Thus the points $(0,0,0,0,\ldots,z)$ are $\infty$-transverse singularities for $F$.
\end{example}
\begin{remark}\label{Rem270} We conclude the section by showing that maximal k-transverse singularities are not stable, which means they can be eliminated by small perturbations of the map. We will show this phenomenon by using one of the LS-maps that generate maximal 1-transverse singularities, seen in Example \ref{Exe267}.
Let us consider the map
\begin{align*}
F_{\varepsilon}:\mathbb{R}^{2}&\rightarrow \mathbb{R}^{2}\\
(t,\xi)&\mapsto(t\xi-\frac{\varepsilon}{2}t^{2},\xi).
\end{align*}

The singular set is given by $S_{1}(F_{\varepsilon}) = \{(t,\xi)\in U: \dfrac{\partial f_{\varepsilon}}{\partial t}(t,\xi)=0\}$ where $f_{\varepsilon}(t,\xi) = t\xi - \dfrac{\varepsilon}{2}t^{2}$. Hence $\dfrac{\partial f_{\varepsilon}}{\partial t}(t,\xi) = \xi - \varepsilon t$ and so $S_{1}(F_{\varepsilon})$ is the straight-line with equation $\xi = \varepsilon t$. Since  $\dfrac{\partial^{2} f_{\varepsilon}}{\partial \xi \partial t}(t,\xi)= 1$ the singular points are 1-transverse; we also remark that they cannot be 2-transverse because the third derivatives are zero (cf. Proposition \ref{Pro262T}(T)). Finally, $\dfrac{\partial^{2} f_{\varepsilon}}{\partial t^{2}}(t,\xi)= \varepsilon$. For $\varepsilon = 0$ the singular set coincides with the t-axis and is made up of maximal 1-transverse singularities, given that $\dfrac{\partial^{2} f_{0}}{\partial t^{2}}(t,\xi)= 0$. When $\varepsilon \neq 0$ the singular points become 1-singularities, i.e. fold points, because $\dfrac{\partial^{2} f_{\varepsilon}}{\partial t^{2}}(t,\xi)\neq 0$ (cf. Proposition \ref{Pro263S}(S)).
\end{remark}

\end{document}